\let\chooseClass1   
\let\chooseClass2   

\ifx\chooseClass1
\IfFileExists{extarticle.cls}{
\documentclass[14pt]{extarticle}
\emergencystretch 7 pt
\setlength{\oddsidemargin}{-17 mm} 
\setlength{\textwidth}{192 mm}     
\setlength{\textheight}{246 mm}    
\setlength{\topmargin}{-31 mm}     
}{
\documentclass[12pt]{article}
\textheight = 8.5in
\textwidth 6.3in
\setlength{\oddsidemargin}{0mm}
\setlength{\topmargin}{-20 mm}
}{}
\fi

\ifx\chooseClass2
\documentclass[12pt]{article}
\fi

\usepackage{amsmath,
            amsfonts,
            amssymb,
            amsthm,
            color
}

\IfFileExists{euscript.sty}{\usepackage[mathcal]{euscript}}{}
\IfFileExists{mathrsfs.sty}{\usepackage{mathrsfs}}{\let\mathscr\mathfrak}

\InputIfFileExists{diagrams.tex}{}{%
\IfFileExists{diagrams.sty}{\usepackage{diagrams}}
{\message{Without diagrams you can not process this file!}}}
\diagramstyle[height=2em,balance,righteqno,PostScript=dvips,nohug]

\usepackage{ifpdf}
\ifpdf
    \IfFileExists{hyperref.sty}{\usepackage[pdftex]{hyperref}}{}
\else
    \IfFileExists{hyperref.sty}{\usepackage[hypertex]{hyperref}}{}
\fi

\IfFileExists{dsfont.sty}%
{\usepackage[sans]{dsfont}\newcommand\1{{\mathds 1}}}%
{\newcommand\1{{1\mkern-5mu {\mathrm I}}}}

\makeatletter
\def\@seccntformat#1{\csname the#1\endcsname.\quad}
\renewcommand\section{\@startsection {section}{1}{\z@}%
                                   {-3.5ex \@plus -1ex \@minus -.2ex}%
                                   {2.3ex \@plus.2ex}%
                                   {\normalfont\large\bfseries}}
\renewcommand\subsection{\@startsection{subsection}{2}{\z@}%
                        {3.25ex plus 1ex minus .2ex}{-.5em}%
                        {\normalfont\normalsize\bfseries}}
\@addtoreset{equation}{section}
\numberwithin{equation}{subsection}
\makeatother

\newtheoremstyle{boldhead}
{\topsep}
{\topsep}
{\slshape}
{}
{\bfseries}
{.}
{ }
{\thmname{#1}\thmnumber{ #2}\thmnote{ (#3)}}

\newtheoremstyle{boldremark}
{\topsep}
{\topsep}
{\upshape}
{}
{\bfseries}
{.}
{ }
{\thmname{#1}\thmnumber{ #2}\thmnote{ (#3)}}

\swapnumbers

\theoremstyle{boldhead}

\newtheorem{theorem}[subsection]{Theorem}
\newtheorem{corollary}[subsection]{Corollary}
\newtheorem*{corollaryNoNumber}{Corollary}
\newtheorem{lemma}[subsection]{Lemma}
\newtheorem{proposition}[subsection]{Proposition}

\theoremstyle{boldremark}

\newtheorem{definition}[subsection]{Definition}
\newtheorem{example}[subsection]{Example}
\newtheorem{examples}[subsection]{Examples}
\newtheorem{remark}[subsection]{Remark}
\newtheorem*{acknowledgement}{Acknowledgements}

\newenvironment{events}{\begin{list}{}{\addtolength{\itemsep}{-0.6em}}}{\end{list}}
\newcommand{\sbull}{\makebox[0mm][r]{$\scriptstyle\bullet$}}
\newcommand{\bull}{{\scriptscriptstyle\bullet}}

\def\rhaha{\raise.24ex\hbox{$\rightharpoonup$}\kern-1em\lower.24ex\hbox{$\rightharpoondown$}}%
\def\lhaha{\raise.24ex\hbox{$\leftharpoonup$}\kern-1em\lower.24ex\hbox{$\leftharpoondown$}}%
\def\dhaha{\downharpoonleft\kern-.22em\downharpoonright\kern.02em}%
\def\uhaha{\upharpoonleft\kern-.22em\upharpoonright\kern.02em}
\newarrowhead{twoharpoons}\rhaha\lhaha\dhaha\uhaha

\newarrow{MapsTo}{mapsto}---{->}
\newarrow{Mono}{boldhook}---{->}
\newarrow{TTo}----{->}
\newarrow{Twoar}===={=>}

\newcommand\NN{{\mathbb N}}
\newcommand\ZZ{{\mathbb Z}}

\newcommand{\ca}{{\mathcal A}}
\newcommand{\cb}{{\mathcal B}}
\newcommand{\cc}{{\mathcal C}}
\newcommand{\cd}{{\mathcal D}}
\newcommand{\ck}{{\mathcal K}}
\newcommand{\cv}{{\mathcal V}}

\newcommand{\fA}{{\mathfrak A}}
\newcommand{\fu}{{\mathscr U}}

\newcommand{\Cat}{{\mathcal C}at}
\newcommand{\CatCat}{{\mathcal C}at\text-{\mathcal C}at}
\newcommand{\Com}{\mathsf C}
\newcommand{\KCat}{{\mathcal K}\text-{\mathcal C}at}
\newcommand{\KCatb}{{\mathcal K}\text-{\mathcal C}at}
\newcommand{\KCatCat}{{\mathcal K}\text-{\mathcal C}at\text-{\mathcal C}at}
\newcommand{\kf}{{\mathsf k}}
\newcommand{\Kht}{\mathsf K}

\newcommand{\hcirc}{\circ_h}
\newcommand{\leftunit}{{\mathbf l}}
\newcommand{\rightunit}{{\mathbf r}}
\newcommand{\uni}{{\mathbf i}}

\newcommand{\sS}[2]{\vphantom{#2}#1 #2}

\newcommand{\n}[1]{\nobreakdash-\hspace{0pt}}
\newcommand{\ainf}[1]{$A_\infty$\nobreakdash-\hspace{0pt}}

\let\kk\Bbbk
\let\emptyset\varnothing
\let\eps\varepsilon
\let\epsilon\varepsilon
\let\ge\geqslant
\let\le\leqslant

\let\tens\otimes

\DeclareMathOperator\Coder{Coder}
\DeclareMathOperator\Cone{Cone}

\DeclareMathOperator\gr{gr}
\DeclareMathOperator\grad{-grad}
\DeclareMathOperator\Hom{Hom}
\DeclareMathOperator\HOM{\sf Hom}
\DeclareMathOperator\id{id}

\DeclareMathOperator\im{Im}
\DeclareMathOperator\modul{-mod}
\DeclareMathOperator\Mor{Mor}
\DeclareMathOperator\Ob{Ob}
\DeclareMathOperator\pr{pr}

\newcommand{\appref}[1]{Appendix~\ref{#1}}
\newcommand{\corref}[1]{Corollary~\ref{#1}}
\newcommand{\defref}[1]{Definition~\ref{#1}}

\newcommand{\exasref}[1]{Examples~\ref{#1}}

\newcommand{\lemref}[1]{Lemma~\ref{#1}}
\newcommand{\propref}[1]{Proposition~\ref{#1}}
\newcommand{\remref}[1]{Remark~\ref{#1}}
\newcommand{\secref}[1]{Section~\ref{#1}}
\newcommand{\thmref}[1]{Theorem~\ref{#1}}

\begin{document}
\bibliographystyle{amsalpha}
\title{Category of $A_\infty$-categories}
\author{Volodymyr Lyubashenko}
\date{lub@imath.kiev.ua
\thanks{Institute of Mathematics,
National Academy of Sciences of Ukraine,
3 Tereshchenkivska st.,
Kyiv-4, 01601 MSP,
Ukraine}
\thanks{The research was supported in part by grant 01.07/132 of
State Fund for Fundamental Research of Ukraine}
}
\maketitle

\begin{abstract}
We define natural \ainf-transformations and construct \ainf-category of
\ainf-functors. The notion of non-strict units in an \ainf-category is
introduced. The 2\n-category of (unital) \ainf-categories, (unital)
functors and transformations is described.
\end{abstract}

\footnotetext{
\ainf-functors, \ainf-transformations, unit \ainf-transformation,
2\n-category.}


\allowdisplaybreaks[1]

The study of higher homotopy associativity conditions for topological
spaces began with Stasheff's article \cite[I]{Stasheff:HomAssoc}. In a
sequel to this paper \cite[II]{Stasheff:HomAssoc} Stasheff defines also
\ainf-algebras and their homotopy-bar constructions. These algebras and
their applications to topology were actively studied, for instance, by
Smirnov~\cite{Smirnov80} and Kadeishvili
\cite{Kadeishvili80,Kadeishvili82}. We adopt some notations of Getzler
and Jones~\cite{GetzlerJones:A-infty}, which reduce the number of signs
in formulas. The notion of an \ainf-category is a natural
generalization of \ainf-algebras. It arose in connection with Floer
homology in Fukaya's work \cite{Fukaya:A-infty,Fukaya:FloerMirror-II}
and was related by Kontsevich to mirror symmetry
\cite{Kontsevich:alg-geom/9411018}. See Keller~\cite{math.RA/9910179}
for a survey on \ainf-algebras and categories.

In the present article we show that given two \ainf-categories $\ca$
and $\cb$, one can construct a third \ainf-category $A_\infty(\ca,\cb)$
whose objects are \ainf-functors $f:\ca\to\cb$, and morphisms are
natural \ainf-transformations between such functors. This result was
also obtained by Fukaya \cite{Fukaya:FloerMirror-II} and by Kontsevich
and Soibelman \cite{KonSoi-book}, independently and, apparently,
earlier. We describe compositions between such categories of
\ainf-functors, which would allow us to construct a 2\n-category of
unital \ainf-categories. The latter notion is our generalization of
strictly unital \ainf-categories (cf. Keller~\cite{math.RA/9910179}).
We also discuss unit elements in unital \ainf-categories, unital
natural \ainf-transformations, and unital \ainf-functors.

\subsection*{Plan of the article with comments and explanations.}
The first section describes some notation, sign conventions,
composition convention, etc. used in the article. The ground
commutative ring $\kk$ is not assumed to be a field. This is suggested
by the development of homological algebra in \cite{Drinf:DGquot}.
Working over a ring $\kk$ instead of a field has strong consequences.
For instance, one may not hope for Kadeishvili's theorem on minimal
models \cite{Kadeishvili82} to hold for all \ainf-algebras over $\kk$.

In the second section we recall or give definitions of the main
objects. A $\kk$\n-quiver is such a graph that the set of arrows
(morphisms) between two vertices (objects) is a $\kk$\n-module
(\defref{def-Quiver}). We view quivers as categories without
multiplication and units. Cocategories are $\kk$\n-quivers and
$\kk$\n-coalgebras with a matrix type decomposition into
$\kk$\n-submodules, indexed by pairs of objects
(\defref{def-Cocategory}). \ainf-categories are defined as a special
kind of differential graded cocategories -- the ones of the form of the
tensor cocategory $T\ca$ of a $\kk$\n-quiver $\ca$
(\defref{def-Ainf-category}). \ainf-functors are homomorphisms of
cocategories that commute with the differential
(\defref{def-ainf-functor}). \ainf-transformations between
\ainf-functors are defined as coderivations (\defref{def-ainf-trans}).
They seem to make \ainf-category theory closer to ordinary category
theory. Notice, however, that \ainf-transformations are analogs of
transformations between ordinary functors, which do not satisfy the
naturality condition. Natural \ainf-transformations are introduced in
\defref{def-globular-set}. \ainf-functors and \ainf-transformations are
determined by their components.

In the third section we study tensor products of cocategories and
homomorphisms between them (\secref{sec-cocat-hom}). We concentrate on
homomorphisms from the tensor product of tensor cocategories to another
tensor cocategory. Given $\kk$\n-quivers $\ca$ and $\cb$, we consider
another $\kk$\n-quiver $\Coder(\ca,\cb)$, whose objects are the
cocategory homomorphisms $T\ca\to T\cb$ and morphisms are coderivations
(\secref{sec-cocat-hom}). We construct a cocategory homomorphism
$\alpha:T\ca\tens T\Coder(\ca,\cb)\to T\cb$ (\corref{cor-action-alpha}
to \propref{pro-theta-D-D-theta-theta}), based on a map
$\theta:T\Coder(\ca,\cb)\to\Hom_\kk(T\ca,T\cb)$
\eqref{eq-theta-kl-sum}. The homomorphism $\alpha$ is universal
(\propref{prop-phi-C1-Cq-psi}), in other words, $T\Coder(\ca,\cb)$ is
the inner hom-object $\HOM(T\ca,T\cb)$ in the monoidal category
generated by tensor cocategories.

This universality is exploited in the fourth section in order to show
that the category of tensor cocategories is enriched in the monoidal
category of graded cocategories. That is: there exists an associative
unital multiplication
\(M:T\Coder(\ca,\cb)\tens T\Coder(\cb,\cc)\to T\Coder(\ca,\cc)\),
which is a cocategory homomorphism (\propref{pro-M-assoc-eta-unit}).
Its explicit description uses the map $\theta$.

The fifth section extends the results of the third section to
differential graded tensor cocategories, that is, to \ainf-categories.
With two \ainf-categories $\ca$, $\cb$ is associated a third
\ainf-category $A_\infty(\ca,\cb)$ (\propref{pro-B-F-KS-LH}). Its
objects are \ainf-functors $\ca\to\cb$, and its morphisms are
coderivations. To reduce the number of signs in the theory we prefer to
work with grading of a graded $\kk$\n-quiver or \ainf-category $\ca$
shifted by 1: $s\ca=\ca[1]$. In this notation the quiver
$A_\infty(\ca,\cb)$ is a full subquiver of the quiver
$s^{-1}\Coder(s\ca,s\cb)$. The proof of \propref{pro-B-F-KS-LH}
consists of constructing a differential $B$ in the tensor cocategory of
$A_\infty(\ca,\cb)$. The explicit formula \eqref{eq-Bn-components} for
$B$ uses the map $\theta$. A cocategory homomorphism from a tensor
product of differential tensor cocategories to a single such cocategory
is called an \ainf-functor in the generalized sense
(\secref{sec-Differentials}). Restricting the cocategory homomorphism
of \corref{cor-action-alpha} we get a homomorphism of differential
graded cocategories $Ts\ca\tens TsA_\infty(\ca,\cb)\to Ts\cb=T(s\cb)$
(\corref{cor-ainf-Ainfty-AB}). Its universality
(\propref{prop-phi-C1-Cq-psi-ainf}) may be interpreted as
$TsA_\infty(\ca,\cb)$ being the inner hom-object $\HOM(Ts\ca,Ts\cb)$ in
the monoidal category generated by differential graded tensor
cocategories.

This universality is used in the sixth section to show that the
category of \ainf-categories is enriched in the monoidal category of
differential graded cocategories. Namely, the multiplication $M$ of
\propref{pro-M-assoc-eta-unit} restricted to
\(M:TsA_\infty(\ca,\cb)\tens TsA_\infty(\cb,\cc)\to
TsA_\infty(\ca,\cc)\)
is an \ainf-functor, that is, it commutes with the differential
(equation~\eqref{eq-1B-B1-M-MB}). By universality
(\propref{prop-phi-C1-Cq-psi-ainf}) $M$ corresponds to a unique
\ainf-functor
\[ A_\infty(\ca,\_): A_\infty(\cb,\cc) \to
A_\infty(A_\infty(\ca,\cb),A_\infty(\ca,\cc)).
\]
We prove that it is strict and describe it in
\propref{pro-Ainf-A--strict}. Natural \ainf-transformations are defined
as cycles in the differential graded quiver of all
\ainf-transformations (\defref{def-globular-set}).

Identifying cohomologous natural \ainf-transformations (that is,
considering cohomology of the quiver of \ainf-transformations) in the
seventh section, we get a non-2-unital 2\n-category $A_\infty$, whose
objects are \ainf-categories, 1\n-morphisms are \ainf-functors, and
2\n-morphisms are equivalence classes of natural \ainf-transformations.
Here non-2-unital means that unit 2\n-morphisms are missing in the
2\n-category $A_\infty$. However, unit 1\n-morphisms are present -- the
identity \ainf-functors. Before constructing $A_\infty$ we construct a
non-2-unital 2\n-category $\ck A_\infty$ enriched in $\ck$ -- the
homotopy category of differential graded complexes of $\kk$\n-modules
(\propref{pro-1uni-n2uni-KAinf}). Morphisms of $\ck$ are chain maps
modulo homotopy. The notion of 2\n-category enriched in a symmetric
monoidal category is discussed in \appref{ap-sec-enrich}. The idea of
the construction is that the binary operation becomes strictly
associative if homotopic chain maps are identified. Similarly with
other identities in a 2\n-category. The non-2-unital 2\n-category
$A_\infty$ is obtained from $\ck A_\infty$ by taking the 0\n-th
cohomology.

\ainf-categories are analogs of non-unital categories -- categories
without unit morphisms. We define a unital \ainf-category $\cc$ so that
its cohomology $H^\bull(\cc)$ is a unital category, and for any
representative $1_X\in\cc^0(X,X)$ of the unit class
$[1_X]\in H^0(\cc(X,X))$ the binary compositions with $1_X$ are
homotopic to identity as chain maps $\cc(X,Y)\to\cc(X,Y)$ or
$\cc(Y,X)\to\cc(Y,X)$ (\defref{def-unital-cat},
\lemref{lem-1ib1-i1b-1}). We prove that for a unital \ainf-category
there exists a natural \ainf-transformation
$\uni^\cc:\id_\cc\to\id_\cc:\cc\to\cc$ of the identity functor, whose
square is equivalent to $\uni^\cc$ (\propref{pro-unital-unit}). It is
called a unit transformation of $\cc$ (\defref{def-unit}) and, indeed,
it is a unit 2\n-morphism in the 2\n-category $A_\infty$. Moreover,
$f\uni^\cc:f\to f$ is a unit 2\n-morphism of an \ainf-functor
$f:\ca\to\cc$ (\corref{cor-rgiB-r-firB-r}). The unit transformation
$\uni^\cc$ is determined uniquely up to an equivalence
(\corref{cor-unit-uniquely}). If $\cc$ is unital, then
$A_\infty(\ca,\cc)$ is unital as well (\propref{pro-C-unit-AAC-unit}).

The full 2-subcategory $\ck\sS{^u}A_\infty$ (resp. $\sS{^u}A_\infty$)
of $\ck A_\infty$ (resp. $A_\infty$), whose objects are unital
\ainf-categories, and 1\n-morphisms are all \ainf-functors, is
{\color{red} almost} 2-unital by \corref{cor-uni-Acat-K-2-cat} (resp.
\corref{cor-uni-Acat-ordin-2-cat}). {\color{red} Namely, all axioms of
a 2\n-category hold true except that the right action of a 1\n-morphism
on a unit 2\n-morphism does not necessarily give a unit 2\n-morphism.}
This makes clear the meaning of the statements `a natural
\ainf-transformation is invertible in $\sS{^u}A_\infty$' and `an
\ainf-functor is an equivalence in $\sS{^u}A_\infty$'
(\corref{cor-uni-Acat-ordin-2-cat}). We show in \propref{pro-rfg-p-r-1}
that a natural \ainf-transformation is invertible if and only if its
0\n-th component is invertible modulo boundary (in the sense of
\secref{sec-Invertible-trans}). The binary composition with a cycle
invertible modulo boundary is homotopy invertible
(\lemref{lem-r-p-inverse}).

In the eighth section we discuss unital \ainf-functors (between unital
\ainf-categories). Their first components map unit elements into unit
elements modulo boundary (\defref{def-Unital-functors}). For a unital
functor $f:\ca\to\cb$ we have an equivalence of natural
\ainf-transformations $\uni^\ca f\equiv f\uni^\cb$
(\propref{pro-unital-functors}). An \ainf-functor isomorphic to a
unital \ainf-functor is unital as well, see \eqref{eq-riBg-rgiC}.
Unital \ainf-categories, unital \ainf-functors and equivalence classes
of natural \ainf-transformations form a 1\n-2-unital 2\n-category
 $\ck A_\infty^u\subset\ck A_\infty$ enriched in $\ck$ and a
conventional 2\n-category $A_\infty^u\subset\sS{^u}A_\infty$
(\defref{def-Ainfuu}), which is a close analog of the 2\n-category of
usual categories. The construction of $A_\infty^u$ is the main point of
the article. It is needed for developing a theory of free
\ainf-categories, since it is expected that their universality
properties are formulated in the language of 2\n-categories.

There is a forgetful 2\n-functor $\kf:A_\infty^u\to\KCat$, which takes
a unital \ainf-category into the same differential quiver equipped with
the binary composition, viewed as a $\ck$\n-category
(\propref{pro-k-Ainf-Kcat}). The 2\n-functor $\kf$ reduces a unital
\ainf-functor to its first component and a natural \ainf-transformation
to its 0\n-th component. It turns out that many properties of an
\ainf-functor are determined by its first component and many properties
of a natural \ainf-transformation are determined by its 0\n-th
component. For instance, if the first component of an \ainf-functor
$\phi$ is homotopy invertible, then any natural \ainf-transformation
$y:f\phi\to g\phi$ is equivalent to $t\phi$ for a unique (up to
equivalence) natural \ainf-transformation $t:f\to g$ (Cancellation
\lemref{lem-cancel-phi}). If an \ainf-functor $\phi:\cc\to\cb$ to a
unital \ainf-category $\cb$ has homotopy invertible first component,
and each object of $\cb$ is ``isomorphic modulo boundary'' to an object
from $\phi(\Ob\cc)$, then $\phi$ is a unital equivalence, and $\cc$ is
unital (\thmref{thm-phi-C-B-hinv-equiv}). An equivalence between unital
\ainf-categories in $\sS{^u}A_\infty$ is always unital
(\corref{cor-equiv-unital}), which is not an immediate consequence of
definitions.

As a first example of a unital \ainf-category we list strictly unital
\ainf-categories (\secref{sec-Strict-unit-cat}), which is a well-known
notion. Other examples of unital \ainf-categories are obtained via
\thmref{thm-phi-C-B-hinv-equiv}. For instance, if an \ainf-functor
$\phi:\cc\to\cb$ to a strictly unital \ainf-category $\cb$ is
invertible, then $\cc$ is unital
(\secref{sec-Other-examples-unit-cat}). We stress again that taking the
0\n-th cohomology of a unital \ainf-category $\cc$ we get a
$\kk$\n-linear category $H^0(\cc)$. This $H^0$ can be viewed as a
2\n-functor (\secref{sec-H0-cohom-cat}).

In \appref{ap-sec-enrich} we define 2\n-categories enriched in a
symmetric monoidal category. Non-2-unital 2\n-categories are described
in \defref{def-non-unit-2cat}. 2\n-unital (usual) 2\n-categories admit
a concise \defref{def-2-cat-pack} and an expanded
\defref{def-unit-2cat}+\ref{def-non-unit-2cat}.

In \appref{ap-Contractibility} we prove that the cone of a homotopical
isomorphism is contractible.

\section{Conventions}\label{Conventions}
We fix a universe $\fu$ \cite[Sections~0,1]{GroVer-Prefai-SGA4a},
\cite{Bourb-Univers-SGA4a}. Many classes and sets in this paper will
mean $\fu$\n-small sets, even if not explicitly mentioned.

$\kk$ denotes a ($\fu$\n-small) unital associative commutative ring. By
abuse of notation it denotes also a chain complex, whose 0\n-th
component is $\kk$, and the other components vanish. A $\kk$\n-module
means a $\fu$\n-small $\kk$\n-module. The tensor product $\tens$
usually means $\tens_\kk$ -- the tensor product of graded
$\kk$\n-modules. It turns the category of graded $\kk$\n-modules into a
closed monoidal category. We will use its standard symmetry
$c:x\tens y\mapsto(-)^{xy}y\tens x=(-)^{\deg x\cdot\deg y}y\tens x$.
This paper contains many signs, and everywhere we abbreviate the usual
$(-1)^{(\deg x)(\deg y)}$ to $(-)^{xy}$. Similarly, $(-)^x$ means
$(-1)^{\deg x}$, or, simply, $(-1)^x$, if $x$ is an integer.

It is easy to understand the line
\[
\ca(X_0,X_1)\tens_\kk\ca(X_1,X_2)\tens_\kk\dots\tens_\kk\ca(X_{n-1},X_n),
\]
and it is much harder to understand the order in
\[
\ca(X_{n-1},X_n)\tens_\kk\dots\tens_\kk\ca(X_1,X_2)\tens_\kk\ca(X_0,X_1).
\]
That is why we use the right operators: the composition of two maps (or
morphisms) $f:X\to Y$ and $g:Y\to Z$ is denoted by $fg:X\to Z$. A map
is written on elements as $f:x\mapsto xf=(x)f$. However, these
conventions are not used systematically, and $f(x)$ might be used
instead.

When $f,g:X\to Y$ are chain maps, $f\sim g$ means that they are
homotopic. We denote by $\ck$ the category of differential graded
$\kk$\n-modules, whose morphisms are chain maps modulo homotopy. A
complex of $\kk$\n-modules $X$ is called \emph{contractible} if
$\id_X\sim0$, in other words, if $X$ is isomorphic to 0 in $\ck$.

If $C$ is a $\ZZ$\n-graded $\kk$\n-module, then its suspension
$sC=C[1]$ is the same $\kk$\n-module with the shifted grading
$(sC)^d=C^{d+1}$. The ``identity'' map $C\to sC$ of degree $-1$ is also
denoted by $s$. We follow the Getzler--Jones sign conventions
\cite{GetzlerJones:A-infty}, which include the idea to apply operations
to complexes with shifted grading, and Koszul's rule:
\[ (x\tens y)(f\tens g)=(-)^{yf}xf\tens yg
=(-1)^{\deg y\cdot\deg f}xf\tens yg.
\]
It takes its origin in Koszul's note~\cite{Koszul-CRAS47-signs}. The
main notions of graded algebra were given their modern names in
H.~Cartan's note~\cite{Cartan-CRAS48-graded}. See Boardman
\cite{Boardman-signs} for operad-like approach to signs as opposed to
closed symmetric monoidal category picture of Mac Lane
\cite{MacLane-Homology} (standard sign commutation rule). Combined
together, these sign conventions make the number of signs in this paper
tolerable.

If $u:A\to C$, $a\mapsto au$, is a chain map, its cone is the complex
$\Cone(u)=C\oplus A[1]$, $\Cone^k(u)=C^k\oplus A^{k+1}$, with the
differential $(c,a)d=(cd^C+au,ad^{A[1]})=(cd^C+au,-ad^A)$.

\section{\texorpdfstring{$A_\infty$-categories, $A_\infty$-functors and
 $A_\infty$-transformations}
{A8-categories, A8-functors and A8-transformations}}
\begin{definition}[Quiver]\label{def-Quiver}
A graded $\kk$\n-quiver $\ca$ consists of the following data: a class
of objects $\Ob\ca$ (a $\fu$\n-small set); a $\ZZ$\n-graded
$\kk$\n-module $\ca(X,Y)=\Hom_\ca(X,Y)$ for each pair of objects $X$,
$Y$. A morphism of $\kk$\n-quivers $f:\ca\to\cb$ is given by a map
$f:\Ob\ca\to\Ob\cb$, $X\mapsto Xf$ and by a $\kk$\n-linear map
$\ca(X,Y)\to\cb(Xf,Yf)$ for each pair of objects $X$, $Y$ of $\ca$.
\end{definition}

To a given graded $\kk$\n-quiver $\ca$ we associate another graded
$\kk$\n-quiver -- its tensor coalgebra $T\ca$, which has the same class
of objects as $\ca$. For each sequence $(X_0,X_1,X_2,\dots,X_n)$ of
objects of $\ca$ there is the $\ZZ$\n-graded $\kk$\n-module
$T^n\ca=
\ca(X_0,X_1)\tens_\kk\ca(X_1,X_2)\tens_\kk\dots\tens_\kk\ca(X_{n-1},X_n)$.
For the sequence $(X_0)$ with $n=0$ it reduces to $T^0\ca=\kk$ in
degree 0. The graded $\kk$\n-module $T\ca(X,Y)=\oplus_{n\ge0}T^n\ca$ is
the sum of the above modules over all sequences with $X_0=X$, $X_n=Y$.
The quiver $T\ca$ is equipped with the cut comultiplication
$\Delta:T\ca(X,Y)\to\oplus_{Z\in\Ob\ca}T\ca(X,Z)\bigotimes_\kk T\ca(Z,Y)$,
$h_1\tens h_2\tens\dots\tens h_n\mapsto \sum_{k=0}^n
h_1\tens\dots\tens h_k\bigotimes h_{k+1}\tens\dots\tens h_n$,
and the counit
$\eps=\bigl(T\ca(X,Y) \rTTo^{\pr_0} T^0\ca(X,Y)\to\kk\bigr)$, where the
last map is $\id_\kk$ if $X=Y$, or $0$ if $X\ne Y$ (and
$T^0\ca(X,Y)=0$). For this article it is the main example of the
following notion:

\begin{definition}[Cocategory]\label{def-Cocategory}
A graded cocategory $\cc$ is a graded $\kk$\n-quiver $\cc$, equipped
with a comultiplication -- a $\kk$\n-linear map
$\Delta^Z_{X,Y}:\cc(X,Y)\to\cc(X,Z)\tens_\kk\cc(Z,Y)$ of degree $0$ for
all triples $X$, $Y$, $Z$ of objects of $\cc$, and with a counit -- a
$\kk$\n-linear map $\eps_X:\cc(X,X)\to\kk$ of degree $0$ for all
objects $X$ of $\cc$, such that the usual associativity equations and
two counit equation hold.
\end{definition}

Associated to each graded cocategory $\cc$ is a graded
$\kk$\n-coalgebra $C=\oplus_{X,Y\in\Ob\cc}\cc(X,Y)$. Vice versa, to a
graded $\kk$\n-coalgebra, decomposed like that into $\kk$\n-submodules
$\cc(X,Y)$, $X,Y\in\Ob\cc$, for some $\fu$\n-small set $\Ob\cc$, so
that
$\Delta(\cc(X,Y))\subset\oplus_{Z\in\Ob\cc}\cc(X,Z)\tens_\kk\cc(Z,Y)$
for all pairs $X$, $Y$ of objects of $\cc$, and $\eps(\cc(X,Y))=0$ for
$X\ne Y$, we associate a graded cocategory.

This interpretation allows one to define a \emph{cocategory
homomorphism} $f:\cc\to\cd$ as a particular case of a coalgebra
homomorphism: a map $f:\Ob\cc\to\Ob\cd$, and $\kk$\n-linear maps
$\cc(X,Y)\to\cd(Xf,Yf)$ for all pairs of objects $X$, $Y$ of $\cc$,
compatible with comultiplication and counit. Given cocategory
homomorphisms $f,g:\cc\to\cd$ we say that a system of $\kk$\n-linear
maps $r:\cc(X,Y)\to\cd(Xf,Yg)$, $X,Y\in\Ob\cc$ is an
$(f,g)$-\emph{coderivation}, if the equation
$r\Delta=\Delta(f\tens r+r\tens g)$ holds.

In particular, these definitions apply to the tensor coalgebras
$Ts\ca=T(s\ca)$ of (suspended) $\kk$\n-quivers $s\ca$. In this case
cocategory homomorphisms and coderivations have a special form as we
shall see below.

\begin{definition}[\ainf-category, Kontsevich
\cite{Kontsevich:alg-geom/9411018}]\label{def-Ainf-category}
An \ainf-category $\ca$ consists of the following data: a graded
$\kk$\n-quiver $\ca$; a differential $b:Ts\ca\to Ts\ca$ of degree 1,
which is a (1,1)-coderivation, such that $(T^0s\ca)b=0$.
\end{definition}

The definition of a (1,1)-coderivation
$b\Delta=\Delta(1\tens b+b\tens1)$ implies that a $\kk$\n-quiver
morphism $b$ is determined by a system of $\kk$\n-linear maps
$b\pr_1:Ts\ca\to s\ca$ with components of degree 1
\[ b_n:s\ca(X_0,X_1)\tens s\ca(X_1,X_2)\tens\dots\tens s\ca(X_{n-1},X_n)
\to s\ca(X_0,X_n), \qquad n\ge1,
\]
via the formula
\begin{equation*}
b_{kl}=(b\big|_{T^ks\ca})\pr_l : T^ks\ca \to T^ls\ca, \qquad
b_{kl} = \sum_{\substack{r+n+t=k\\r+1+t=l}}
1^{\tens r}\tens b_n\tens1^{\tens t}.
\end{equation*}
Notice that the last condition of the definition implies $b_0=0$. In
particular, $b_{k0}=0$, and $k<l$ implies $b_{kl}=0$. Since $b^2$ is a
(1,1)-coderivation of degree 2, the equation $b^2=0$ is equivalent to
its particular case $b^2\pr_1=0$, that is, for all $k>0$
\begin{equation}
\sum_{r+n+t=k} (1^{\tens r}\tens b_n\tens1^{\tens t})b_{r+1+t} = 0 :
T^ks\ca \to s\ca.
\label{eq-b-b-0}
\end{equation}

Using another, more traditional, form of components of $b$:
\[ m_n = \bigl(\ca^{\tens n} \rTTo^{s^{\tens n}} (s\ca)^{\tens n}
\rTTo^{b_n} s\ca \rTTo^{s^{-1}} \ca\bigr) \]
we rewrite \eqref{eq-b-b-0} as follows:
\begin{equation}
\sum_{r+n+t=k} (-)^{t+rn}(1^{\tens r}\tens m_n\tens1^{\tens t})m_{r+1+t}
= 0 : T^k\ca \to \ca.
\label{eq-m-sbs-1}
\end{equation}
Notice that this equation differs in sign from \cite{math.RA/9910179},
because we are using right operators!

\begin{definition}[\ainf-functor, e.g. Keller \cite{math.RA/9910179}]
\label{def-ainf-functor}
An \ainf-functor $f:\ca\to\cb$ consists of the following data:
\ainf-categories $\ca$ and $\cb$, a cocategory homomorphism
$f:Ts\ca\to Ts\cb$ of degree 0, which commutes with the differential $b$
 {\color{blue}such that $(T^0s\ca)f\subset T^0s\cb$}.
\end{definition}

The definition of a cocategory homomorphism $f\Delta=\Delta(f\tens f)$,
$f\epsilon=\epsilon$ implies that $f$ is determined by a map
$f:\Ob\ca\to\Ob\cb$, $X\mapsto Xf$ and a system of $\kk$\n-linear maps
$f\pr_1:Ts\ca\to s\cb$ with components of degree 0
\[ f_n:s\ca(X_0,X_1)\tens s\ca(X_1,X_2)\tens\dots\tens s\ca(X_{n-1},X_n)
\to s\cb(X_0f,X_nf),
\]
$n\ge1$, (note that $f_0=0$) via the formula
\begin{equation}
f_{kl}=(f\big|_{T^ks\ca})\pr_l : T^ks\ca \to T^ls\cb, \qquad
f_{kl} = \sum_{i_1+\dots+i_l=k}
f_{i_1} \tens f_{i_2} \tens\dots\tens f_{i_l}.
\label{eq-f-kl-components}
\end{equation}

In particular, $f_{00}=\id:\kk\to\kk$, and $k<l$ implies $f_{kl}=0$.
Since $fb$ and $bf$ are both $(f,f)$-coderivations of degree 1, the
equation $fb=bf$ is equivalent to its particular case
$fb\pr_1=bf\pr_1$, that is, for all $k>0$
\begin{equation}
\sum_{l>0;i_1+\dots+i_l=k}
(f_{i_1} \tens f_{i_2} \tens\dots\tens f_{i_l}) b_l =
\sum_{r+n+t=k} (1^{\tens r}\tens b_n\tens1^{\tens t}) f_{r+1+t} :
T^ks\ca \to s\cb.
\label{eq-fff-b-b-f}
\end{equation}

Using $m_n$ we rewrite \eqref{eq-fff-b-b-f} as follows:
\begin{gather*}
\sum^{l>0}_{i_1+\dots+i_l=k} (-)^\sigma
(f_{i_1} \tens f_{i_2} \tens\dots\tens f_{i_l}) m_l =
\sum_{r+n+t=k} (-)^{t+rn}(1^{\tens r}\tens m_n\tens1^{\tens t}) f_{r+1+t}
: T^k\ca \to \cb, \\
\sigma = (i_2-1) + 2(i_3-1) +\dots+ (l-2)(i_{l-1}-1) + (l-1)(i_l-1).
\end{gather*}
Notice that this equation differs in sign from \cite{math.RA/9910179},
because we are using right operators.

\begin{example}
An \ainf-category with one object is called an \emph{\ainf-algebra}
(Stasheff~\cite{Stasheff:HomAssoc}). An \ainf-functor between
\ainf-algebras is called an \emph{\ainf-homomorphism}
(Kadeishvili~\cite{Kadeishvili82}). These notions are psychologically
easier to deal with, than the general case. The following notion of
\ainf-transformations also makes sense for \ainf-algebras, however,
such a context seems too narrow for it, because an \ainf-transformation
is an analog of a transformation between ordinary functors without the
naturality condition. Needless to say, in ordinary category theory
there is no reason to consider unnatural transformations. The reasons
to do it for \ainf-version are given in \secref{sec-Acat-of-Afuns}.
\end{example}

\begin{definition}[\ainf-transformation]\label{def-ainf-trans}
An \ainf-transformation $r:f\to g:\ca\to\cb$ of degree $d$ (pre natural
transformation in terms of \cite{Fukaya:FloerMirror-II}) consists of
the following data: \ainf-categories $\ca$ and $\cb$; \ainf-functors
$f,g:\ca\to\cb$; an $(f,g)$-coderivation $r:Ts\ca\to Ts\cb$ of degree
$d$.
\end{definition}

The definition of an $(f,g)$-coderivation
$r\Delta=\Delta(f\tens r+r\tens g)$ implies that $r$ is determined by a
system of $\kk$\n-linear maps $r\pr_1:Ts\ca\to s\cb$ with components of
degree $d$
\[ r_n:s\ca(X_0,X_1)\tens s\ca(X_1,X_2)\tens\dots\tens s\ca(X_{n-1},X_n)
\to s\cb(X_0f,X_ng),
\]
$n\ge0$, via the formula
\begin{align}
r_{kl} &= (r\big|_{T^ks\ca})\pr_l : T^ks\ca \to T^ls\cb, \notag \\
r_{kl} &= \sum_{\substack{q+1+t=l\\i_1+\dots+i_q+n+j_1+\dots+j_t=k}}
f_{i_1} \tens\dots\tens f_{i_q} \tens r_n\tens
g_{j_1} \tens\dots\tens g_{j_t}.
\label{eq-r-kl-components}
\end{align}
Note that $r_0$ is a system of $\kk$\n-linear maps
$\sS{_X}r_0:\kk\to s\cb(Xf,Xg)$, $X\in\Ob\ca$. In fact, the terms
`$A_\infty$\n-transformation' and `coderivation' are synonyms.

In particular, $r_{0l}$ vanishes unless $l=1$, and $r_{01}=r_0$. The
component $r_{kl}$ vanishes unless $1\le l\le k+1$.

\begin{examples}\label{exas-r-T1-T2}
1) The restriction of an \ainf-transformation $r$ to $T^1$ is
\[ r\big|_{T^1s\ca} = r_1\oplus[(f_1\tens r_0)+(r_0\tens g_1)], \]
where $r_1:s\ca(X,Y)\to s\cb(Xf,Yg)$,
\begin{align*}
f_1\tens r_0 &: s\ca(X,Y) = s\ca(X,Y)\tens\kk \rTTo^{f_1\tens r_0}
s\cb(Xf,Yf)\tens s\cb(Yf,Yg), \\
r_0\tens g_1 &: s\ca(X,Y) = \kk\tens s\ca(X,Y) \rTTo^{r_0\tens g_1}
s\cb(Xf,Xg)\tens s\cb(Xg,Yg).
\end{align*}

2) The restriction of an \ainf-transformation $r$ to $T^2$ is
\begin{multline*}
r\big|_{T^2s\ca} = r_2\oplus
[(f_2\tens r_0)+(f_1\tens r_1)+(r_1\tens g_1)+(r_0\tens g_2)] \oplus \\
\oplus [(f_1\tens f_1\tens r_0)+(f_1\tens r_0\tens g_1)
+(r_0\tens g_1\tens g_1)],
\end{multline*}
where $r_2:s\ca(X,Y)\tens s\ca(Y,Z)\to s\cb(Xf,Zg)$,
\begin{align*}
f_2\tens r_0 &: s\ca(X,Y)\tens s\ca(Y,Z)\tens\kk \to
s\cb(Xf,Zf)\tens s\cb(Zf,Zg), \\
f_1\tens r_1 &: s\ca(X,Y)\tens s\ca(Y,Z) \to s\cb(Xf,Yf)\tens s\cb(Yf,Zg), \\
r_1\tens g_1 &: s\ca(X,Y)\tens s\ca(Y,Z) \to s\cb(Xf,Yg)\tens s\cb(Yg,Zg), \\
r_0\tens g_2 &: \kk\tens s\ca(X,Y)\tens s\ca(Y,Z) \to
s\cb(Xf,Xg)\tens s\cb(Xg,Zg), \\
f_1\tens f_1\tens r_0 &: s\ca(X,Y)\tens s\ca(Y,Z)\tens\kk \to
s\cb(Xf,Yf)\tens s\cb(Yf,Zf)\tens s\cb(Zf,Zg), \\
f_1\tens r_0\tens g_1 &: s\ca(X,Y)\tens\kk\tens s\ca(Y,Z) \to
s\cb(Xf,Yf)\tens s\cb(Yf,Yg)\tens s\cb(Yg,Zg), \\
r_0\tens g_1\tens g_1 &: \kk\tens s\ca(X,Y)\tens s\ca(Y,Z) \to
s\cb(Xf,Xg)\tens s\cb(Xg,Yg)\tens s\cb(Yg,Zg).
\end{align*}
\end{examples}

The $\kk$\n-module of $(f,g)$-coderivations $r$ is
$\prod_{n=0}^\infty V_n$, where
{\color{blue} the product is taken in the category of
graded $\kk$\n-modules,}
\begin{equation}
V_n = \prod_{X_0,\dots,X_n\in\Ob\ca} \Hom_\kk\bigl(
s\ca(X_0,X_1)\tens\dots\tens s\ca(X_{n-1},X_n),s\cb(X_0f,X_ng)\bigr)
\label{eq-Vn-module-rn}
\end{equation}
is the graded $\kk$\n-module of $n$\n-th components $r_n$. It is
equipped with the differential $d:V_n\to V_n$, given by the following
formula
\begin{equation}
r_nd = r_nb_1 - (-)^{r_n} \sum_{\alpha+1+\beta=n}
(1^{\tens\alpha}\tens b_1\tens1^{\tens\beta}) r_n.
\label{eq-rnd-rnb1-b1rn}
\end{equation}

\section{Coderivations and cocategory homomorphisms}
Let $\ca$, $\cb$ be graded $\kk$\n-quivers, and let
$f^0,f^1,\dots,f^n:T\ca\to T\cb$ be cocategory homomorphisms. Consider
$n$ coderivations $r_1$, \dots, $r_n$ as in
\[ f^0 \rTTo^{r^1} f^1 \rTTo^{r^2} \dots f^{n-1} \rTTo^{r^n} f^n :
T\ca \to T\cb.
\]
We construct the following system of $\kk$\n-linear maps from these
data:
$\theta=(r^1\tens\dots\tens r^n)\theta:T\ca(X,Y)\to T\cb(Xf^0,Yf^n)$ of
degree $\deg r^1+\dots+\deg r^n$. Its components
$\theta_{kl}=\theta\big|_{T^k\ca}\pr_l:T^k\ca\to T^l\cb$ are given by
the following formula
\begin{equation}
\theta_{kl} = \sum f^0_{i^0_1}\tens\dots\tens f^0_{i^0_{m_0}}\tens
r^1_{j_1}\tens f^1_{i^1_1}\tens\dots\tens f^1_{i^1_{m_1}}\tens\dots\tens
r^n_{j_n}\tens f^n_{i^n_1}\tens\dots\tens f^n_{i^n_{m_n}},
\label{eq-theta-kl-sum}
\end{equation}
where summation is taken over all terms with
\[ m_0+m_1+\dots+m_n+n=l, \quad i^0_1+\dots+i^0_{m_0}+j_1+
i^1_1+\dots+i^1_{m_1}+\dots+j_n+i^n_1+\dots+i^n_{m_n} = k. \]
The component $\theta_{kl}$ vanishes unless $n\le l\le k+n$. If $n=0$,
we set $()\theta=f^0$. If $n=1$, the formula gives $(r^1)\theta=r^1$.

\begin{proposition}\label{pro-theta-D-D-theta-theta}
For each $n\ge0$ the map $\theta$ satisfies the equation
\begin{equation}
(r^1\tens r^2\tens\dots\tens r^n)\theta\Delta =
\Delta \sum_{k=0}^n (r^1\tens\dots\tens r^k)\theta \tens
(r^{k+1}\tens\dots\tens r^n)\theta.
\label{eq-theta-Delta}
\end{equation}
\end{proposition}

\begin{proof}
Let us write down the required equation using $\bigotimes$ to separate
the two copies of $T\cb$ in $T\cb\bigotimes T\cb$, and $\tens$ to
denote the multiplication in $T\cb$. We have to prove that for each
$n,x,y,z\ge0$
\begin{multline*}
\bigl(T^x\ca \rTTo^{(r^1\tens\dots\tens r^n)\theta_{x,y+z}}
T^{y+z}\cb \rTTo^{\Delta_{yz}} T^y\cb\bigotimes T^z\cb\bigr) \\
= \sum_{k=0}^n \sum_{u+v=x} \bigl(T^x\ca \rTTo^{\Delta_{uv}}
T^u\ca\bigotimes T^v\ca \rTTo^{(r^1\tens\dots\tens r^k)\theta_{uy}
\bigotimes(r^{k+1}\tens\dots\tens r^n)\theta_{vz}}
T^y\cb\bigotimes T^z\cb \bigr).
\end{multline*}
Substituting \eqref{eq-theta-kl-sum} for $\theta$ in the above equation
we come to an identity, which is proved by inspection. Indeed, skipping
all the intermediate steps, we get the following equation:
\begin{multline*}
\sum_{k=0}^n
\sum_{\substack{m_0+m_1+\dots+m_n+n=y+z\\
m_0+m_1+\dots+m_{k-1}+k-1<y\le m_0+m_1+\dots+m_k+k\\
i^0_1+\dots+i^0_{m_0}+j_1+i^1_1+\dots+i^1_{m_1}+\dots+j_n+i^n_1+\dots+i^n_{m_n}=x}}
\\
f^0_{i^0_1}\tens\dots\tens f^0_{i^0_{m_0}}\tens\dots\tens
r^k_{j_k}\tens f^k_{i^k_1}\tens\dots\tens f^k_{i^k_w}\bigotimes
f^k_{i^k_{w+1}}\tens\dots\tens f^k_{i^k_{m_k}}\tens r^{k+1}_{j_{k+1}}
\tens\dots\tens f^n_{i^n_1}\tens\dots\tens f^n_{i^n_{m_n}}
\\
= \sum_{k=0}^n \sum_{u+v=x}
\sum_{\substack{m_0+m_1+\dots+m_{k-1}+k+w=y\\
i^0_1+\dots+i^0_{m_0}+\dots+j_k+i^k_1+\dots+i^k_w=u}}
\;\;\sum_{\substack{t+m_{k+1}+\dots+m_n+n-k=z\\
l_1+\dots+l_t+j_{k+1}+\dots+i^n_1+\dots+i^n_{m_n}=v}}
\\
f^0_{i^0_1}\tens\dots\tens f^0_{i^0_{m_0}}\tens\dots\tens
r^k_{j_k}\tens f^k_{i^k_1}\tens\dots\tens f^k_{i^k_w}\bigotimes
f^k_{l_1}\tens\dots\tens f^k_{l_t}\tens r^{k+1}_{j_{k+1}}
\tens\dots\tens f^n_{i^n_1}\tens\dots\tens f^n_{i^n_{m_n}}.
\end{multline*}
In the left hand side $w$ denotes the expression
$y-(m_0+m_1+\dots+m_{k-1}+k)$ and lies in the interval $0\le w\le m_k$.
Identifying $m_k$ in the left hand side with $w+t$ in the right hand
side, we deduce that the both sides are equal.
\end{proof}
 {\color{blue}
An easy proof follows due to $f^p\Delta=\Delta(f^p\tens f^p)$ from the
formula
\[ \theta_{kl}
= \sum_{i_0+j_1+i_1+\dots+j_n+i_n=k}^{m_0+m_1+\dots+m_n+n=l}
f^0_{i_0m_0}\tens r^1_{j_1}\tens f^1_{i_1m_1}\tens\dots\tens
r^n_{j_n}\tens f^n_{i_nm_n}: T^k\ca \to T^l\cb,
\]
which is a version of \eqref{eq-theta-kl-sum}, where
$f^p_{im}:T^i\ca\to T^m\cb$ are matrix elements of $f^p$.
 }

\subsection{Cocategory homomorphisms.}\label{sec-cocat-hom}
Graded $\kk$\n-coalgebras form a symmetric monoidal category. The
tensor product $C\tens_\kk D$ of $\kk$\n-coalgebras $C$, $D$ is
equipped with the comultiplication
\( C\tens D \rTTo^{\Delta\tens\Delta} C\tens C\tens D\tens D
\rTTo^{1\tens c\tens1} C\tens D\tens C\tens D \),
using the standard symmetry $c$ of graded $\kk$\n-modules. Since graded
cocategories are, in fact, graded coalgebras with a special
decomposition, they also form a symmetric monoidal category gCoCat. If
$\cc$ and $\cd$ are cocategories, then the class of objects of their
tensor product $\cc\tens\cd$ is $\Ob\cc\times\Ob\cd$, and
$\cc\tens\cd(X\times U,Y\times W)=\cc(X,Y)\tens_\kk\cd(U,W)$.

Let $\phi:T\ca\tens T\cc\to T\cb$ be a cocategory homomorphism of
degree 0. It is determined uniquely by its composition with $\pr_1$,
that is, by a family $\phi\pr_1=(\phi_{nm})_{n,m\ge0}$,
$\phi_{nm}:T^n\ca\tens T^m\cc\to\cb$, $\phi_{00}=0$, with the same
underlying map of objects $\Ob\ca\times\Ob\cc\to\Ob\cb$. Indeed, for
given families of composable arrows
$f^0 \rTTo^{p^1} f^1 \rTTo^{p^2} \dots f^{n-1} \rTTo^{p^n} f^n$ of
$\ca$ and
$g^0 \rTTo^{t^1} g^1 \rTTo^{t^2} \dots g^{m-1} \rTTo^{t^m} g^m$ of
$\cc$ we have
\begin{multline}
(p^1\tens\dots\tens p^n\tens t^1\tens\dots\tens t^m)\phi \\
= \sum_{\substack{i_1+\dots+i_k=n\\j_1+\dots+j_k=m}} (-)^\sigma
(p^1\tens\dots\tens p^{i_1}\tens t^1\tens\dots\tens t^{j_1})\phi_{i_1j_1}
\tens (p^{i_1+1}\tens\dots\tens p^{i_1+i_2}\tens
t^{j_1+1}\tens\dots\tens t^{j_1+j_2})\phi_{i_2j_2} \\
\tens\dots\tens (p^{i_1+\dots+i_{k-1}+1}\tens\dots\tens
p^{i_1+\dots+i_k} \tens t^{j_1+\dots+j_{k-1}+1} \tens\dots\tens
t^{j_1+\dots+j_k})\phi_{i_kj_k}.
\label{eq-pt-phi=phi-nm}
\end{multline}
The sign depends on the parity of an integer
\begin{multline}
\sigma = (t^1+\dots+t^{j_1}) (p^{i_1+1}+\dots+p^{i_1+\dots+i_k}) +
(t^{j_1+1}+\dots+t^{j_1+j_2}) (p^{i_1+i_2+1}+\dots+p^{i_1+\dots+i_k}) \\
+\dots+ (t^{j_1+\dots+j_{k-2}+1}+\dots+t^{j_1+\dots+j_{k-1}})
(p^{i_1+\dots+i_{k-1}+1}+\dots+p^{i_1+\dots+i_k}),
\label{eq-sign-sigma}
\end{multline}
where each coderivation has to be replaced with its degree. Recall that
in our notation (\secref{Conventions}) we abbreviate
$(-1)^{(\deg t)(\deg p)}$ to $(-)^{tp}$. By definition the homomorphism
$\phi$ satisfies the equation
\begin{diagram}[LaTeXeqno]
T\ca\tens T\cc & \rTTo^\phi & T\cb & \rTTo^\Delta & T\cb\tens T\cb \\
\dTTo<{\Delta\tens\Delta} &&&& \uTTo>{\phi\tens\phi} \\
T\ca\tens T\ca\tens T\cc\tens T\cc && \rTTo^{1\tens c\tens1} &&
T\ca\tens T\cc\tens T\ca\tens T\cc
\label{dia-phi-homomorphism}
\end{diagram}

Introduce $\kk$\n-linear maps
$(t^1\tens\dots\tens t^m)\chi:T\ca\to T\cb$ by the formula
$a[(t^1\tens\dots\tens t^m)\chi]=(a\tens t^1\tens\dots\tens t^m)\phi$,
$a\in T\ca$. Then the above equation is equivalent to
\begin{equation}
(t^1\tens t^2\tens\dots\tens t^m)\chi\Delta =
\Delta \sum_{k=0}^m (t^1\tens\dots\tens t^k)\chi \tens
(t^{k+1}\tens\dots\tens t^m)\chi.
\label{eq-chi-Delta}
\end{equation}
for all $m\ge0$.

When $\ca$, $\cb$ are graded $\kk$\n-quivers, we define a new
$\kk$\n-quiver $\Coder(\ca,\cb)$, whose objects are cocategory
homomorphisms $f:T\ca\to T\cb$. These homomorphisms are determined by a
system $f\pr_1=(f_n)_{n\ge1}$ of morphisms of $\kk$\n-quivers
$f_n:T^n\ca\to\cb$ of degree $0$ with the same underlying map
$\Ob\ca\to\Ob\cb$, see \eqref{eq-f-kl-components}. The $\kk$\n-module
of morphisms between $f,g:T\ca\to T\cb$ consists of
$(f,g)$-coderivations:
\[ [\Coder(\ca,\cb)(f,g)]^d =
\{ r:T\ca\to T\cb \mid r\Delta = \Delta(f\tens r + r\tens g),
\quad \deg r = d \}, \quad d\in\ZZ.
\]
Such a coderivation $r$ is determined by a system of $\kk$\n-linear maps
$r\pr_1=(r_n)_{n\ge0}$, $r_n:T^n\ca(X,Y)\to\cb(Xf,Yg)$ of
degree $d$ as in \eqref{eq-r-kl-components}.

\begin{corollary}[to \propref{pro-theta-D-D-theta-theta}]
\label{cor-action-alpha}
A map $\alpha:T\ca\tens T\Coder(\ca,\cb)\to T\cb$,
$a\tens r^1\tens\dots\tens r^n\mapsto a[(r^1\tens\dots\tens r^n)\theta]$,
is a cocategory homomorphism of degree 0.
\end{corollary}

\begin{proof}
Equation~\eqref{eq-theta-Delta} means that
equation~\eqref{eq-chi-Delta} holds for $\chi=\theta$, which is
equivalent to \eqref{dia-phi-homomorphism} for $\phi=\alpha$,
$\cc=\Coder(\ca,\cb)$.
\end{proof}

\begin{proposition}\label{prop-phi-C1-Cq-psi}
For any cocategory homomorphism
$\phi:T\ca\tens T\cc^1\tens T\cc^2\tens\dots\tens T\cc^q\to T\cb$ of
degree 0 there is a unique cocategory homomorphism
$\psi:T\cc^1\tens T\cc^2\tens\dots\tens T\cc^q\to T\Coder(\ca,\cb)$
of degree 0, such that
\[ \phi = \bigl(T\ca\tens T\cc^1\tens T\cc^2\tens\dots\tens T\cc^q
\rTTo^{1\tens\psi} T\ca\tens T\Coder(\ca,\cb) \rTTo^\alpha T\cb\bigr).
\]
\end{proposition}

\begin{proof}
Let us start with a simple case $q=1$, $\cc=\cc^1$. Each object $g$ of
$\cc$ induces a cocategory morphism $g\psi:a\mapsto(a\tens g)\phi$. We
set $\psi_0=0$. Each element $p\in\cc(g,h)$ induces a coderivation
$(p)\psi_1=p\psi:a\mapsto(a\tens p)\phi$. Suppose that $\psi_i$ are
already found for $0\le i<n$. Then we find $\psi_n$ from the sought
identity $\chi=\psi\theta$. Namely, for
$g^0 \rTTo^{p^1} g^1 \rTTo^{p^2} \dots g^{n-1} \rTTo^{p^n} g^n$ we have
to satisfy the identity
\[ (p^1\tens\dots\tens p^n)\chi = (p^1\tens\dots\tens p^n)\psi_n +
\sum_{l=2}^n \sum_{i_1+\dots+i_l=n} [(p^1\tens\dots\tens p^n).
(\psi_{i_1}\tens\psi_{i_2}\tens\dots\tens\psi_{i_l})]\theta,
\]
which expresses $\psi$ via its components $\psi_k$. Notice that the
unknown $\psi_n$ occurs only in the singled out summand, corresponding
to $l=1$. The factors $\psi_i$ in the sum are already known, since
$i<n$. So we define $(p^1\tens\dots\tens p^n)\psi_n:T\ca\to T\cb$ as
the difference of $(p^1\tens\dots\tens p^n)\chi$ and the sum in the
right hand side. Assume that $\psi$ is a cocategory homomorphism up to
the level $n$, that is,
\begin{equation}
(p^1\tens\dots\tens p^m)\psi = \sum_{l=1}^m \sum_{i_1+\dots+i_l=m}
(p^1\tens\dots\tens p^m).
(\psi_{i_1}\tens\psi_{i_2}\tens\dots\tens\psi_{i_l})
\label{eq-p1-pm-psi}
\end{equation}
for all $0\le m\le n$. Taking into account equations
\eqref{eq-theta-Delta}, we see that \eqref{eq-chi-Delta} is equivalent
to an equation of the form
\[ (p^1\tens\dots\tens p^n)\psi_n\Delta =
\Delta [g^0\psi\tens(p^1\tens\dots\tens p^n)\psi_n +
(p^1\tens\dots\tens p^n)\psi_n\tens g^n\psi + \mu]. \]
Moreover, if $(p^1\tens\dots\tens p^n)\psi_n$ were a
$(g^0\psi,g^n\psi)$-coderivation, it would imply \eqref{eq-chi-Delta}
by \secref{sec-cocat-hom}. We deduce that, indeed, the above $\mu=0$,
and $(p^1\tens\dots\tens p^n)\psi_n$ is a
$(g^0\psi,g^n\psi)$-coderivation. Thus, we have found a unique
$(p^1\tens\dots\tens p^n)\psi_n\in\Coder(\ca,\cb)$ and
\eqref{eq-p1-pm-psi} for $m=n$ defines uniquely an element
$(p^1\tens\dots\tens p^n)\psi\in T\Coder(\ca,\cb)$.

The case $q>1$ is similar to the case $q=1$, however, the reasoning is
slightly obstructed by a big amount of indices. So we explain in
detail the case $q=2$ only, and in the general case no new phenomena
occur. Further we shall use the obtained formulas in the case $q=2$.

Let $\cc=\cc^1$, $\cd=\cc^2$. We consider a cocategory homomorphism
$\phi:T\ca\tens T\cc\tens T\cd\to T\cb$, $a\mapsto a[(c\tens d)\chi]$.
We have to obtain from it a unique cocategory homomorphism
$\psi:T\cc\tens T\cd\to T\Coder(\ca,\cb)$.

A pair of objects $f\in\Ob\cc$, $g\in\Ob\cd$ induces a cocategory
morphism $(f,g)\psi:a\mapsto(a\tens1_f\tens1_g)\phi$. We set
$\psi_{00}=0$. An object $f\in\Ob\cc$ and an element $t\in\cd(g^0,g^1)$
induce an $((f,g^0)\psi,(f,g^1)\psi)$-coderivation
$(f\tens t)\psi_{01}=(f\tens t)\psi:a\mapsto(a\tens1_f\tens t)\phi$. An
element $p\in\cc(f^0,f^1)$ and an object $g\in\Ob\cd$ induce an
$((f^0,g)\psi,(f^1,g)\psi)$-coderivation
$(p\tens g)\psi_{10}=(p\tens g)\psi:a\mapsto(a\tens p\tens1_g)\phi$.
Suppose that $\psi_{ij}$ are already found for $0\le i\le n$,
$0\le j\le m$, $(i,j)\ne(n,m)$. Then we find $\psi_{nm}$ from the
sought identity $\chi=\psi\theta$. Namely, for
$f^0 \rTTo^{p^1} f^1 \rTTo^{p^2} \dots f^{n-1} \rTTo^{p^n} f^n$ in
$\cc$ and
$g^0 \rTTo^{t^1} g^1 \rTTo^{t^2} \dots g^{m-1} \rTTo^{t^m} g^m$ in
$\cd$ we have to satisfy the identity
\begin{multline}
(p^1\tens\dots\tens p^n\tens t^1\tens\dots\tens t^m)\chi =
(p^1\tens\dots\tens p^n\tens t^1\tens\dots\tens t^m)\psi_{nm} \\
+\sum_{\substack{i_1+\dots+i_k=n\\j_1+\dots+j_k=m}}^{k>1}(-)^\sigma\bigl[
(p^1\tens\dots\tens p^{i_1}\tens t^1\tens\dots\tens t^{j_1})\psi_{i_1j_1}
\tens (p^{i_1+1}\tens\dots\tens p^{i_1+i_2}\tens
t^{j_1+1}\tens\dots\tens t^{j_1+j_2})\psi_{i_2j_2} \\
\tens\dots\tens (p^{i_1+\dots+i_{k-1}+1}\tens\dots\tens
p^{i_1+\dots+i_k}\tens t^{j_1+\dots+j_{k-1}+1} \tens\dots\tens
t^{j_1+\dots+j_k})\psi_{i_kj_k}\bigr]\theta,
\label{eq-pt-chi=psi-theta}
\end{multline}
which is nothing else but \eqref{eq-pt-phi=phi-nm}. The sign is
determined by \eqref{eq-sign-sigma}. All terms of the sum are already
known. So we define a map
$(p^1\tens\dots\tens p^n\tens t^1\tens\dots\tens t^m)\psi_{nm}:T\ca\to T\cb$
as the difference of the left hand side and the sum in the right hand
side. The fact that $\phi$ is a homomorphism is equivalent to the
identity for the map $\chi$ for all $n,m\ge0$:
\begin{multline}
(p^1\tens\dots\tens p^n\tens t^1\tens\dots\tens t^m)\chi\Delta =
\Delta\sum_{k=0}^n\sum_{l=0}^m (-)^{(p^{k+1}+\dots+p^n)(t^1+\dots+t^l)} \\
(p^1\tens\dots\tens p^k\tens t^1\tens\dots\tens t^l)\chi \tens
(p^{k+1}\tens\dots\tens p^n\tens t^{l+1}\tens\dots\tens t^m)\chi,
\label{eq-p-t-chi-delta}
\end{multline}
similarly to \secref{sec-cocat-hom}. From it we get an equation for
$(p^1\tens\dots\tens p^n\tens t^1\tens\dots\tens t^m)\psi_{nm}$
\begin{multline*}
(p^1\tens\dots\tens p^n\tens t^1\tens\dots\tens t^m)\psi_{nm}\Delta
= \Delta\bigl[(f^0,g^0)\psi\tens
(p^1\tens\dots\tens p^n\tens t^1\tens\dots\tens t^m)\psi_{nm} \\
+ (p^1\tens\dots\tens p^n\tens t^1\tens\dots\tens t^m)\psi_{nm}
\tens(f^n,g^m)\psi + \mu\bigr].
\end{multline*}
Notice that if $\psi$ is indeed a homomorphism and $\psi_{nm}$ is its
component, then $\phi$ is a homomorphism, hence,
\eqref{eq-p-t-chi-delta} holds. Thus, the above equation with $\mu=0$
implies \eqref{eq-p-t-chi-delta} (and it happens only for one value of
$\mu$). Since we know that \eqref{eq-p-t-chi-delta} holds, it implies
$\mu=0$. Therefore,
$(p^1\tens\dots\tens p^n\tens t^1\tens\dots\tens t^m)\psi_{nm}$ is a
$((f^0,g^0)\psi,(f^n,g^m)\psi)$-coderivation. Thus, we have found a
unique
$(p^1\tens\dots\tens p^n\tens t^1\tens\dots\tens t^m)\psi_{nm}
\in\Coder(\ca,\cb)$, and
\begin{multline*}
(p^1\tens\dots\tens p^n\tens t^1\tens\dots\tens t^m)\psi =
(p^1\tens\dots\tens p^n\tens t^1\tens\dots\tens t^m)\psi_{nm} \\
+ \sum_{\substack{i_1+\dots+i_k=n\\j_1+\dots+j_k=m}}^{k>1} (-)^\sigma
(p^1\tens\dots\tens p^{i_1}\tens t^1\tens\dots\tens t^{j_1})\psi_{i_1j_1}
\tens (p^{i_1+1}\tens\dots\tens p^{i_1+i_2}\tens
t^{j_1+1}\tens\dots\tens t^{j_1+j_2})\psi_{i_2j_2} \\
\tens\dots\tens (p^{i_1+\dots+i_{k-1}+1}\tens\dots\tens
p^{i_1+\dots+i_k}\tens t^{j_1+\dots+j_{k-1}+1} \tens\dots\tens
t^{j_1+\dots+j_k})\psi_{i_kj_k}
\end{multline*}
defines uniquely an element of $T\Coder(\ca,\cb)$. The above formula
implies that $\psi$ is a homomorphism.

A generalization to $q>2$ is straightforward.
\end{proof}

We interpret the above proposition as the existence of inner
hom-objects $\HOM(T\ca,T\cb)=T\Coder(\ca,\cb)$ in the monoidal category
of cocategories of the form $T\cc^1\tens T\cc^2\tens\dots\tens T\cc^r$.

\section{A category enriched in cocategories}
\label{sec-enr-cat-grad-centr-bimod}
Let us show that the category of tensor coalgebras of graded
$\kk$\n-quivers is enriched in gCoCat.

Let $\ca$, $\cb$, $\cc$ be graded $\kk$\n-quivers. Consider the
cocategory homomorphism given by the upper right path in the diagram
\begin{diagram}[LaTeXeqno]
T\ca\tens T\Coder(\ca,\cb)\tens T\Coder(\cb,\cc) & \rTTo^{\alpha\tens1}
& T\cb\tens T\Coder(\cb,\cc) \\
\dTTo<{1\tens M} & = & \dTTo>\alpha \\
T\ca\tens T\Coder(\ca,\cc) & \rTTo^\alpha & T\cc
\label{dia-def-M}
\end{diagram}
By \propref{prop-phi-C1-Cq-psi} there is a graded cocategory morphism
of degree 0
\[ M: T\Coder(\ca,\cb)\tens T\Coder(\cb,\cc)
\to T\Coder(\ca,\cc). \]

Denote by $\1$ a graded 1-object-0-morphisms $\kk$\n-quiver, that is,
$\Ob\1=\{*\}$, $\1(*,*)=0$. Then $T\1=\kk$ is a unit object of the
monoidal category of graded cocategories. Denote by
$\rightunit:T\ca\tens T\1\to T\ca$ and
$\leftunit:T\1\tens T\ca\to T\ca$ the corresponding natural cocategory
isomorphisms. By \propref{prop-phi-C1-Cq-psi} there exists a unique
cocategory morphism $\eta_\ca:T\1\to T\Coder(\ca,\ca)$, such that
\[ \rightunit = \bigl(T\ca\tens T\1 \rTTo^{1\tens\eta_\ca}
T\ca\tens T\Coder(\ca,\ca) \rTTo^\alpha T\ca\bigr).
\]
Namely, the object $*\in\Ob\1$ goes to the identity homomorphism
$\id_\ca:\ca\to\ca$, which acts as the identity map on objects, and has
only one non-vanishing component $(\id_\ca)_1=\id:\ca(X,Y)\to\ca(X,Y)$.

\begin{proposition}[See also Kontsevich and Soibelman
\cite{KonSoi-book}]\label{pro-M-assoc-eta-unit}
The multiplication $M$ is associative and $\eta$ is its two-sided unit:
\begin{diagram}[nobalance]
T\Coder(\ca,\cb)\tens T\Coder(\cb,\cc)\tens T\Coder(\cc,\cd)
& \rTTo^{M\tens1} & T\Coder(\ca,\cc)\tens T\Coder(\cc,\cd) \\
\dTTo<{1\tens M} && \dTTo>M \\
T\Coder(\ca,\cb)\tens T\Coder(\cb,\cd) & \rTTo^M
& T\Coder(\ca,\cd)
\end{diagram}
\end{proposition}

\begin{proof}
The cocategory homomorphism
\begin{multline*}
T\ca\tens T\Coder(\ca,\cb)\tens T\Coder(\cb,\cc)\tens T\Coder(\cc,\cd) \\
\rTTo^{\alpha\tens1\tens1} T\cb\tens T\Coder(\cb,\cc)\tens T\Coder(\cc,\cd)
\rTTo^{\alpha\tens1} T\cc\tens T\Coder(\cc,\cd) \rTTo^\alpha T\cd
\end{multline*}
can be written down as
\begin{equation*}
(\alpha\tens1\tens1)(1\tens M)\alpha
= (1\tens1\tens M)(\alpha\tens1)\alpha
= (1\tens1\tens M)(1\tens M)\alpha,
\end{equation*}
or as
\begin{equation*}
(1\tens M\tens1)(\alpha\tens1)\alpha
= (1\tens M\tens1)(1\tens M)\alpha.
\end{equation*}
The uniqueness part of \propref{prop-phi-C1-Cq-psi} implies that
$(1\tens M)M=(M\tens1)M$.

Similarly one proves that $\eta$ is a unit for $M$.
\end{proof}

By \eqref{eq-f-kl-components} we find
\[ \text{gCoCat}(T\1,T\Coder(\ca,\cb)) =
\text{Maps}(\{*\},\Ob\Coder(\ca,\cb)) = \text{gCoCat}(T\ca,T\cb).
\]
Thus we can interpret \propref{pro-M-assoc-eta-unit} as saying that the
category of tensor coalgebras of graded $\kk$\n-quivers admits an
enrichment in gCoCat.

Let us find explicit formulas for $M$. It is defined on objects as
composition: if $f:\ca\to\cb$ and $g:\cb\to\cc$ are cocategory
morphisms, then $(f,g)M=fg:\ca\to\cc$. On coderivations $M$ is
specified by its composition with
$\pr_1:T\Coder(\ca,\cc)\to\Coder(\ca,\cc)$. Let us write in this
section $(h,k)$ as a shorthand for $\Coder(\ca,\cc)(h,k)$, the
$\kk$\n-module of $(h,k)$-coderivations. The components of $M$ are
\begin{gather*}
M_{nm}=M\big|_{T^n\tens T^m}\pr_1:
T^n\Coder(\ca,\cb)\tens T^m\Coder(\cb,\cc) \to \Coder(\ca,\cc), \\
M_{nm}: (f^0,f^1)\tens\dots\tens(f^{n-1},f^n)
\tens(g^0,g^1)\tens\dots\tens(g^{m-1},g^m) \to (f^0g^0,f^ng^m),
\end{gather*}
where $f^0,\dots,f^n:\ca\to\cb$ and $g^0,\dots,g^m:\cb\to\cc$ are
cocategory morphisms. We have $M_{00}=0$.

According to proof of \propref{prop-phi-C1-Cq-psi} the component
$M_{nm}$ is determined recursively from
equation~\eqref{eq-pt-chi=psi-theta}:
\begin{multline}
(p^1\tens\dots\tens p^n)\theta(t^1\tens\dots\tens t^m)\theta =
(p^1\tens\dots\tens p^n\tens t^1\tens\dots\tens t^m)M_{nm} \\
+ \sum_{\substack{i_1+\dots+i_k=n\\j_1+\dots+j_k=m}}^{k>1}
\hspace*{-2mm} (-)^\sigma \bigl[
(p^1\tens\dots\tens p^{i_1}\tens t^1\tens\dots\tens t^{j_1})M_{i_1j_1}
\tens (p^{i_1+1}\tens\dots\tens p^{i_1+i_2}\tens
t^{j_1+1}\tens\dots\tens t^{j_1+j_2})M_{i_2j_2} \\
\tens\dots\tens (p^{i_1+\dots+i_{k-1}+1}\tens\dots\tens
p^{i_1+\dots+i_k}\tens t^{j_1+\dots+j_{k-1}+1} \tens\dots\tens
t^{j_1+\dots+j_k})M_{i_kj_k}\bigr]\theta.
\label{eq-theta-theta=M-nm}
\end{multline}
Since $(p^1\tens\dots\tens p^n\tens t^1\tens\dots\tens t^m)M_{nm}$ is a
coderivation, it is determined by its composition with projection
$\pr_1$. Composing \eqref{eq-theta-theta=M-nm} with $\pr_1$ we get
\begin{equation}
(p^1\tens\dots\tens p^n)\theta(t^1\tens\dots\tens t^m)\theta\pr_1 =
(p^1\tens\dots\tens p^n\tens t^1\tens\dots\tens t^m)M_{nm}\pr_1.
\label{eq-theta-theta-pr1}
\end{equation}
Therefore, if $m=0$ and $n$ is positive, $M_{n0}$ is given by the
formula:
\begin{align*}
M_{n0}: (f^0,f^1)\tens\dots\tens(f^{n-1},f^n)\tens\kk_{g^0} &\to
(f^0g^0,f^ng^0), \\
r^1\tens\dots\tens r^n\tens1 &\mapsto
(r^1\tens\dots\tens r^n\mid g^0)M_{n0},
\end{align*}
\[ [(r^1\tens\dots\tens r^n\mid g^0)M_{n0}]\pr_1 =
(r^1\tens\dots\tens r^n)\theta g^0\pr_1, \]
where $\mid$ separates the arguments in place of $\tens$.
If $m=1$, then $M_{n1}$ is given by the formula:
\begin{align*}
M_{n1}: (f^0,f^1)\tens\dots\tens(f^{n-1},f^n)\tens(g^0,g^1) &\to
(f^0g^0,f^ng^1), \\
r^1\tens\dots\tens r^n\tens t^1 &\mapsto
(r^1\tens\dots\tens r^n\tens t^1)M_{n1},
\end{align*}
\[ [(r^1\tens\dots\tens r^n\tens t^1)M_{n1}]\pr_1 =
(r^1\tens\dots\tens r^n)\theta t^1\pr_1. \]
Explicitly we write
\begin{align}
[(r^1\tens\dots\tens r^n\mid g^0)M_{n0}]_k &=
\sum_l (r^1\tens\dots\tens r^n)\theta_{kl} g^0_l,
\label{eq-Mn0-k-component} \\
[(r^1\tens\dots\tens r^n\tens t^1)M_{n1}]_k &=
\sum_l (r^1\tens\dots\tens r^n)\theta_{kl} t^1_l.
\label{eq-Mn1-k-component}
\end{align}
Finally, $M_{nm}=0$ for $m>1$, since the left hand side of
\eqref{eq-theta-theta-pr1} vanishes.

\begin{examples}
1) The component $M_{01}$ is the composition:
$(f^0\mid t^1)M_{01}=f^0t^1$.

2) The component $M_{10}$ is the composition:
$(r^1\mid g^0)M_{10}=r^1g^0$.

3) If $r:f\to g:\ca\to\cb$ and $p:h\to k:\cb\to\cc$ are
\ainf-transformations, then $(r\tens p)M_{11}:fh\to gk:\ca\to\cc$ has
the following components:
\begin{align*}
[(r\tens p)M_{11}]_0 &= r_0p_1, \\
[(r\tens p)M_{11}]_1 &= r_1p_1 + (f_1\tens r_0)p_2 + (r_0\tens g_1)p_2,
\text{ etc.}
\end{align*}

4) If $f \rTTo^r g \rTTo^p h:\ca\to\cb$ are \ainf-transformations, and
$k:\cb\to\cc$ is an \ainf-functor, then
$(r\tens p\mid k)M_{20}:fk\to hk:\ca\to\cc$ has the following
components:
\begin{align*}
[(r\tens p\mid k)M_{20}]_0 &= (r_0\tens p_0)k_2, \\
[(r\tens p\mid k)M_{20}]_1 &= (r_1\tens p_0)k_2 + (r_0\tens p_1)k_2 \\
&\quad+ (r_0\tens p_0\tens h_1)k_3 + (r_0\tens g_1\tens p_0)k_3 +
(f_1\tens r_0\tens p_0)k_3,
\text{ etc.}
\end{align*}

5) If $f \rTTo^r g \rTTo^p h:\ca\to\cb$ and $t:k\to l:\cb\to\cc$ are
\ainf-transformations, then
$(r\tens p\tens t)M_{21}:fk\to hl:\ca\to\cc$ has the following
components:
\begin{align*}
[(r\tens p\tens t)M_{21}]_0 &= (r_0\tens p_0)t_2, \\
[(r\tens p\tens t)M_{21}]_1 &= (r_1\tens p_0)t_2 + (r_0\tens p_1)t_2 \\
&\quad+ (r_0\tens p_0\tens h_1)t_3 + (r_0\tens g_1\tens p_0)t_3 +
(f_1\tens r_0\tens p_0)t_3,
\text{ etc.}
\end{align*}
\end{examples}

\section{\texorpdfstring{$A_\infty$-category of $A_\infty$-functors}
{A8-category of A8-functors}}
 \label{sec-Acat-of-Afuns}
Let us construct a new \ainf-category $A_\infty(\ca,\cb)$ out of given
two $\ca$ and $\cb$. Its underlying graded $\kk$\n-quiver is a full
subquiver of $s^{-1}\Coder(s\ca,s\cb)$. The objects of
$A_\infty(\ca,\cb)$ are \ainf-functors $f:\ca\to\cb$. Given two such
functors $f,g:\ca\to\cb$ we define the graded $\kk$\n-module
$A_\infty(\ca,\cb)(f,g)$ as the space of all \ainf-transformations
$r:f\to g$, namely,
\[ [A_\infty(\ca,\cb)(f,g)]^{d+1} =
\{ r:f\to g \mid \text{ \ainf-transformation }
r:Ts\ca\to Ts\cb \text{ has degree } d \}.
\]
In this section we use the notation
$(f,g)=sA_\infty(\ca,\cb)(f,g)=\Coder(s\ca,s\cb)(f,g)$ for the sake of
brevity. The degree of $r$ as an element of $(f,g)$ will be exactly
$d$:
\[ (f,g)^d = \{ r:f\to g \mid \text{ \ainf-transformation }
r:Ts\ca\to Ts\cb \text{ has degree } d \}. \]
We will use only this (natural) degree of $r$ in order to permute it
with other things by Koszul's rule.

Notice that even if $\ca$, $\cb$ have one object (and are
\ainf-algebras), the quiver $A_\infty(\ca,\cb)$ has several objects.
Thus theory of \ainf-algebras leads to the theory of \ainf-categories.

\begin{proposition}[See also Fukaya \cite{Fukaya:FloerMirror-II},
Kontsevich and Soibelman \cite{KonSoi-AinfCat-NCgeom,KonSoi-book} and
Le\-f\`evre-Ha\-se\-ga\-wa \cite{Lefevre-Ainfty-these}]
\label{pro-B-F-KS-LH}
Let $\ca$, $\cb$ be \ainf-categories. Then there exists a unique
(1,1)-coderivation $B:TsA_\infty(\ca,\cb)\to TsA_\infty(\ca,\cb)$ of
degree 1, such that $B_0=0$ and
\begin{equation}
(r^1\tens\dots\tens r^n)\theta b = [(r^1\tens\dots\tens r^n)B]\theta
+(-)^{r^1+\dots+r^n}b(r^1\tens\dots\tens r^n)\theta
\label{eq-theta-b-B-theta}
\end{equation}
for all $n\ge0$,
$r^1\tens\dots\tens r^n\in(f^0,f^1)\tens\dots\tens(f^{n-1},f^n)$. It
satisfies $B^2=0$, thus, it gives an \ainf-structure of
$A_\infty(\ca,\cb)$.
\end{proposition}

\begin{proof}
For $n=0$ \eqref{eq-theta-b-B-theta} reads as $f^0b=(f^0)B+bf^0$,
hence, $(f^0)B=f^0b-bf^0=0$. In particular, we may set $B_0=0$. Assume
that the coderivation components $B_j$ for $j<n$ are already found, so
that \eqref{eq-theta-b-B-theta} is satisfied up to $n-1$ arguments. Let
us determine a $\kk$\n-linear map
$(r^1\tens\dots\tens r^n)B_n:Ts\ca\to Ts\cb$ from
equation~\eqref{eq-theta-b-B-theta}, rewritten as follows:
\begin{multline}
(r^1\tens\dots\tens r^n)B_n = (r^1\tens\dots\tens r^n)\theta b
-(-)^{r^1+\dots+r^n}b(r^1\tens\dots\tens r^n)\theta \\
-\sum_{q+j+t=n}^{j<n}
[(r^1\tens\dots\tens r^n)(1^{\tens q}\tens B_j\tens1^{\tens t})]\theta.
\label{eq-Bn-bb-sum}
\end{multline}

Let us show that $(r^1\tens\dots\tens r^n)B_n$ is a
$(f^0,f^n)$-coderivation. Indeed,
\begin{align*}
(r^1 & \tens\dots\tens r^n)B_n\Delta_\cb =
(r^1\tens\dots\tens r^n)\theta b\Delta_\cb
-(-)^{r^1+\dots+r^n}b(r^1\tens\dots\tens r^n)\theta\Delta_\cb \\
&\qquad -\sum_{q+j+t=n}^{j<n} [(r^1\tens\dots\tens r^n)
(1^{\tens q}\tens B_j\tens1^{\tens t})]\theta\Delta_\cb \\
&= (r^1\tens\dots\tens r^n)\theta\Delta_\cb(1\tens b + b\tens1) \\
&\qquad -(-)^{r^1+\dots+r^n}b\Delta_\ca\sum_{k=0}^n
(r^1\tens\dots\tens r^k)\theta\tens(r^{k+1}\tens\dots\tens r^n)\theta \\
&\qquad -\sum_{q+j+t=n}^{j<n} \Delta_\ca
\{[(r^1\tens\dots\tens r^n)(1^{\tens q}\tens B_j\tens1^{\tens t})]
\Delta_{A_\infty(\ca,\cb)}(\theta\tens\theta)\} \\
&= \Delta_\ca\Bigl\{ \sum_{k=0}^n
(r^1\tens\dots\tens r^k)\theta\tens(r^{k+1}\tens\dots\tens r^n)\theta
(1\tens b + b\tens1) \\
&\qquad -(-)^{r^1+\dots+r^n} (1\tens b + b\tens1) \sum_{k=0}^n
(r^1\tens\dots\tens r^k)\theta\tens(r^{k+1}\tens\dots\tens r^n)\theta \\
&\qquad -\sum_{k+v+j+t=n}^{j<n} [(r^1\tens\dots\tens r^n)
(1^{\tens k}\bigotimes1^{\tens v}\tens B_j\tens1^{\tens t})]
(\theta\bigotimes\theta) \\
&\qquad -\sum_{q+j+w+u=n}^{j<n} [(r^1\tens\dots\tens r^n)
(1^{\tens q}\tens B_j\tens1^{\tens w}\bigotimes1^{\tens u})]
(\theta\bigotimes\theta) \Bigr\} \\
&= \Delta_\ca\Bigl\{ \sum_{k=0}^n
(r^1\tens\dots\tens r^k)\theta\tens(r^{k+1}\tens\dots\tens r^n)\theta b \\
&\qquad -\sum_{k=0}^n (-)^{r^{k+1}+\dots+r^n}
(r^1\tens\dots\tens r^k)\theta\tens b(r^{k+1}\tens\dots\tens r^n)\theta \\
&\qquad -\sum_{k+v+j+t=n}^{j<n} (r^1\tens\dots\tens r^k)\theta\tens
(r^{k+1}\tens\dots\tens r^n)(1^{\tens v}\tens B_j\tens1^{\tens t})\theta \\
&\qquad + \sum_{k=0}^n (-)^{r^{k+1}+\dots+r^n}
(r^1\tens\dots\tens r^k)\theta b\tens(r^{k+1}\tens\dots\tens r^n)\theta \\
&\qquad -\sum_{k=0}^n (-)^{r^1+\dots+r^n}
b(r^1\tens\dots\tens r^k)\theta\tens(r^{k+1}\tens\dots\tens r^n)\theta \\
&\qquad -\sum_{k=0}^n (-)^{r^{k+1}+\dots+r^n} \sum_{q+j+w=k}^{j<n}
(r^1\tens\dots\tens r^k)(1^{\tens q}\tens B_j\tens1^{\tens w})\theta
\tens(r^{k+1}\tens\dots\tens r^n)\theta \Bigr\} \\
&= \Delta_\ca\bigl[ f^0\tens(r^1\tens\dots\tens r^n)B_n
+ (r^1\tens\dots\tens r^n)B_n\tens f^n \bigr].
\end{align*}
The last three sums cancel out for all $k<n$ due to
\eqref{eq-Bn-bb-sum}, and for $k=n$ they give
$(r^1\tens\dots\tens r^n)B_n\tens f^n$ due to the same equation.
Similarly for the previous three sums. Therefore, \eqref{eq-Bn-bb-sum}
is, indeed, a recursive definition of components $B_n$ of a
coderivation $B$. The uniqueness of $B$ is obvious.

Clearly, $B^2:TsA_\infty(\ca,\cb)\to TsA_\infty(\ca,\cb)$ is a
(1,1)-coderivation of degree 2. From \eqref{eq-theta-b-B-theta} we find
\begin{multline*}
[(r^1\tens\dots\tens r^n)B^2]\theta =
[(r^1\tens\dots\tens r^n)B]\theta b
-(-)^{r^1+\dots+r^n+1}b[(r^1\tens\dots\tens r^n)B]\theta \\
= (r^1\tens\dots\tens r^n)\theta b^2
-(-)^{r^1+\dots+r^n}b[(r^1\tens\dots\tens r^n)\theta]b \\
-(-)^{r^1+\dots+r^n+1}b[(r^1\tens\dots\tens r^n)\theta]b
-b^2(r^1\tens\dots\tens r^n)\theta = 0.
\end{multline*}
Composing this equation with $\pr_1:Ts\cb\to s\cb$ we get
\[ 0 = [(r^1\tens\dots\tens r^n)B^2]\theta\pr_1 =
(r^1\tens\dots\tens r^n)[B^2]_n\pr_1. \]
Therefore, all components of the $(f^0,f^n)$-coderivation
$(r^1\tens\dots\tens r^n)[B^2]_n$ vanish. We deduce that the
coderivations $(r^1\tens\dots\tens r^n)[B^2]_n$ vanish, hence, all
$[B^2]_n=0$. Finally, $B^2=0$.
\end{proof}

Let us find explicitly the components of $B$, composing
\eqref{eq-theta-b-B-theta} with $\pr_1:Ts\cb\to s\cb$:
\begin{gather*}
B_1:(f,g)\to(f,g), \quad r\mapsto (r)B_1 = [r,b] = rb-(-)^rbr, \\
B_n:(f^0,f^1)\tens\dots\tens(f^{n-1},f^n)\to(f^0,f^n), \;
r^1\tens\dots\tens r^n\mapsto (r^1\tens\dots\tens r^n)B_n,
\text{ for } n>1,
\end{gather*}
where the last transformation is defined by its composition with
$\pr_1$:
\[ [(r^1\tens\dots\tens r^n)B_n]\pr_1
= [(r^1\tens\dots\tens r^n)\theta]b\pr_1. \]
In the other terms, for $n>1$
\begin{equation}
[(r^1\tens\dots\tens r^n)B_n]_k
= \sum_l (r^1\tens\dots\tens r^n)\theta_{kl}b_l.
\label{eq-Bn-components}
\end{equation}
Since $B^2=0$, we have, in particular,
\[ \sum_{r+n+t=k} (1^{\tens r}\tens B_n\tens1^{\tens t})B_{r+1+t} = 0 :
T^ksA_\infty(\ca,\cb) \to sA_\infty(\ca,\cb). \]

\begin{examples}
1) When $n=1$, $r:f\to g:\ca\to\cb$, we find the components of the
\ainf-transformation $(r)B_1:f\to g:\ca\to\cb$ as follows (see
\exasref{exas-r-T1-T2}):
\begin{align*}
[(r)B_1]_0 &= r_0b_1, \\
[(r)B_1]_1 &=
r_1b_1 + (f_1\tens r_0)b_2 + (r_0\tens g_1)b_2 -(-)^rb_1r_1, \\
[(r)B_1]_2 &= r_2b_1 + (f_2\tens r_0)b_2 + (f_1\tens r_1)b_2 +
(r_1\tens g_1)b_2 + (r_0\tens g_2)b_2 + (f_1\tens f_1\tens r_0)b_3 \\
&\hspace*{-2mm} +(f_1\tens r_0\tens g_1)b_3 +(r_0\tens g_1\tens g_1)b_3
-(-)^rb_2r_1 -(-)^r(1\tens b_1)r_2 -(-)^r(b_1\tens1)r_2.
\end{align*}

2) When $n=2$, $f \rTTo^r g \rTTo^p h:\ca\to\cb$, we find the
components of the \ainf-transformation $(r\tens p)B_2:f\to h:\ca\to\cb$
as follows:
\begin{align*}
[(r\tens p)B_2]_0 &= (r_0\tens p_0)b_2, \\
[(r\tens p)B_2]_1 &= (r_1\tens p_0)b_2 + (r_0\tens p_1)b_2 \\
&\quad+ (r_0\tens p_0\tens h_1)b_3 + (r_0\tens g_1\tens p_0)b_3 +
(f_1\tens r_0\tens p_0)b_3, \\
[(r\tens p)B_2]_2 &= (r_2\tens p_0)b_2 + (r_1\tens p_1)b_2
+ (r_0\tens p_2)b_2 + (r_1\tens p_0\tens h_1)b_3
+ (r_0\tens p_1\tens h_1)b_3 \\
&\quad+ (r_1\tens g_1\tens p_0)b_3 + (r_0\tens g_1\tens p_1)b_3
+ (f_1\tens r_1\tens p_0)b_3 + (f_1\tens r_0\tens p_1)b_3 \\
&\quad+ (r_0\tens p_0\tens h_2)b_3 + (r_0\tens g_2\tens p_0)b_3
+ (f_2\tens r_0\tens p_0)b_3 \\
&\quad+ (r_0\tens p_0\tens h_1\tens h_1)b_4
+ (r_0\tens g_1\tens g_1\tens p_0)b_4
+ (r_0\tens g_1\tens p_0\tens h_1)b_4 \\
&\quad+ (f_1\tens r_0\tens p_0\tens h_1)b_4
+ (f_1\tens r_0\tens g_1\tens p_0)b_4
+ (f_1\tens f_1\tens r_0\tens p_0)b_4.
\end{align*}

3) When $n=3$, $f \rTTo^r g \rTTo^p h \rTTo^t k:\ca\to\cb$, we find the
components of the \ainf-transformation
$(r\tens p\tens t)B_3:f\to k:\ca\to\cb$ as follows:
\begin{align*}
[(r\tens p\tens t)B_3]_0 &= (r_0\tens p_0\tens t_0)b_3, \\
[(r\tens p\tens t)B_3]_1 &= (r_1\tens p_0\tens t_0)b_3
+ (r_0\tens p_1\tens t_0)b_3 + (r_0\tens p_0\tens t_1)b_3 \\
&\quad+ (r_0\tens p_0\tens t_0\tens k_1)b_4
+ (r_0\tens p_0\tens h_1\tens t_0)b_4 \\
&\quad+ (r_0\tens g_1\tens p_0\tens t_0)b_4
+ (f_1\tens r_0\tens p_0\tens t_0)b_4.
\end{align*}
\end{examples}

\subsection{Differentials.}\label{sec-Differentials}
Let $\ca$, $\cc^1$, $\cc^2$, \dots, $\cc^q$, $\cb$ be \ainf-categories.
Let
$\phi:Ts\ca\tens Ts\cc^1\tens Ts\cc^2\tens\dots\tens Ts\cc^q\to Ts\cb$,
$(a\tens c^1\tens c^2\tens\dots\tens c^q)\mapsto
a.(c^1\tens c^2\tens\dots\tens c^q)\chi$
be a cocategory homomorphism of degree 0. If the homomorphism $\phi$
commutes with the differential:
\[ \phi b =
\bigl(\sum_{r+t=q}1^{\tens r}\tens b\tens1^{\tens t} \bigr)\phi,
\]
then $\phi$ is called an \ainf-functor (in a generalized sense,
extending \defref{def-ainf-functor}). This condition is fulfilled if
and only if $\chi$ commutes with the differential:
\begin{multline}
(c^1\tens c^2\tens\dots\tens c^q)\chi b \\
= \sum_{k=1}^q (-)^{c^{k+1}+\dots+c^q}
(c^1\tens\dots\tens c^kb\tens\dots\tens c^q)\chi
+ (-)^{c^1+\dots+c^q} b(c^1\tens c^2\tens\dots\tens c^q)\chi.
\label{eq-chi-b-1b1-b-chi}
\end{multline}
In particular, for $q=1$ we get the equation
\[ (c)\chi b = (cb)\chi + (-)^c b(c)\chi. \]

\begin{corollary}[to \propref{pro-B-F-KS-LH}]\label{cor-ainf-Ainfty-AB}
There is a unique \ainf-category structure for $A_\infty(\ca,\cb)$,
such that the action homomorphism
$\alpha:Ts\ca\tens TsA_\infty(\ca,\cb)\to Ts\cb$ is an \ainf-functor.
\end{corollary}

\begin{proof}
The homomorphism
$\alpha:a\tens r^1\tens\dots\tens r^n\mapsto
a[(r^1\tens\dots\tens r^n)\theta]$
of \corref{cor-action-alpha} uses $\chi=\theta$. Hence, $\alpha$ is an
\ainf-functor if and only if $(r)\theta b=(rB)\theta+(-)^rb(r)\theta$
for $r=r^1\tens\dots\tens r^n$, $n\ge0$, that is, if and only if
equations~\eqref{eq-theta-b-B-theta} hold.
\end{proof}

\begin{proposition}\label{prop-phi-C1-Cq-psi-ainf}
For any \ainf-functor
$\phi:Ts\ca\tens Ts\cc^1\tens Ts\cc^2\tens\dots\tens Ts\cc^q\to Ts\cb$
there is a unique \ainf-functor
$\psi:Ts\cc^1\tens Ts\cc^2\tens\dots\tens Ts\cc^q\to TsA_\infty(\ca,\cb)$,
such that
\begin{equation}
\phi = \bigl(Ts\ca\tens Ts\cc^1\tens Ts\cc^2\tens\dots\tens Ts\cc^q
\rTTo^{1\tens\psi} Ts\ca\tens TsA_\infty(\ca,\cb) \rTTo^\alpha
Ts\cb\bigr).
\label{eq-phi-1psi-alpha}
\end{equation}
\end{proposition}

\begin{proof}
If $f^i\in\Ob\cc^i$, then \eqref{eq-chi-b-1b1-b-chi} implies that the
cocategory homomorphism
$(f^1,f^2,\dots$, \linebreak[1]$f^q)\psi=(f^1,f^2,\dots,f^q)\chi$
of degree 0 commutes with the differential $b$. Hence, it is an
\ainf-functor, that is, an object of $A_\infty(\ca,\cb)$. By
\propref{prop-phi-C1-Cq-psi} there exists a unique cocategory
homomorphism
$\psi:Ts\cc^1\tens Ts\cc^2\tens\dots\tens Ts\cc^q\to TsA_\infty(\ca,\cb)$,
of degree 0, such that \eqref{eq-phi-1psi-alpha} holds. We have to
prove that $\psi$ commutes with the differential.

Using \eqref{eq-theta-b-B-theta} and \eqref{eq-chi-b-1b1-b-chi} we find
for $\chi=\psi\theta$ the equation
\begin{multline}
(c^1\tens c^2\tens\dots\tens c^q)\psi B\theta =
(c^1\tens c^2\tens\dots\tens c^q)\psi\theta b
-(-)^{c^1+\dots+c^q} b(c^1\tens c^2\tens\dots\tens c^q)\psi\theta \\
= (c^1\tens c^2\tens\dots\tens c^q)
\sum_{k=1}^q (1^{\tens k-1}\tens b\tens1^{\tens q-k})\psi\theta.
\label{eq-psi-B-theta-1b1-psi-theta}
\end{multline}
Notice that $\psi B$ and
$\kappa\overset{\text{def}}=
\sum_{k=1}^q (1^{\tens k-1}\tens b\tens1^{\tens q-k})\psi$
are both $(\psi,\psi)$-coderivations. Let $c^i$ be an element of
$T^{n^i}s\cc^i=s\cc^i(f^i,-)\tens\dots\tens s\cc^i(-,g^i)$ for
$1\le i\le q$. Composing \eqref{eq-psi-B-theta-1b1-psi-theta} with
$\pr_1:Ts\cb\to s\cb$ we get
\begin{multline*}
(c^1\tens c^2\tens\dots\tens c^q)[\psi B]_{n^1\dots n^q}\pr_1 =
(c^1\tens c^2\tens\dots\tens c^q)\psi B\theta\pr_1 \\
= (c^1\tens c^2\tens\dots\tens c^q)\kappa\theta\pr_1 =
(c^1\tens c^2\tens\dots\tens c^q)\kappa_{n^1\dots n^q}\pr_1.
\end{multline*}
Since all components of
$((f^1,f^2,\dots,f^q)\psi,(g^1,g^2,\dots,g^q)\psi)$-coderivations
$(c^1\tens c^2\tens\dots\tens c^q)[\psi B]_{n^1\dots n^q}$ and
$(c^1\tens c^2\tens\dots\tens c^q)\kappa_{n^1\dots n^q}$ coincide,
these coderivations coincide as well. Therefore, all the components of
$(\psi,\psi)$-coderivations $\psi B$ and $\kappa$ coincide. We conclude
that these coderivations coincide as well:
$\psi B=\kappa=\sum_{k=1}^q (1^{\tens k-1}\tens b\tens1^{\tens q-k})\psi$.
\end{proof}

\section{\texorpdfstring{Enriched category of $A_\infty$-categories}
{Enriched category of A8-categories}}
\begin{definition}[Differential graded cocategory]
A differential graded cocategory $\cc$ is a graded cocategory equipped
with a (1,1)-coderivation
$b=\bigl(b:\cc(X,Y)\to\cc(X,Y)\bigr)_{X,Y\in\Ob\cc}$ of degree 1, such
that $b^2=0$.
\end{definition}

As in \secref{def-Cocategory} a differential graded cocategory can be
identified with a differential graded $\kk$\n-coalgebra, decomposed in
a special way. An example of a differential graded cocategory is given
by $Ts\ca$, where $\ca$ is an \ainf-category.

Differential graded cocategories form a symmetric monoidal category
dgCoCat. If $\cc$ and $\cd$ are differential graded cocategories then
their tensor product is the graded cocategory $\cc\tens\cd$, equipped
with the differential $1\tens b+b\tens1$. We want to show that the
category of \ainf-categories is enriched in dgCoCat.

Let $\ca$, $\cb$, $\cc$ be \ainf-categories. There is a graded
cocategory morphism of degree 0
\[ M: TsA_\infty(\ca,\cb)\tens TsA_\infty(\cb,\cc)
\to TsA_\infty(\ca,\cc), \]
defined in \secref{sec-enr-cat-grad-centr-bimod} via
diagram~\eqref{dia-def-M}. Since all cocategory morphisms $\alpha$,
$\alpha\tens1$ in this diagram commute with the differential by
\corref{cor-ainf-Ainfty-AB}, the cocategory morphism $M$ also commutes
with the differential:
\begin{equation}
(1\tens B+B\tens1)M=MB
\label{eq-1B-B1-M-MB}
\end{equation}
by \propref{prop-phi-C1-Cq-psi-ainf}. Therefore, $M$ is an
\ainf-functor. The unit $\eta_\ca:T\1\to T\Coder(s\ca,s\ca)$,
$1\mapsto\id_\ca$ also is an \ainf-functor for trivial reasons. The set
\[ \text{dgCoCat}(T\1,TsA_\infty(\ca,\cb)) =
\text{Maps}(\{*\},\Ob A_\infty(\ca,\cb)) = \text{dgCoCat}(Ts\ca,Ts\cb)
\]
is the set of \ainf-functors $\ca\to\cb$. We summarize the above
statements as follows: the category of \ainf-categories is enriched in
dgCoCat. Moreover, it is enriched in the monoidal subcategory of
dgCoCat generated by $Ts\cc$, where $\cc$ are \ainf-categories.

Let us apply \propref{prop-phi-C1-Cq-psi-ainf} to the \ainf-functor
$M$. From that result we deduce the existence of a unique \ainf-functor
\[ A_\infty(\ca,\_): A_\infty(\cb,\cc) \to
A_\infty(A_\infty(\ca,\cb),A_\infty(\ca,\cc)), \]
such that
\begin{multline}
M = \bigl[ TsA_\infty(\ca,\cb)\tens TsA_\infty(\cb,\cc)
\rTTo^{1\tens A_\infty(\ca,\_)} \\
TsA_\infty(\ca,\cb)\tens TsA_\infty(A_\infty(\ca,\cb),A_\infty(\ca,\cc))
\rTTo^\alpha TsA_\infty(\ca,\cc)\bigr].
\label{eq-M-1A-alpha}
\end{multline}
Let us find the components of $A_\infty(\ca,\_)$.

\begin{proposition}\label{pro-Ainf-A--strict}
The \ainf-functor $A_\infty(\ca,\_)$ is strict. It maps an object of
$A_\infty(\cb,\cc)$, an \ainf-functor $g:\cb\to\cc$, to the object of
the target \ainf-category
$(1\tens g)M:A_\infty(\ca,\cb)\to A_\infty(\ca,\cc)$ (which is also an
\ainf-functor). The first component $A_\infty(\ca,\_)_1$ maps an
element $t$ of $sA_\infty(\cb,\cc)(g,h)$, a $(g,h)$-coderivation
$t:Ts\cb\to Ts\cc$, to the $((1\tens g)M,(1\tens h)M)$-coderivation
$(1\tens t)M:TsA_\infty(\ca,\cb)\to TsA_\infty(\ca,\cc)$, an element of
$sA_\infty(A_\infty(\ca,\cb),A_\infty(\ca,\cc))
((1\tens g)M,(1\tens h)M)$.
\end{proposition}

\begin{proof}
Clearly, $A_\infty(\ca,\_)$ gives the mapping of objects
$g\mapsto(1\tens g)M$ as described. To prove that
$A_\infty(\ca,\_)_1:t\mapsto(1\tens t)M$ and $A_\infty(\ca,\_)_k=0$ for
$k>1$ it suffices to substitute a cocategory homomorphism with those
components into \eqref{eq-M-1A-alpha} and to check this equation (see
\propref{prop-phi-C1-Cq-psi}). Let
\begin{align}
p^1\tens\dots\tens p^n &\in sA_\infty(\ca,\cb)(f^0,f^1)\tens\dots\tens
sA_\infty(\ca,\cb)(f^{n-1},f^n), \notag \\
t^1\tens\dots\tens t^m &\in sA_\infty(\cb,\cc)(g^0,g^1)\tens\dots\tens
sA_\infty(\cb,\cc)(g^{m-1},g^m).
\label{eq-p1-pn-t1-tm-in}
\end{align}
The equation to check is
\[ (p^1\tens\dots\tens p^n\tens t^1\tens\dots\tens t^m)M =
(p^1\tens\dots\tens p^n).[(t^1\tens\dots\tens t^m)A_\infty(\ca,\_)]\theta,
\]
that is,
\[ (p^1\tens\dots\tens p^n\tens t^1\tens\dots\tens t^m)M =
(p^1\tens\dots\tens p^n).[(1\tens t^1)M\tens\dots\tens(1\tens t^m)M]\theta.
\]
The left hand side is a cocategory homomorphism. Let us prove that the
right hand side
\[ (p^1\tens\dots\tens p^n\tens t^1\tens\dots\tens t^m)L
\overset{\text{def}}=
(p^1\tens\dots\tens p^n).[(1\tens t^1)M\tens\dots\tens(1\tens t^m)M]\theta
\]
is also a cocategory homomorphism. Indeed,
\begin{align*}
& (p^1\tens\dots\tens p^n\tens t^1\tens\dots\tens t^m)L\Delta \\
&= (p^1\tens\dots\tens p^n).
[(1\tens t^1)M\tens\dots\tens(1\tens t^m)M]\theta\Delta \\
&= (p^1\tens\dots\tens p^n)\Delta \sum_{k=0}^m
[(1\tens t^1)M\tens\dots\tens(1\tens t^k)M]\theta \\
&\hspace*{52mm} \tens[(1\tens t^{k+1})M\tens\dots\tens(1\tens t^m)M]\theta \\
&= \sum_{l=0}^n \sum_{k=0}^m (-)^{(p^{l+1}+\dots+p^n)(t^1+\dots+t^k)}
(p^1\tens\dots\tens p^l).
[(1\tens t^1)M\tens\dots\tens(1\tens t^k)M]\theta \\
&\hspace*{52mm} \tens(p^{l+1}\tens\dots\tens p^n).
[(1\tens t^{k+1})M\tens\dots\tens(1\tens t^m)M]\theta \\
&= [(p^1\tens\dots\tens p^n)\Delta\tens(t^1\tens\dots\tens t^m)\Delta]
(1\tens c\tens1)(L\tens L)
\end{align*}
by \propref{pro-theta-D-D-theta-theta}.

Let us prove that the components of $M$ and $L$ coincide. For $m=0$ and
any $n>0$ we have
\[ (p^1\tens\dots\tens p^n\mid g^0)L_{n0}
= (p^1\tens\dots\tens p^n).(1\tens g^0)M\pr_1
= (p^1\tens\dots\tens p^n\mid g^0)M_{n0}, \]
hence, $L_{n0}=M_{n0}$. For $m=1$ and any $n\ge0$ we have
\[ (p^1\tens\dots\tens p^n\tens t^1)L_{n1}
= (p^1\tens\dots\tens p^n).(1\tens t^1)M\pr_1
= (p^1\tens\dots\tens p^n\tens t^1)M_{n1}, \]
hence, $L_{n1}=M_{n1}$. For $m>1$ and any $n\ge0$ we have $L_{nm}=0$
and $M_{nm}=0$. Therefore, $L=M$ and the proposition is proved.
\end{proof}

\begin{corollary}
For all $m>0$ and all $t^1\tens\dots\tens t^m$ as in
\eqref{eq-p1-pn-t1-tm-in} we have
\begin{equation}
[(1\tens t^1)M\tens\dots\tens(1\tens t^m)M]\tilde{B}_m
= [1\tens(t^1\tens\dots\tens t^m)B_m]M,
\label{eq-ttMB-ttBM}
\end{equation}
where $\tilde{B}$ denotes the differential in
$TsA_\infty(A_\infty(\ca,\cb),A_\infty(\ca,\cc))$.
\end{corollary}

Indeed, the general property of an \ainf-functor
$A_\infty(\ca,\_)\tilde{B}=BA_\infty(\ca,\_)$ reduces to the above
formula, since $A_\infty(\ca,\_)$ is strict.

In the following definition we introduce \ainf-analogs of natural
transformations.

\begin{definition}[$\omega$-globular set of $A_\infty$-categories]
\label{def-globular-set}
A \emph{natural \ainf-transformation} $r:f\to g:\ca\to\cb$ (natural
transformation in terms of \cite{Fukaya:FloerMirror-II}) is an
\ainf-transformation of degree $-1$ such that $rb+br=0$ (that is,
$(r)B_1=0$). The $\omega$\n-globular set \cite{Bat-weak} $A_\omega$ of
\ainf-categories is defined as follows: objects (0\n-morphisms) are
\ainf-categories $\ca$; 1\n-morphisms are \ainf-functors $f:\ca\to\cb$;
2\n-morphisms are natural \ainf-transformations $r:f\to g:\ca\to\cb$;
3\n-morphisms $\lambda:r\to s:f\to g:\ca\to\cb$ are
$(f,g)$-coderivations of degree $-2$, such that $r-s=[\lambda,b]$; for
$n\ge3$ $n$\n-morphisms
$\lambda_n:\lambda_{n-1}\to\mu_{n-1}:\dots:r\to s:f\to g:\ca\to\cb$ are
$(f,g)$-coderivations of degree $1-n$, such that
$\lambda_{n-1}-\mu_{n-1}=[\lambda_n,b]$ (notice that the both sides are
$(f,g)$-coderivations of degree $2-n$).
\end{definition}

\begin{remark}\label{rem-AinfA-Aomega-Aomega}
Let us notice that the \ainf-functor $A_\infty(\ca,\_)$ from
\propref{pro-Ainf-A--strict} defines a map of the $\omega$\n-globular
set $A_\omega$ into itself. Indeed, objects $\cb$ of $A_\omega$ are
mapped into objects $A_\infty(\ca,\cb)$, 1\n-morphisms
$g:\cb\to\cc$ are mapped into 1\n-morphisms
$(1\tens g)M:A_\infty(\ca,\cb)\to A_\infty(\ca,\cc)$ and the first
component $A_\infty(\ca,\_)_1$ maps $(g,h)$-coderivations into
$((1\tens g)M,(1\tens h)M)$-coderivations. Moreover, if the equation
$\lambda_{n-1}-\mu_{n-1}=\lambda_nB_1$ holds for $(g,h)$-coderivations,
then
$(1\tens\lambda_{n-1})M-(1\tens\mu_{n-1})M=(1\tens\lambda_n)M\tilde{B}_1$
by \eqref{eq-ttMB-ttBM}, so the sources and the targets are preserved.
\end{remark}

It might be useful to turn the $\omega$\n-globular set $A_\omega$ into
a weak non-unital $\omega$\n-category in the sense of some of the
existing definitions of the latter. Plenty of such definitions
including \cite{Bat-weak} are listed in Leinster's survey
\cite{Leinster:DefSurvey}. We do not try to proceed in this direction.
Instead we truncate the $\omega$\n-globular set to a 2\n-globular set
(that is, we deal with 0-, 1- and 2\n-morphisms) and we make a
2\n-category out of it.

\section{\texorpdfstring{2-categories of $A_\infty$-categories}
 {2-categories of A8-categories}}
Let $\ck$ denote the category $\Kht(\kk\modul)=H^0(\Com(\kk\modul))$ of
differential graded complexes of $\kk$\n-modules, whose morphisms are
chain maps modulo homotopy. Equipped with the usual tensor product, the
unit object $\kk$ and the standard symmetry, $\ck$ becomes a
$\kk$\n-linear closed monoidal symmetric category. The inner hom-object
is the usual $\Hom_\kk^\bull(\text-,\text-)$. There is a notion of a
category $\cc$ enriched in $\ck$ ($\ck$\n-categories, $\ck$\n-functors,
$\ck$\n-natural transformations), see Kelly \cite{KellyGM:bascec}: for all
objects $X$, $Y$ of $\cc$ $\cc(X,Y)$ is an object of $\ck$. There is a
similar notion of a 2\n-category enriched in $\ck$, or a
$\ck$-2-category: it consist of a class of objects $\Ob\cc$, a class of
1\n-morphisms $\cc(X,Y)$ for each pair of objects $X$, $Y$ of $\ca$, an
object of 2\n-morphisms $\cc(X,Y)(f,g)\in\Ob\ck$ for each pair of
1\n-morphisms $f,g\in\cc(X,Y)$ and other data. We shall consider
1\n-unital, non-2-unital $\ck$-2-categories. They are equipped with the
following operations: associative composition of 1\n-morphisms,
commuting left and right associative actions of 1\n-morphisms on
2\n-morphisms (these actions are morphisms in $\ck$), 1\n-units,
associative vertical composition of 2\n-morphisms (a morphism in $\ck$)
compatible with the left and right actions of 1\n-morphisms on
2\n-morphisms and such that the both ways to obtain the horizontal
composition coincide. Precise definitions are given in
\appref{ap-sec-enrich}.

\begin{proposition}\label{pro-1uni-n2uni-KAinf}
The following data define a 1\n-unital, non-2-unital $\ck$-2-category
$\ck A_\infty$:
\begin{events}
\item[\sbull] objects are \ainf-categories;

\item[\sbull] 1\n-morphisms are \ainf-functors;

\item[\sbull] an object of 2\n-morphisms between $f,g:\ca\to\cb$ is
$(A_\infty(\ca,\cb)(f,g),m_1)\in\Ob\ck$, $m_1=sB_1s^{-1}$;

\item[\sbull] the composition of 1\n-morphisms is the composition of
\ainf-functors;

\item[\sbull] unit 1\n-morphisms are identity \ainf-functors;

\item[\sbull] the right action of a 1\n-morphism $k:\cb\to\cc$ on
2\n-morphisms is the chain map
$(A_\infty(\ca,\cb)(f,g),m_1)\to(A_\infty(\ca,\cc)(fk,gk),m_1)$,
$rs^{-1}\mapsto(rs^{-1})\cdot k=(rk)s^{-1}$, where $r$ is an
$(f,g)$-coderivation;

\item[\sbull] the left action of a 1\n-morphism $e:\cd\to\ca$ on
2\n-morphisms is the chain map
$(A_\infty(\ca,\cb)(f,g),m_1)\to(A_\infty(\cd,\cb)(ef,eg),m_1)$,
$rs^{-1}\mapsto e\cdot(rs^{-1})=(er)s^{-1}$, where $r$ is an
$(f,g)$-coderivation;

\item[\sbull] the vertical composition is the chain map
$m_2=(s\tens s)B_2s^{-1}:A_\infty(\ca,\cb)(f,g)\tens
A_\infty(\ca,\cb)(g,h)\to A_\infty(\ca,\cb)(f,h)$.
\end{events}
\end{proposition}

\begin{proof}
Clearly, the composition of 1\n-morphisms and the actions of
1\n-morphisms on 2\n-morphisms are associative. The right and the left
actions are unital, and commute with each other. The equation
$-(1\tens m_1+m_1\tens1)m_2+m_2m_1=0$ (see \eqref{eq-m-sbs-1} for
$k=2$) shows that $m_2$ is a chain map. The equation
\begin{equation}
m_3m_1 + (1\tens1\tens m_1+1\tens m_1\tens1+m_1\tens1\tens1)m_3
- (m_2\tens1)m_2 + (1\tens m_2)m_2 = 0
\label{eq-m2-associative-K}
\end{equation}
(see \eqref{eq-m-sbs-1} for $k=3$) shows that $m_2$ is associative in
$\ck$.

Let us check that the vertical composition is compatible with the
actions of 1\n-morphisms on 2\n-morphisms. Applying
equation~\eqref{eq-1B-B1-M-MB} to
$r\tens p\tens1\in sA_\infty(\ca,\cb)(f,g)\tens
sA_\infty(\ca,\cb)(g,h)\tens\kk \subset
T^2sA_\infty(\ca,\cb)(f,h)\tens T^0sA_\infty(\cb,\cc)(k,k)$
we find that
\begin{equation}
(r\tens p\mid k)M_{20}B_1+(rk\tens pk)B_2 =
[(r\tens p)(1\tens B_1+B_1\tens1)\mid k]M_{20} + (r\tens p)B_2k.
\label{eq-rkpkB2-M20-rpB2k}
\end{equation}
One deduces that
$(rs^{-1}\cdot k\tens ps^{-1}\cdot k)m_2=(rs^{-1}\tens ps^{-1})m_2\cdot k$
in $\ck$. Applying equation~\eqref{eq-1B-B1-M-MB} to
$1\tens r\tens p\in\kk\tens sA_\infty(\ca,\cb)(f,g)\tens
sA_\infty(\ca,\cb)(g,h)\subset
T^0sA_\infty(\cd,\ca)(e,e)\tens T^2sA_\infty(\ca,\cb)(f,h)$
we find that
\begin{multline}
(er\tens ep)B_2 = (e\mid r\tens p)M_{02}B_1+(er\tens ep)B_2 \\
= [e\mid(r\tens p)(1\tens B_1+B_1\tens1)]M_{02} + e(r\tens p)B_2
= e(r\tens p)B_2
\label{eq-krkpB2-M02-krpB2}
\end{multline}
(notice that $M_{02}=0$). Therefore,
$(e\cdot rs^{-1}\tens e\cdot ps^{-1})m_2=e\cdot(rs^{-1}\tens ps^{-1})m_2$.

Now let us prove distributivity. Applying
equation~\eqref{eq-1B-B1-M-MB} to
$r\tens p\in sA_\infty(\ca,\cb)(f,g)\tens sA_\infty(\cb,\cc)(h,k)$ we
find that
\[ (rh\tens gp)B_2 + (-)^{rp}(fp\tens rk)B_2 + (r\tens p)M_{11}B_1
= (r\tens p)(1\tens B_1+B_1\tens1)M_{11}. \]
Thus, $(rh\tens gp)B_2+(-)^{rp}(fp\tens rk)B_2=0$ in $\ck$. We deduce
that modulo homotopy
\begin{multline*}
(rs^{-1}\cdot h\tens g\cdot ps^{-1})m_2s =
(rhs^{-1}\tens gps^{-1})(s\tens s)B_2 = (-)^{p+1}(rh\tens gp)B_2 \\
= (-)^{rp+p}(fp\tens rk)B_2 = (-)^{rp+p}(fps^{-1}s\tens rks^{-1}s)B_2 =
(-)^{rp+p+r+1}(f\cdot ps^{-1}\tens rs^{-1}\cdot k)m_2s.
\end{multline*}
Therefore,
$(rs^{-1}\cdot h\tens g\cdot ps^{-1})m_2
=(-)^{(r+1)(p+1)}(f\cdot ps^{-1}\tens rs^{-1}\cdot k)m_2$
in $\ck$, as stated.
\end{proof}

The 0\n-th cohomology functor $H^0=\ck(\kk,\_):\ck\to\kk\modul$,
$X\mapsto H^0(X)=\ck(\kk,X)$ is lax monoidal symmetric, since the
complex $\kk$ concentrated in degree 0 is the unit object of $\ck$. It
determines a functor $H^0:\KCatb\to\kk\text-{\mathcal C}at$. To a
$\ck$\n-category $\cc$ it assigns a $\kk$\n-linear category $H^0(\cc)$
with the same class of objects, and for each pair $X$, $Y$ of objects
of $\cc$ the $\kk$\n-module $H^0(\cc)(X,Y)=H^0(\cc(X,Y))$. The functor
$H^0:\KCatb\to\kk\text-{\mathcal C}at$ is also lax monoidal symmetric.
Therefore, there is a functor $\KCatCat\to\CatCat$, again denoted $H^0$
by abuse of notation, and the corresponding functor
$\ck\text-2\text-\Cat^{nu}\to2\text-\Cat^{nu}$. See
\appref{ap-sec-enrich} for the definition of 1\n-unital, non-2-unital
$\ck$- or $\kk$- 2\n-categories. To $\ck A_\infty$ the functor assigns
a 1\n-unital, non-2-unital $\kk$\n-linear 2\n-category $A_\infty$. Let
us describe it in detail. Objects of $A_\infty$ are \ainf-categories,
1\n-morphisms are \ainf-functors, and 2\n-morphisms are elements of
\[ H^0(A_\infty(\ca,\cb)(f,g),m_1) \rTTo^s_\sim
H^{-1}(sA_\infty(\ca,\cb)(f,g),B_1),
\]
that is, equivalence classes of natural \ainf-transformations
$r:f\to g:\ca\to\cb$. Natural \ainf-transformations
$r,t:f\to g:\ca\to\cb$ are \emph{equivalent}, if they are connected by
a 3\n-morphism $\lambda:r\to t$, that is, $r-t=\lambda B_1$. Both
compositions of 1\n-morphisms with 2\n-morphisms
$\Mor_2(\ca,\cb)\times\Mor_1(\cb,\cc)\to\Mor_2(\ca,\cc)$,
$(r,h)\mapsto rh$ and
$\Mor_1(\ca,\cb)\times\Mor_2(\cb,\cc)\to\Mor_2(\ca,\cc)$,
$(f,p)\mapsto fp$ are compositions of $\kk$\n-linear maps
$Ts\ca\to Ts\cb\to Ts\cc$. The vertical composition $m_2$ of
2\n-morphisms, translated to $H^{-1}(sA_\infty)$ assigns the natural
\ainf-transformation $r\cdot p=(r\tens p)B_2$ to natural
\ainf-transformations $r:f\to g$ and $p:g\to h$. Indeed,
$(rs^{-1}\tens ps^{-1})m_2s
=(rs^{-1}\tens ps^{-1})(s\tens s)B_2=(r\tens p)B_2$.
Compatibility of these constructions with the equivalence relation is
obvious from the construction, and can be verified directly.

\begin{remark}\label{rem-AA-A-A-2-fun}
Let $\ca$ be an \ainf-category. It determines a map
$A_\omega\to A_\omega$, described in \remref{rem-AinfA-Aomega-Aomega}.
This map restricts to a map $A_\infty(\ca,\_):A_\infty\to A_\infty$. It
takes an \ainf-category $\cb$ to the \ainf-category
$A_\infty(\ca,\cb)$, an \ainf-functor $g:\cb\to\cc$ to the
\ainf-functor $(1\tens g)M:A_\infty(\ca,\cb)\to A_\infty(\ca,\cc)$, and
an equivalence class of a natural \ainf-transformation
$t:g\to h:\ca\to\cb$ to the equivalence class of the natural
\ainf-transformation
$(1\tens t)M:(1\tens g)M\to(1\tens h)M:
A_\infty(\ca,\cb)\to A_\infty(\ca,\cc)$,
see \remref{rem-AinfA-Aomega-Aomega}. Let us prove that
$A_\infty(\ca,\_):A_\infty\to A_\infty$ is a strict 2\n-functor.
Indeed, it preserves the composition of 1\n-morphisms,
$(1\tens f)M(1\tens g)M=(1\tens fg)M$, and the both compositions of
1\n-morphisms and 2\n-morphisms, $(1\tens f)M(1\tens t)M=(1\tens ft)M$,
$(1\tens t)M(1\tens f)M=(1\tens tf)M$, due to associativity of $M$. The
vertical composition of 2\n-morphisms is preserved due to
\eqref{eq-ttMB-ttBM} for $m=2$:
\[ [(1\tens t^1)M\tens(1\tens t^2)M]\tilde{B}_2
= [1\tens(t^1\tens t^2)B_2]M. \]
\end{remark}

\begin{definition}[Unital $A_\infty$-categories]\label{def-unital-cat}
Let $\cc$ be an \ainf-category. It is called \emph{unital} if for each
object $X$ of $\cc$ there is a \emph{unit element} -- a $\kk$\n-linear
map $\sS{_X}\uni^\cc_0:\kk\to(s\cc)^{-1}(X,X)$ such that
$\sS{_X}\uni^\cc_0b_1=0$,
$(\sS{_X}\uni^\cc_0\tens\sS{_X}\uni^\cc_0)b_2-\sS{_X}\uni^\cc_0\in\im
b_1$,
and for all pairs $X$, $Y$ of objects of $\cc$ the chain maps
$(1\tens\sS{_Y}\uni^\cc_0)b_2$,
$(\sS{_X}\uni^\cc_0\tens1)b_2:s\cc(X,Y)\to s\cc(X,Y)$ are homotopy
invertible.
\end{definition}

In particular, an \ainf-algebra $\cc$ is unital if it has an element
$\uni^\cc_0\in(s\cc)^{-1}$ such that $\uni^\cc_0b_1=0$,
$(\uni^\cc_0\tens\uni^\cc_0)b_2-\uni^\cc_0\in\im b_1$, and the chain
maps $(1\tens\uni^\cc_0)b_2$, $(\uni^\cc_0\tens1)b_2:s\cc\to s\cc$ are
homotopy invertible. Our definition differs from that of a homological
unit (e.g. \cite{Lefevre-Ainfty-these}) by the last invertibility
condition. It produces rather a homotopical unit:

\begin{lemma}\label{lem-1ib1-i1b-1}
Let $\sS{_X}\uni^\cc_0$ be as in \defref{def-unital-cat} of a unital
\ainf-category $\cc$, then for each pair $X$, $Y$ of objects of $\cc$
we have
\begin{align*}
(1\tens\sS{_Y}\uni^\cc_0)b_2 \sim 1 &: s\cc(X,Y) \to s\cc(X,Y), \\
(\sS{_X}\uni^\cc_0\tens1)b_2 \sim-1 &: s\cc(X,Y) \to s\cc(X,Y).
\end{align*}
\end{lemma}

\begin{proof}
For each object $X$ of $\cc$ there is a $\kk$\n-linear map
$\sS{_X}v_0:\kk\to(s\cc)^{-2}(X,X)$ such that
$(\sS{_X}\uni^\cc_0\tens\sS{_X}\uni^\cc_0)b_2-\sS{_X}\uni^\cc_0
=\sS{_X}v_0b_1$.
Hence,
\begin{align*}
& (\sS{_X}\uni^\cc_0\tens1)b_2(\sS{_X}\uni^\cc_0\tens1)b_2
= (\sS{_X}\uni^\cc_0\tens\sS{_X}\uni^\cc_0\tens1)(1\tens b_2)b_2
\overset{\eqref{eq-b-b-0}}{\underset{k=3}\sim}
- (\sS{_X}\uni^\cc_0\tens\sS{_X}\uni^\cc_0\tens1)(b_2\tens1)b_2 \\
&= - [(\sS{_X}\uni^\cc_0\tens\sS{_X}\uni^\cc_0)b_2\tens1]b_2
= -(\sS{_X}\uni^\cc_0\tens1)b_2 -(\sS{_X}v_0b_1\tens1)b_2 \\
&= -(\sS{_X}\uni^\cc_0\tens1)b_2 + b_1(\sS{_X}v_0\tens1)b_2
+ (\sS{_X}v_0\tens1)b_2b_1
\sim -(\sS{_X}\uni^\cc_0\tens1)b_2,
\end{align*}
\begin{align*}
& (1\tens\sS{_Y}\uni^\cc_0)b_2(1\tens\sS{_Y}\uni^\cc_0)b_2
= - (1\tens\sS{_Y}\uni^\cc_0\tens\sS{_Y}\uni^\cc_0)(b_2\tens1)b_2
\overset{\eqref{eq-b-b-0}}{\underset{k=3}\sim}
(1\tens\sS{_Y}\uni^\cc_0\tens\sS{_Y}\uni^\cc_0)(1\tens b_2)b_2 \\
&= [1\tens(\sS{_Y}\uni^\cc_0\tens\sS{_Y}\uni^\cc_0)b_2]b_2
= (1\tens\sS{_Y}\uni^\cc_0)b_2 + (1\tens\sS{_Y}v_0b_1)b_2 \\
&= (1\tens\sS{_Y}\uni^\cc_0)b_2 - b_1(1\tens\sS{_Y}v_0)b_2
- (1\tens\sS{_Y}v_0)b_2b_1
\sim (1\tens\sS{_Y}\uni^\cc_0)b_2.
\end{align*}
We see that $-(\sS{_X}\uni^\cc_0\tens1)b_2$ and
$(1\tens\sS{_Y}\uni^\cc_0)b_2$ are invertible idempotents in $\ck$.
Therefore, these maps are both homotopic to the identity map.
\end{proof}

This lemma shows that a unital \ainf-algebra may be defined as an
\ainf-algebra $\cc$, such that the graded associative $\kk$\n-algebra
$H^\bull(\cc,m_1)$ has a unit $1\in H^0(\cc,m_1)$ and for some/any
representative $1^\cc\in\cc^0$ of the class $1\in H^0(\cc,m_1)$ the
chain maps $(\id\tens1^\cc)m_2$, $(1^\cc\tens\id)m_2:\cc\to\cc$ are
homotopic to $\id_\cc$. A unit element $\uni^\cc_0\in(s\cc)^{-1}$
corresponds to a unit $1^\cc\in\cc^0$ via $1^\cc s=\uni^\cc_0$.

\begin{proposition}\label{pro-unital-unit}
Let $\cc$ be a unital \ainf-category. Then the collection
$\sS{_X}\uni^\cc_0$ extends to a natural \ainf-transformation
$\uni^\cc:\id_\cc\to\id_\cc:\cc\to\cc$ such that
$(\uni^\cc\tens\uni^\cc)B_2\equiv\uni^\cc$.
\end{proposition}

\begin{proof}
Let $\kk$\n-linear maps $\sS{_X}v_0:\kk\to(s\cc)^{-2}(X,X)$ satisfy the
equations
$(\sS{_X}\uni^\cc_0\tens\sS{_X}\uni^\cc_0)b_2-\sS{_X}\uni^\cc_0
=\sS{_X}v_0b_1$.
We will prove that given $\sS{_X}\uni^\cc_0$, $\sS{_X}v_0$ (with
$\sS{_X}\uni^\cc_0b_1=0$) are 0\n-th components of a natural
\ainf-transformation $\uni^\cc$ and a 3\n-morphism $v$ as follows:
\begin{gather*}
\uni^\cc: \id_\cc\to\id_\cc: \cc\to\cc, \\
v: (\uni^\cc\tens\uni^\cc)B_2 \to \uni^\cc:\id_\cc\to\id_\cc:\cc\to\cc.
\end{gather*}
That is, we will prove the existence of $(1,1)$-coderivations
$\uni^\cc,v:Ts\cc\to Ts\cc$ of degree $-1$ (resp. $-2$) such that
\begin{align*}
\uni^\cc b + b\uni^\cc &= 0, \\
(\uni^\cc\tens\uni^\cc)B_2 - \uni^\cc &= vb - bv.
\end{align*}
We already have the 0\n-th components $\uni^\cc_0$, $v_0$. Let us
construct the other components of $\uni^\cc$ and $v$ by induction.
Given a positive $n$, assume that we have already found components
$\uni^\cc_m$, $v_m$ of the sought $\uni^\cc$, $v$ for $m<n$, such
that the equations
\begin{equation}
(\uni^\cc b)_m + (b\uni^\cc)_m = 0:
s\cc(X_0,X_1)\tens\dots\tens s\cc(X_{m-1},X_m) \to s\cc(X_0,X_m),
\label{eq-ibm-bim-0}
\end{equation}
\begin{multline}
[(\uni^\cc\tens\uni^\cc)B_2]_m - \uni^\cc_m = (vb)_m - (bv)_m: \\
s\cc(X_0,X_1)\tens\dots\tens s\cc(X_{m-1},X_m) \to s\cc(X_0,X_m)
\label{eq-iiB2m-im-vbm-bvm}
\end{multline}
are satisfied for all $m<n$. Here
\( (f)_m = \bigl(T^ms\cc \rMono Ts\cc \rTTo^f Ts\cc
\rTTo^{\pr_1} s\cc \bigr) \)
for an arbitrary morphism of quivers $f:Ts\cc\to Ts\cc$.

Introduce $(1,1)$-coderivations $\tilde{\uni},v:Ts\cc\to Ts\cc$ of
degree $-1$ and $-2$ by their components
$(\uni^\cc_0,\uni^\cc_1,\dots$,\linebreak[1]$\uni^\cc_{n-1},0,0,\dots)$
(resp. $(v_0,v_1,\dots,v_{n-1},0,0,\dots)$). Define
$(1,1)$-coderivations $\lambda=\tilde{\uni}b+b\tilde{\uni}$ of degree
$0$ and
$\nu=(\tilde{\uni}\tens\tilde{\uni})B_2-\tilde{\uni}-\tilde{v}B_1$ of
degree $-1$. Then equations \eqref{eq-ibm-bim-0},
\eqref{eq-iiB2m-im-vbm-bvm} imply that $\lambda_m=0$, $\nu_m=0$ for
$m<n$. The identity $\lambda b-b\lambda=0$ implies that
\[ \lambda_nd = \lambda_nb_1 -
\sum_{q+1+t=n}(1^{\tens q}\tens b_1\tens1^{\tens t})\lambda_n = 0. \]
The $n$\n-th component of the identity
\begin{equation*}
\nu B_1 = (\tilde{\uni}\tens\tilde{\uni})B_2B_1 - \tilde{\uni}B_1 =
- (\tilde{\uni}\tens\tilde{\uni})(1\tens B_1+B_1\tens1)B_2 - \lambda
= - (\tilde{\uni}\tens\lambda)B_2 + (\lambda\tens\tilde{\uni})B_2
- \lambda
\end{equation*}
gives
\begin{align*}
\nu_nd &= \nu_nb_1
+ \sum_{q+1+t=n}(1^{\tens q}\tens b_1\tens1^{\tens t})\nu_n
= - (\uni^\cc_0\tens\lambda_n)b_2 + (\lambda_n\tens\uni^\cc_0)b_2
- \lambda_n \\
&= - \lambda_n(\uni^\cc_0\tens1)b_2 + \lambda_n(1\tens\uni^\cc_0)b_2
- \lambda_n = - \lambda_nu'.
\end{align*}
Here the chain map
\[ u' = (\sS{_{X_0}}\uni^\cc_0\tens1)b_2 -
(1\tens\sS{_{X_n}}\uni^\cc_0)b_2 +1: s\cc(X_0,X_n) \to s\cc(X_0,X_n) \]
is homotopic to $-1$ by \lemref{lem-1ib1-i1b-1}. Hence, the map
\[ u = \Hom(N,u'): \Hom^\bull(N,s\cc(X_0,X_n)) \to
\Hom^\bull(N,s\cc(X_0,X_n)), \quad \lambda_n \mapsto \lambda_nu' \]
is also homotopic to $-1$ for each complex of $\kk$\n-modules $N$, in
particular, for
$N=s\cc(X_0,X_1)\tens_\kk\dots\tens_\kk s\cc(X_{n-1},X_n)$. Therefore,
the complex $\Cone(u)$ is contractible by \lemref{lem-contractible}.
Since $-\lambda_nd=0$ and $\nu_nd+\lambda_nu=0$, the element
\[ (\nu_n,\lambda_n) \in \Hom_\kk^{-1}(N,s\cc(X_0,X_n)) \oplus
\Hom_\kk^0(N,s\cc(X_0,X_n)) = \Cone^{-1}(u) \]
is a cycle. Due to acyclicity of $\Cone(u)$ this element
is a boundary of some element
\[ (v_n,\uni^\cc_n) \in \Hom_\kk^{-2}(N,s\cc(X_0,X_n)) \oplus
\Hom_\kk^{-1}(N,s\cc(X_0,X_n)) = \Cone^{-2}(u), \]
that is, $v_nd+\uni^\cc_nu=\nu_n$ and $-\uni^\cc_nd=\lambda_n$. These
equations can be rewritten as follows:
\begin{equation*}
- \uni^\cc_nb_1
- \sum_{q+1+t=n}(1^{\tens q}\tens b_1\tens1^{\tens t})\uni^\cc_n
= (\tilde{\uni}b)_n + (b\tilde{\uni})_n,
\end{equation*}
\begin{multline*}
v_nb_1 - \sum_{q+1+t=n}(1^{\tens q}\tens b_1\tens1^{\tens t})v_n
- (\uni^\cc_0\tens\uni^\cc_n)b_2 - (\uni^\cc_n\tens\uni^\cc_0)b_2
+ \uni^\cc_n \\
= [(\tilde{\uni}\tens\tilde{\uni})B_2]_n - (\tilde{v}b)_n
+ (b\tilde{v})_n.
\end{multline*}
In other words, $(1,1)$-coderivations with components
$(\uni^\cc_0,\dots,\uni^\cc_{n-1},\uni^\cc_n,0,\dots)$,
$(v_0,\dots,v_{n-1}$, $v_n,0,\dots)$ satisfy equations
\eqref{eq-ibm-bim-0}, \eqref{eq-iiB2m-im-vbm-bvm} for $m\le n$. The
construction of $\uni^\cc$, $v$ goes on inductively.
\end{proof}

\begin{definition}\label{def-unit}
A \emph{unit transformation} of an \ainf-category $\cc$ is a natural
\ainf-transformation $\uni^\cc:\id_\cc\to\id_\cc:\cc\to\cc$ such that
$(\uni^\cc\tens\uni^\cc)B_2\equiv\uni^\cc$, and for each pair $X$, $Y$
of objects of $\cc$ the chain maps $(1\tens\sS{_Y}\uni^\cc_0)b_2$,
$(\sS{_X}\uni^\cc_0\tens1)b_2:s\cc(X,Y)\to s\cc(X,Y)$ are homotopy
invertible.
\end{definition}

We have shown in \propref{pro-unital-unit} that an \ainf-category $\cc$
is unital if and only if it has a unit transformation. Similar
(although not identical to our) definitions of units and unital
\ainf-categories are proposed by Kontsevich and
Soibelman~\cite{KonSoi-AinfCat-NCgeom} and
Lef\`evre-Hasegawa~\cite{Lefevre-Ainfty-these}.

\begin{proposition}\label{pro-C-unit-AAC-unit}
Let $\ca$, $\cc$ be \ainf-categories. If $\cc$ is unital, then
$A_\infty(\ca,\cc)$ is unital as well.
\end{proposition}

\begin{proof}
We claim that
\[ (1\tens\uni^\cc)M: (1\tens\id_\cc)M \to (1\tens\id_\cc)M:
A_\infty(\ca,\cc) \to A_\infty(\ca,\cc) \]
is a unit of $A_\infty(\ca,\cc)$. Indeed,
$(1\tens\id_\cc)M=\id_{A_\infty(\ca,\cc)}$, and there is a 3\n-morphism
$v:(\uni^\cc\tens\uni^\cc)B_2\to\uni^\cc:\id_\cc\to\id_\cc:\cc\to\cc$.
Hence, $(\uni^\cc\tens\uni^\cc)B_2=\uni^\cc+vB_1$, and by
\eqref{eq-ttMB-ttBM}
\begin{multline*}
[(1\tens\uni^\cc)M\tens(1\tens\uni^\cc)M]\tilde{B}_2
= [1\tens(\uni^\cc\tens\uni^\cc)B_2]M \\
= (1\tens\uni^\cc)M + (1\tens vB_1)M
= (1\tens\uni^\cc)M + (1\tens v)M\tilde{B}_1
\equiv (1\tens\uni^\cc)M.
\end{multline*}

Let $f:\ca\to\cc$ be an \ainf-functor. The 0\n-th component of
$(1\tens\uni^\cc)M$ is
$\sS{_f}[(1\tens\uni^\cc)M]_0:\kk\to sA_\infty(\ca,\cc)(f,f)$,
$1\mapsto f\uni^\cc$. It remains to prove that for each pair of
\ainf-functors $f,g:\ca\to\cc$ the maps
\begin{align*}
\bigl(1\tens\sS{_g}[(1\tens\uni^\cc)M]_0\bigr)B_2 &:
sA_\infty(\ca,\cc)(f,g) \to sA_\infty(\ca,\cc)(f,g),
\quad r \mapsto (r\tens g\uni^\cc)B_2, \\
\bigl(\sS{_f}[(1\tens\uni^\cc)M]_0\tens1\bigr)B_2 &:
sA_\infty(\ca,\cc)(f,g) \to sA_\infty(\ca,\cc)(f,g),
\quad r \mapsto r(f\uni^\cc\tens1)B_2,
\end{align*}
are homotopy invertible.

Let us define a decreasing filtration of the complex
$(sA_\infty(\ca,\cc)(f,g),B_1)$. For $n\in\ZZ_{\ge0}$, we set
\begin{align*}
\Phi_n &=
\{ r\in sA_\infty(\ca,\cc)(f,g) \mid \forall l<n \quad r_l = 0 \} \\
&= \{ r\in sA_\infty(\ca,\cc)(f,g) \mid
\forall l<n \quad (T^ls\ca)r = 0 \}.
\end{align*}
Clearly, $\Phi_n$ is stable under $B_1=[\_,b]$, and we have
\[ sA_\infty(\ca,\cc)(f,g) = \Phi_0 \supset \Phi_1 \supset\dots\supset
\Phi_n \supset \Phi_{n+1} \supset\dots\quad. \]
Due to \eqref{eq-Bn-components} and \eqref{eq-theta-kl-sum} the chain
maps $a=(1\tens g\uni^\cc)B_2$, $c=(f\uni^\cc\tens1)B_2$ preserve the
subcomplex $\Phi_n$. By \defref{def-ainf-trans}
$sA_\infty(\ca,\cc)(f,g)=\prod_{n=0}^\infty V_n$, where $V_n$ is the
 {\color{blue}graded}
$\kk$\n-module \eqref{eq-Vn-module-rn} of $n$\n-th components $r_n$ of
$(f,g)$-coderivations $r$,
 {\color{blue} and the product is taken in the category of graded
 $\kk$\n-modules}.
The filtration consists of
 {\color{blue}graded}
$\kk$\n-submodules
$\Phi_n=0\times\dots\times0\times\prod_{m=n}^\infty V_m$.

The graded complex associated with this filtration is
$\oplus_{n=0}^\infty V_n$, and the differential $d:V_n\to V_n$ induced
by $B_1$ is given by formula \eqref{eq-rnd-rnb1-b1rn}. The associated
endomorphisms $\gr a$, $\gr c$ of $\oplus_{n=0}^\infty V_n$ are given
by the formulas
\begin{align*}
(r_n)\gr a = (r_n\tens\sS{_g}\uni^\cc_0)b_2 &=
\prod_{X_0,\dots,X_n\in\Ob\ca}
(\sS{_{X_0,\dots,X_n}}r_n\tens\sS{_{X_ng}}\uni^\cc_0)b_2, \\
(r_n)\gr c = r_n(\sS{_f}\uni^\cc_0\tens1)b_2 &=
\prod_{X_0,\dots,X_n\in\Ob\ca}
\sS{_{X_0,\dots,X_n}}r_n(\sS{_{X_0f}}\uni^\cc_0\tens1)b_2,
\end{align*}
$r_n\in V_n$, as formulas \eqref{eq-theta-kl-sum},
\eqref{eq-Bn-components} show. Due to \lemref{lem-1ib1-i1b-1} for each
pair $X$, $Y$ of objects of $\cc$ the chain maps
$(1\tens\sS{_Y}\uni^\cc_0)b_2$,
$-(\sS{_X}\uni^\cc_0\tens1)b_2:s\cc(X,Y)\to s\cc(X,Y)$ are homotopic to
the identity map, that is, $(1\tens\sS{_Y}\uni^\cc_0)b_2=1+hd+dh$,
$(\sS{_X}\uni^\cc_0\tens1)b_2=-1+h'd+dh'$ for some $\kk$\n-linear maps
$h,h':s\cc(X,Y)\to s\cc(X,Y)$ of degree $-1$. Let us choose such
homotopies $\sS{_{X_0,X_n}}h$,
$\sS{_{X_0,X_n}}h':s\cc(X_0f,X_ng)\to s\cc(X_0f,X_ng)$ for each pair
$X_0$, $X_n$ of objects of $\ca$. Denote by
$H,H':\prod_{n=0}^\infty V_n\to\prod_{n=0}^\infty V_n$ the diagonal
maps
$\sS{_{X_0,\dots,X_n}}r_n
\mapsto\sS{_{X_0,\dots,X_n}}r_n\sS{_{X_0,X_n}}h$,
$\sS{_{X_0,\dots,X_n}}r_n
\mapsto\sS{_{X_0,\dots,X_n}}r_n\sS{_{X_0,X_n}}h'$.
Then $\gr a=1+Hd+dH$, $\gr c=-1+H'd+dH'$. The chain maps $a-HB_1-B_1H$,
$c-H'B_1-B_1H'$, being restricted to maps
$\oplus_{m=0}^\infty V_m\to\prod_{m=0}^\infty V_m$ give upper
triangular $\NN\times\NN$ matrices which, in turn, determine the whole
map. Thus, $a-HB_1-B_1H=1+N$, $c-H'B_1-B_1H'=-1+N'$, where the
$\NN\times\NN$ matrices $N$, $N'$ are strictly upper triangular.
Therefore, $1+N$ and $-1+N'$ are invertible (since their inverse maps
$\sum_{i=0}^\infty(-N)^i$ and $-\sum_{i=0}^\infty(N')^i$ make sense).
Hence, $a=(1\tens g\uni^\cc)B_2$ and $c=(f\uni^\cc\tens1)B_2$ are
invertible in $\ck$.
\end{proof}

\begin{corollary}\label{cor-1giB-1-fi1B--1}
Let $f,g:\ca\to\cc$ be \ainf-functors. If $\cc$ is unital, then
\begin{align*}
(1\tens g\uni^\cc)B_2\sim1 &: sA_\infty(\ca,\cc)(f,g) \to
sA_\infty(\ca,\cc)(f,g), \text{ and} \\
(f\uni^\cc\tens1)B_2\sim-1 &: sA_\infty(\ca,\cc)(f,g) \to
sA_\infty(\ca,\cc)(f,g).
\end{align*}
\end{corollary}

\begin{proof}
Follows from \propref{pro-C-unit-AAC-unit} and \lemref{lem-1ib1-i1b-1}.
\end{proof}

\begin{corollary}\label{cor-rgiB-r-firB-r}
Let $r:f\to g:\ca\to\cc$ be a natural \ainf-transformation. If $\cc$ is
unital, then
\[ (r\tens g\uni^\cc)B_2 \equiv r, \qquad
(f\uni^\cc\tens r)B_2 \equiv r. \]
\end{corollary}

\begin{proof}
By \corref{cor-1giB-1-fi1B--1} there are homotopies
$h,h':sA_\infty(\ca,\cc)(f,g)\to sA_\infty(\ca,\cc)(f,g)$, which give
\begin{align*}
(r\tens g\uni^\cc)B_2 = r(1\tens g\uni^\cc)B_2 &= r + rB_1h + rhB_1
= r + (rh)B_1 \equiv r, \\
(f\uni^\cc\tens r)B_2 = - r(f\uni^\cc\tens 1)B_2 &= r + rB_1h' + rh'B_1
= r + (rh')B_1 \equiv r. \qquad \qed
\end{align*}
\renewcommand{\qedsymbol}{}
\end{proof}

\begin{corollary}\label{cor-unit-uniquely}
The unit transformation of a unital category is determined uniquely up
to equivalence.
\end{corollary}

Indeed, take $f=\id_\cc$ and notice that any two unit transformations
$\uni^\cc$ and $\sS'\uni^\cc$ of $\cc$ satisfy
$\sS'\uni^\cc\equiv\sS'\uni^\cc\cdot\uni^\cc\equiv\uni^\cc$.

\begin{corollary}\label{cor-uni-Acat-K-2-cat}
The full $\ck$\n-2-subcategory $\ck\sS{^u}A_\infty$ of non-2-unital
$\ck$\n-2-category $\ck A_\infty$, whose objects are unital
\ainf-categories and the other data are as in $\ck A_\infty$, is a
{\color{red} 1\n-unital left-2-unital} $\ck$\n-2-category. The unit
2\n-endomorphism of an \ainf-functor $f:\ca\to\cb$ is the homotopy
class of the chain map
\[ 1_f:\kk\to(A_\infty(\ca,\cb)(f,f),m_1),
\qquad 1\mapsto (f\uni^\cb)s^{-1}. \]
\end{corollary}

\begin{proof}
The composition
\begin{gather*}
A_\infty(\ca,\cb)(f,g) \rTTo^{1\tens1_g}
A_\infty(\ca,\cb)(f,g)\tens A_\infty(\ca,\cb)(g,g) \rTTo^{m_2}
A_\infty(\ca,\cb)(f,g), \\
rs^{-1} \mapsto rs^{-1}\tens(g\uni^\cb)s^{-1} \mapsto
(rs^{-1}\tens(g\uni^\cb)s^{-1})(s\tens s)B_2s^{-1}
= (r\tens g\uni^\cb)B_2s^{-1},
\end{gather*}
is homotopic to the identity map by \corref{cor-1giB-1-fi1B--1}.
Similarly, the composition
\begin{gather*}
A_\infty(\ca,\cb)(f,g) \rTTo^{1_f\tens1}
A_\infty(\ca,\cb)(f,f)\tens A_\infty(\ca,\cb)(f,g) \rTTo^{m_2}
A_\infty(\ca,\cb)(f,g), \\
rs^{-1} \mapsto (f\uni^\cb)s^{-1}\tens rs^{-1} \mapsto
(-)^{r-1}(f\uni^\cb\tens r)B_2s^{-1} = -r(f\uni^\cb\tens 1)B_2s^{-1},
\end{gather*}
is homotopic to the identity map.

{\color{blue} If \(e:\cd\to\ca\) is an \ainf-functor, then
\(e\cdot1_f=e\cdot(f\uni^\cb)s^{-1}=(ef\uni^\cb)s^{-1}=1_{ef}\).}
\end{proof}

\begin{corollary}\label{cor-uni-Acat-ordin-2-cat}
The full 2\n-subcategory $\sS{^u}A_\infty$ of non-2-unital 2\n-category
$A_\infty$, which consists of unital \ainf-categories, all
\ainf-functors between them, and equivalence classes of all natural
\ainf-transformations is a {\color{red} 1\n-unital left-2-unital}
2\n-category. The unit 2\n-endomorphism of an \ainf-functor
$f:\ca\to\cc$ is $f\uni^\cc$. The notions of an isomorphism between
\ainf-functors, an equivalence between \ainf-categories, etc. make
sense in $\sS{^u}A_\infty$. For instance, $r:f\to g:\ca\to\cb$ is an
\emph{isomorphism} if there is a natural \ainf-transformation
 $p:g\to f$, such that $(r\tens p)B_2\equiv f\uni^\cb$ and
$(p\tens r)B_2\equiv g\uni^\cb$. An \ainf-functor $f:\ca\to\cb$ is an
\emph{equivalence} if there exists an \ainf-functor $g:\cb\to\ca$ and
isomorphisms $\id_\ca\to fg$ and $\id_\cb\to gf$.
\end{corollary}

\begin{proof}
Follows from \corref{cor-uni-Acat-K-2-cat}.
\end{proof}

\subsection{Invertible transformations.}\label{sec-Invertible-trans}
Let $\cb$, $\cc$ be \ainf-categories, and let $f,g:\Ob\cc\to\Ob\cb$ be
maps. Assume that $\cb$ is unital and that for each object $X$ of $\cc$
there are $\kk$\n-linear maps
\begin{alignat*}3
\sS{_X}r_0 &: \kk \to (s\cb)^{-1}(Xf,Xg), &\qquad
\sS{_X}p_0 &: \kk \to (s\cb)^{-1}(Xg,Xf), \\
\sS{_X}w_0 &: \kk \to (s\cb)^{-2}(Xf,Xf), &\qquad
\sS{_X}v_0 &: \kk \to (s\cb)^{-2}(Xg,Xg),
\end{alignat*}
such that
\begin{gather}
\sS{_X}r_0b_1=0, \qquad \sS{_X}p_0b_1=0, \notag \\
(\sS{_X}r_0\tens\sS{_X}p_0)b_2-\sS{_{Xf}}\uni^\cb_0 = \sS{_X}w_0b_1,
\label{eq-rp-i-pr-i} \\
(\sS{_X}p_0\tens\sS{_X}r_0)b_2 - \sS{_{Xg}}\uni^\cb_0 = \sS{_X}v_0b_1.
\notag
\end{gather}

\begin{lemma}\label{lem-r-p-inverse}
Let the above assumptions hold. Then for all objects $X$ of $\cc$ and
$Y$ of $\cb$ the chain maps
\begin{alignat*}3
(r_0\tens1)b_2 &: s\cb(Xg,Y) \to s\cb(Xf,Y) &\text{ and }
(p_0\tens1)b_2 &: s\cb(Xf,Y) \to s\cb(Xg,Y), \\
(1\tens r_0)b_2 &: s\cb(Y,Xf) \to s\cb(Y,Xg) &\text{ and }
(1\tens p_0)b_2 &: s\cb(Y,Xg) \to s\cb(Y,Xf)
\end{alignat*}
are homotopy inverse to each other.
\end{lemma}

\begin{proof}
We have
\begin{align}
& (r_0\tens1)b_2(p_0\tens1)b_2 = (p_0\tens r_0\tens1)(1\tens b_2)b_2
\notag \\
&= -(p_0\tens r_0\tens1)[(b_2\tens1)b_2+b_3b_1+(1\tens1\tens b_1)b_3]
\label{eq-r1b-p1b-pr1b1b} \\
&= -(\sS{_{Xg}}\uni^\cb_0\tens1)b_2 -(v_0b_1\tens1)b_2
-(p_0\tens r_0\tens1)b_3b_1 -b_1(p_0\tens r_0\tens1)b_3 \notag \\
&\sim 1 + b_1(v_0\tens1)b_2 + (v_0\tens1)b_2b_1
\sim 1: s\cb(Xg,Y) \to s\cb(Xg,Y). \notag
\end{align}
For symmetry reasons also
\[ (p_0\tens1)b_2(r_0\tens1)b_2 \sim 1: s\cb(Xf,Y) \to s\cb(Xf,Y). \]
Therefore, $(r_0\tens1)b_2$ and $(p_0\tens1)b_2$ are homotopy inverse
to each other.

Similarly, $(1\tens r_0)b_2$ and $(1\tens p_0)b_2$ are homotopy inverse
to each other.
\end{proof}

\begin{proposition}\label{pro-rfg-p-r-1}
Let $r:f\to g:\cc\to\cb$ be a natural \ainf-transformation, where $\cb$
is unital, and let $p_0$, $v_0$, $w_0$ be as in
\secref{sec-Invertible-trans} so that equations~\eqref{eq-rp-i-pr-i}
hold. Then $p_0$, $w_0$ extend to a natural \ainf-transformation
$p:g\to f:\cc\to\cb$, a 3\n-morphism $w$, and there is a 3\n-morphism
$t$ as follows:
\begin{align}
w: (r\tens p)B_2 \to f\uni^\cb: f \to f &: \cc \to \cb,
\label{eq-w-rp-fi-ff-BC} \\
t: (p\tens r)B_2 \to g\uni^\cb: g \to g &: \cc \to \cb.
\label{eq-t-pr-gi-gg-BC}
\end{align}
In particular, $r$ is invertible and $p=r^{-1}$ in $A_\infty$.
\end{proposition}

\begin{proof}
Let us drop equation~\eqref{eq-t-pr-gi-gg-BC} and prove the existence
of $p$ and $w$, satisfying \eqref{eq-w-rp-fi-ff-BC}. We have to satisfy
the equations
\begin{align}
pb+bp &= 0, \label{eq-pb+bp-0} \\
(r\tens p)B_2 - f\uni^\cb &= [w,b]. \label{eq-rp-B2-fi-wb}
\end{align}
Let us construct the components of $p$ and $w$ by induction. Given a
positive $n$, assume that we have already found components $p_m$, $w_m$
of the sought $p$, $w$ for $m<n$, such that the equations
\begin{equation}
(pb)_m + (bp)_m = 0: s\cc(X_0,X_1)\tens\dots\tens s\cc(X_{m-1},X_m)
\to s\cb(X_0g,X_mf),
\label{eq-pbm-bpm-0}
\end{equation}
\begin{multline}
[(r\tens p)B_2]_m - (f\uni^\cb)_m = (wb)_m - (bw)_m: \\
s\cc(X_0,X_1)\tens\dots\tens s\cc(X_{m-1},X_m)\to s\cb(X_0f,X_mf)
\label{eq-rpB2m-fim-wbm-bwm}
\end{multline}
are satisfied for all $m<n$. Introduce a $(g,f)$-coderivation
$\tilde{p}:Ts\cc\to Ts\cb$ of degree $-1$ by its components
$(p_0,p_1,\dots,p_{n-1},0,0,\dots)$ and an $(f,f)$-coderivation
$\tilde{w}:Ts\cc\to Ts\cb$ of degree $-2$ by its components
$(w_0,w_1,\dots,w_{n-1},0,0,\dots)$. Define a $(g,f)$-coderivation
$\lambda=\tilde{p}b+b\tilde{p}$ of degree $0$ and an
$(f,f)$-coderivation $\nu=(r\tens\tilde{p})B_2-f\uni^\cb-[\tilde{w},b]$
of degree $-1$. Then equations \eqref{eq-pb+bp-0},
\eqref{eq-rp-B2-fi-wb} imply that $\lambda_m=0$, $\nu_m=0$ for $m<n$.
The identity $\lambda b-b\lambda=0$ implies that
\[ \lambda_nd = \lambda_nb_1 -
\sum_{q+1+t=n}(1^{\tens q}\tens b_1\tens1^{\tens t})\lambda_n = 0. \]
The identity
\begin{multline*}
[\nu,b] = \nu B_1
= (r\tens\tilde{p})B_2B_1 - (f\uni^\cb)B_1 - \tilde{w}B_1B_1
= - (r\tens\tilde{p})(1\tens B_1+B_1\tens1)B_2 \\
= - (r\tens\tilde{p}B_1)B_2 = - (r\tens\lambda)B_2
\end{multline*}
implies that
\[ \nu_nb_1 + \sum_{q+1+t=n}(1^{\tens q}\tens b_1\tens1^{\tens t})\nu_n
= - (r_0\tens\lambda_n)b_2,
\]
that is, $\nu_nd=-\lambda_nu$. Here the map $u=\Hom(N,(r_0\tens1)b_2)$
for $N=s\cc(X_0,X_1)\tens_\kk\dots\tens_\kk s\cc(X_{n-1},X_n)$ is
homotopy invertible, and the complex $\Cone(u)$ is contractible by
\lemref{lem-contractible}. Hence, the cycle
\[ (\nu_n,\lambda_n) \in \Hom_\kk^{-1}(N,s\cb(X_0f,X_nf)) \oplus
\Hom_\kk^0(N,s\cb(X_0g,X_nf)) = \Cone^{-1}(u)
\]
is the boundary of some element
\[ (w_n,p_n) \in \Hom_\kk^{-2}(N,s\cb(X_0f,X_nf)) \oplus
\Hom_\kk^{-1}(N,s\cb(X_0g,X_nf)) = \Cone^{-2}(u), \]
that is, $w_nd+p_nu=\nu_n$ and $-p_nd=\lambda_n$. In other words,
equations \eqref{eq-pbm-bpm-0}, \eqref{eq-rpB2m-fim-wbm-bwm} are
satisfied for $m=n$, and we prove the statement by induction.

For similar reasons using \lemref{lem-r-p-inverse} there exists a
natural \ainf-transformation $q:g\to f:\cc\to\cb$ with $q_0=p_0$ and a
3\n-morphism
\[ v: (q\tens r)B_2 \to g\uni^\cb: g \to g: \cc \to \cb \]
with given $v_0$. Since $r$ has a left inverse and a right inverse, it
is invertible in $A_\infty$ and $p$ is equivalent to $q$. Hence,
$(p\tens r)B_2$ is equivalent to $(q\tens r)B_2$, and there exists $t$
of \eqref{eq-t-pr-gi-gg-BC}.
\end{proof}

\section{\texorpdfstring{Unital $A_\infty$-functors}{Unital A8-functors}}
\begin{definition}\label{def-Unital-functors}
Let $\ca$, $\cb$ be unital \ainf-categories. An \ainf-functor
$f:\ca\to\cb$ is called \emph{unital} if for all objects $X$ of $\ca$
we have $\sS{_X}\uni^\ca_0f_1-\sS{_{Xf}}\uni^\cb_0\in\im b_1$.
\end{definition}

For instance, an \ainf-homomorphism $f:\ca\to\cb$ of \ainf-algebras is
unital if the cycles $\uni^\ca_0f_1$, $\uni^\cb_0\in(s\cb)^{-1}$ are
cohomologous in $(s\cb,b_1)$. We may say that a unital \ainf-functor
(or \ainf-homomorphism) preserves the cohomology classes of unit
elements.

\begin{proposition}\label{pro-unital-functors}
Let $\ca$, $\cb$ be unital \ainf-categories. An \ainf-functor
$f:\ca\to\cb$ is unital if and only if $\uni^\ca f\equiv f\uni^\cb$.
\end{proposition}

\begin{proof}
If $\uni^\ca f=f\uni^\cb+vB_1$, then
$\sS{_X}\uni^\ca_0f_1=\sS{_X}(\uni^\ca f)_0=\sS{_X}(f\uni^\cb+vB_1)_0
=\sS{_{Xf}}\uni^\cb_0+\sS{_X}v_0b_1$,
hence, $f$ is unital.

Assume now that $f$ is unital. We want to find a 3\n-morphism
\[ v: \uni^\ca f\to f\uni^\cb: f\to f: \ca\to\cb, \]
that is, an $(f,f)$-coderivation $v$ of degree $-2$ such that
\begin{equation}
vB_1 = \uni^\ca f - f\uni^\cb.
\label{eq-vB1-iAf-fiB}
\end{equation}
We subject it to an additional condition described below. Consider
3\n-morphisms
\begin{alignat*}2
x &: (\uni^\ca\tens\uni^\ca)B_2 \to \uni^\ca\, &: \id_\ca \to \id_\ca
&: \ca \to \ca, \\
y &: (\uni^\cb\tens\uni^\cb)B_2 \to \uni^\cb   &: \id_\cb \to \id_\cb
&: \cb \to \cb,
\end{alignat*}
so that
\[ xB_1 = (\uni^\ca\tens\uni^\ca)B_2 - \uni^\ca, \qquad
yB_1 = (\uni^\cb\tens\uni^\cb)B_2 - \uni^\cb.
\]
The following equations between $(f,f)$-coderivations of degree $-2$
are due to \eqref{eq-rkpkB2-M20-rpB2k}, \eqref{eq-krkpB2-M02-krpB2}:
\begin{gather*}
(xf)B_1 = (\uni^\ca\tens\uni^\ca)B_2f - \uni^\ca f
= (\uni^\ca f\tens\uni^\ca f)B_2
+ (\uni^\ca\tens\uni^\ca\mid f)M_{20}B_1 - \uni^\ca f, \\
(fy)B_1 = f(\uni^\cb\tens\uni^\cb)B_2 - f\uni^\cb
= (f\uni^\cb\tens f\uni^\cb)B_2 - f\uni^\cb.
\end{gather*}
Combining them with \eqref{eq-vB1-iAf-fiB} we find that
\begin{align*}
&(xf)B_1 - (\uni^\ca\tens\uni^\ca\mid f)M_{20}B_1 + vB_1 \\
&= (\uni^\ca f\tens\uni^\ca f)B_2 - f\uni^\cb \\
&= (\uni^\ca f\tens vB_1)B_2 + (\uni^\ca f\tens f\uni^\cb)B_2
- f\uni^\cb \\
&= - (\uni^\ca f\tens v)B_2B_1 + (vB_1\tens f\uni^\cb)B_2
+ (f\uni^\cb\tens f\uni^\cb)B_2 - f\uni^\cb \\
&= - (\uni^\ca f\tens v)B_2B_1 + (v\tens f\uni^\cb)B_2B_1 + (fy)B_1.
\end{align*}
Now we may formulate the problem: we are looking for $v$ as above and
an $(f,f)$-coderivation $w$ of degree $-3$, such that
\[ wB_1 = xf - (\uni^\ca\tens\uni^\ca\mid f)M_{20} + v
+ (\uni^\ca f\tens v)B_2 - (v\tens f\uni^\cb)B_2 - fy,
\]
in other terms, a 4\n-morphism
\begin{multline*}
w: xf - (\uni^\ca\tens\uni^\ca\mid f)M_{20} + v \to
fy - (\uni^\ca f\tens v)B_2 + (v\tens f\uni^\cb)B_2: \\
(\uni^\ca f\tens\uni^\ca f)B_2 \to f\uni^\cb: f \to f: \ca \to \cb.
\end{multline*}
Using the chain map
\[ u = (1\tens f\uni^\cb)B_2 -1 - (\uni^\ca f\tens1)B_2:
(sA_\infty(\ca,\cb)(f,f),B_1) \to (sA_\infty(\ca,\cb)(f,f),B_1),
\]
we may rewrite our system of equations as follows:
\begin{align}
-vB_1 &= f\uni^\cb - \uni^\ca f, \notag \\
wB_1 + vu &= xf - (\uni^\ca\tens\uni^\ca\mid f)M_{20} - fy.
\label{eq-diff-vw-fiif-xfiifMfy}
\end{align}
In other words, we look for an element
\[ (w,v)\in
[sA_\infty(\ca,\cb)(f,f)]^{-3}\oplus[sA_\infty(\ca,\cb)(f,f)]^{-2}
= \Cone^{-3}(u), \]
whose boundary is
\begin{multline*}
(xf -(\uni^\ca\tens\uni^\ca\mid f)M_{20}-fy, f\uni^\cb -\uni^\ca f) \\
\in [sA_\infty(\ca,\cb)(f,f)]^{-2}\oplus[sA_\infty(\ca,\cb)(f,f)]^{-1}
= \Cone^{-2}(u).
\end{multline*}
Let us prove that $u$ is homotopy invertible. Since
\[ \sS{_X}\uni^\ca_0f_1
=\sS{_{Xf}}\uni^\cb_0+\sS{_X}zb_1:\kk\to(s\cb)^{-1}(Xf,Xf), \]
for some $\sS{_X}z$, the cycles $\sS{_X}r_0=\sS{_X}\uni^\ca_0f_1$ and
$\sS{_X}p_0=\sS{_{Xf}}\uni^\cb_0$ satisfy
conditions~\eqref{eq-rp-i-pr-i} for $g=f:\Ob\ca\to\Ob\cb$, that is,
\begin{gather*}
(\sS{_X}\uni^\ca_0f_1\tens\sS{_{Xf}}\uni^\cb_0)b_2-\sS{_{Xf}}\uni^\cb_0
=(\sS{_{Xf}}\uni^\cb_0\tens\sS{_{Xf}}\uni^\cb_0)b_2-\sS{_{Xf}}\uni^\cb_0
+ (\sS{_X}z\tens\sS{_{Xf}}\uni^\cb_0)b_2b_1 \in\im b_1, \\
(\sS{_{Xf}}\uni^\cb_0\tens\sS{_X}\uni^\ca_0f_1)b_2-\sS{_{Xf}}\uni^\cb_0
=(\sS{_{Xf}}\uni^\cb_0\tens\sS{_{Xf}}\uni^\cb_0)b_2-\sS{_{Xf}}\uni^\cb_0
- (\sS{_{Xf}}\uni^\cb_0\tens\sS{_X}z)b_2b_1 \in\im b_1.
\end{gather*}
Hence, the natural \ainf-transformation $r=\uni^\ca f:f\to f:\ca\to\cb$
is invertible by \propref{pro-rfg-p-r-1}. In detail, there exists a
natural \ainf-transformation $p:f\to f:\ca\to\cb$ and 3\n-morphisms
$q$, $t$ such that
\begin{equation}
(\uni^\ca f\tens p)B_2 - f\uni^\cb = qB_1, \qquad
(p\tens\uni^\ca f)B_2 - f\uni^\cb = tB_1.
\label{eq-ifpB2-fiqB1-pifB2-fitB1}
\end{equation}
These equations are interpreted as equations~\eqref{eq-rp-i-pr-i} for
the following data. Let $\cc=\1$ be a 1-object-0-morphisms
\ainf-category, $\Ob\cc=\{*\}$, $\cc(*,*)=0$. Consider a map
$\Ob\cc\to\Ob A_\infty(\ca,\cb)$, $*\mapsto f$, and elements
$\uni^\ca f,p\in[sA_\infty(\ca,\cb)(f,f)]^{-1}$,
$q,t\in[sA_\infty(\ca,\cb)(f,f)]^{-2}$. Equations~\eqref{eq-rp-i-pr-i}
for these data are precisely \eqref{eq-ifpB2-fiqB1-pifB2-fitB1}, since
$\sS{_f}\uni^{A_\infty(\ca,\cb)}_0$\linebreak[1]${}=f\uni^\cb$. By
\lemref{lem-r-p-inverse} we deduce that
\[ (\uni^\ca f\tens1)B_2: (sA_\infty(\ca,\cb)(f,f),B_1)
\to (sA_\infty(\ca,\cb)(f,f),B_1),
\]
is homotopy invertible. Since the map $(1\tens f\uni^\cb)B_2-1$ is
homotopic to 0 by \corref{cor-1giB-1-fi1B--1}, we deduce that $u$ is
homotopy invertible. Therefore, $\Cone(u)$ is contractible by
\lemref{lem-contractible}.

To prove the existence of $(w,v)$ satisfying
\eqref{eq-diff-vw-fiif-xfiifMfy} it suffices to show that
$(xf-(\uni^\ca\tens\uni^\ca\mid f)M_{20}-fy,
f\uni^\cb-\uni^\ca f)\in\Cone^{-2}(u)$
is a cycle. And indeed,
\begin{align*}
& [xf-(\uni^\ca\tens\uni^\ca\mid f)M_{20}-fy]B_1
+ (f\uni^\cb-\uni^\ca f)u \\
&= (xB_1)f - (\uni^\ca\tens\uni^\ca\mid f)M_{20}B_1 - f(yB_1) \\
&\qquad + [(f\uni^\cb-\uni^\ca f)\tens f\uni^\cb]B_2 - f\uni^\cb
+ \uni^\ca f + [\uni^\ca f\tens(f\uni^\cb-\uni^\ca f)]B_2 \\
&= (\uni^\ca\tens\uni^\ca)B_2f - \uni^\ca f
- (\uni^\ca\tens\uni^\ca\mid f)M_{20}B_1 - f(\uni^\cb\tens\uni^\cb)B_2
+ f\uni^\cb \\
&\qquad + (f\uni^\cb\tens f\uni^\cb)B_2
- (\uni^\ca f\tens f\uni^\cb)B_2 - f\uni^\cb + \uni^\ca f
+ (\uni^\ca f\tens f\uni^\cb)B_2 - (\uni^\ca f\tens\uni^\ca f)B_2 \\
&= - [(\uni^\ca\tens\uni^\ca)(1\tens B_1+B_1\tens1)\mid f]M_{20} = 0
\end{align*}
due to \eqref{eq-rkpkB2-M20-rpB2k} and \eqref{eq-krkpB2-M02-krpB2}.
Clearly, $(f\uni^\cb-\uni^\ca f)B_1=0$, so the proposition is proven.
\end{proof}

\begingroup
\color{blue}
\begin{corollaryNoNumber}
The $\ck$\n-2-subcategory $\ck A_\infty^u$ of 1\n-unital left-2-unital
$\ck$\n-2-category $\ck\sS{^u}A_\infty$, whose objects are unital
\ainf-categories, 1\n-morphisms are unital \ainf-functors and other
data are as in $\ck\sS{^u}A_\infty$, is a 1\n-2-unital
$\ck$\n-2-category.
\end{corollaryNoNumber}

\begin{proof}
Clearly, the composition of unital functors is unital.

If \(\cd\xrightarrow g\ca\xrightarrow f\cb\) are \ainf-functors and $f$
is unital, then the chain maps
\begin{align*}
1_g\cdot f =(g\uni^\ca s^{-1})\cdot f =(g\uni^\ca f)s^{-1} &: \kk \to
(A_\infty(\cd,\cb)(gf,gf),m_1),
\\
1_{gf} =(gf\uni^\cb)s^{-1} &: \kk \to (A_\infty(\cd,\cb)(gf,gf),m_1)
\end{align*}
are equal in $\ck$. In fact, there is \(v\in sA_\infty(\ca,\cb)(f,f)\)
such that \(g\uni^\ca f=gf\uni^\cb+(gv)B_1\).
\end{proof}
\endgroup

If $\cb$, $\cc$ are unital \ainf-categories, $r:f\to g:\cb\to\cc$ is an
isomorphism of \ainf-functors and $f$ is unital, then $g$ is unital as
well. Indeed, distributivity law in $\sS{^u}A_\infty$ implies
\begin{equation*}
\cb
\pile{\rTTo^f_{r\Downarrow}\\ \rTTo~g\\ \rTTo^{\uni^\cb g\Downarrow}_g}
\cc = \biggl(\cb
\pile{\rTTo^{\id_\cb}_{\uni^\cb\Downarrow}\\ \rTTo_{\id_\cb}}
\cb \pile{\rTTo^f_{r\Downarrow}\\ \rTTo_g} \cc\biggr) = \cb
\pile{\rTTo^f_{\uni^\cb f\Downarrow\equiv f\uni^\cc}\\ \rTTo~f\\
\rTTo^{r\Downarrow}_g}
\cc = \biggl(\cb \pile{\rTTo^f_{r\Downarrow}\\ \rTTo_g} \cc
\pile{\rTTo^{\id_\cc}_{\uni^\cc\Downarrow}\\ \rTTo_{\id_\cc}}
\cc\biggr) = \cb
\pile{\rTTo^f_{r\Downarrow}\\ \rTTo~g\\ \rTTo^{g\uni^\cc\Downarrow}_g}
\cc,
\end{equation*}
\begin{equation}
\text{or}\qquad r\cdot(\uni^\cb g) = \uni^\cb\hcirc r = (\uni^\cb f)\cdot r
= (f\uni^\cc)\cdot r = r\hcirc \uni^\cc = r\cdot(g\uni^\cc),
\label{eq-riBg-rgiC}
\end{equation}
where $\cdot$ and $\hcirc$ denote the vertical and the horizontal
compositions of 2\n-morphisms, hence,
 \((r\tens\uni^\cb g)B_2\equiv(r\tens g\uni^\cc)B_2\).
{\color{blue} Multiplying on the left with an \ainf-transformation
inverse to $r$ we get
 \(\uni^\cb g\equiv(g\uni^\cc\tens\uni^\cb g)B_2
 \equiv(g\uni^\cc\tens g\uni^\cc)B_2\equiv g\uni^\cc\).}

\begin{definition}\label{def-Ainfuu}
The 2\n-category $A_\infty^u$ is a 2\n-subcategory of
$\sS{^u}A_\infty$, whose class of objects consists of all unital
\ainf-categories, 1\n-morphisms are all unital \ainf-functors, and
2\n-morphisms are equivalence classes of all natural
\ainf-transformations between such functors.
\end{definition}

{\color{blue} $A_\infty^u$ is a usual 2\n-category by the above
corollary.}

\begin{proposition}\label{pro-AA-Au-Au-2-functor}
Let $\ca$ be an \ainf-category. The strict 2\n-functor
$A_\infty(\ca,\_)$ maps a unital \ainf-category $\cc$ to the unital
\ainf-category $A_\infty(\ca,\cc)$, and a unital functor to a unital
functor. Its restrictions
$A_\infty(\ca,\_):\sS{^u}A_\infty\to\sS{^u}A_\infty$,
$A_\infty(\ca,\_):A_\infty^u\to A_\infty^u$ preserve 1\n-units and
2\n-units, thus, they are 2\n-functors in the usual sense.
\end{proposition}

\begin{proof}
\propref{pro-C-unit-AAC-unit} shows that $A_\infty(\ca,\cc)$ is unital,
if $\cc$ is unital. If $g:\cb\to\cc$ is a unital \ainf-functor between
unital \ainf-categories $\cb$ and $\cc$, then
$\uni^\cb g\equiv g\uni^\cc$ implies
\begin{equation}
(1\tens\uni^\cb)M(1\tens g)M = (1\tens\uni^\cb g)M
\equiv (1\tens g\uni^\cc)M = (1\tens g)M(1\tens\uni^\cc)M,
\label{eq-iMgM-gMiM}
\end{equation}
hence, $(1\tens g)M$ is unital. The fact, that $A_\infty(\ca,\_)$
preserves 1\n-units and 2\n-units is already proven in
\propref{pro-C-unit-AAC-unit}.
\end{proof}

\subsection{Categories modulo homotopy.}
$\ck$\n-categories form a 2\n-category $\KCat$. We consider also
non-unital $\ck$\n-categories. They form a 2\n-category $\KCat^{nu}$
without 2\n-units (but with 1\n-units -- identity functors).

\begin{proposition}\label{pro-k-Ainf-Kcat}
There is a strict 2\n-functor $\kf:A_\infty\to\KCat^{nu}$ of
non-2-unital 2\n-categories, which assigns to an \ainf-category $\cc$
the $\ck$\n-category $\kf\cc$ with the same class of objects
$\Ob\kf\cc=\Ob\cc$, the same graded $\kk$\n-module of morphisms
$\kf\cc(X,Y)=\cc(X,Y)$, equipped with the differential $m_1$.
Composition in $\kf\cc$ is given by (the homotopy equivalence class of)
$m_2:\cc(X,Y)\tens\cc(Y,Z)\to\cc(X,Z)$. To an \ainf-functor
$f:\ca\to\cb$ is assigned $\kf f:\kf\ca\to\kf\cb$ such that
$\Ob\kf f=f:\Ob\ca\to\Ob\cb$, and for each pair of objects $X$, $Y$ of
$\ca$ we have $\kf f=sf_1s^{-1}:\ca(X,Y)\to\cb(Xf,Yf)$. To a natural
\ainf-transformation $r:f\to g:\ca\to\cb$ is assigned
$\kf r=r_0s^{-1}:\kf f\to\kf g$, that is, for each object $X$ of $\ca$
the component $\sS{_X}\kf r$ is the homotopy equivalence class of chain
map $\sS{_X}r_0s^{-1}:\kk\to\cb(Xf,Xg)$. Unital \ainf-categories and
unital \ainf-functors are mapped by $\kf$ to unital $\ck$\n-categories
and unital $\ck$\n-functors. The restriction
$\kf:A_\infty^u\to\KCat$ is a 2\n-functor, which preserves 1\n-units
and 2\n-units.
\end{proposition}

\begin{proof}
The identity \eqref{eq-m2-associative-K} shows that $m_2$ is
associative in $\ck$. The identity
\[(\kf f\tens\kf f)m_2 + (s\tens s)f_2s^{-1}m_1
+ (1\tens m_1+m_1\tens1)(s\tens s)f_2s^{-1} - m_2\kf f = 0 \]
shows that $\kf f$ preserves the multiplication in $\ck$.

The map $\sS{_X}\kf r$ is a chain map since
$\sS{_X}r_0s^{-1}m_1=\sS{_X}r_0b_1s^{-1}=0$. If
$r\equiv p:f\to g:\ca\to\cb$, then
$\sS{_X}r_0=\sS{_X}p_0+\sS{_X}v_0b_1$ for some
$\sS{_X}v_0\in(s\cb)^{-2}(Xf,Xg)$, therefore,
$\sS{_X}r_0s^{-1}=\sS{_X}p_0s^{-1}+(\sS{_X}v_0s^{-1})m_1$ and chain
maps $\sS{_X}r_0s^{-1}$ and $\sS{_X}p_0s^{-1}$ are homotopic to each
other, that is, $\kf r=\kf p$. The identity
\begin{multline*}
0 = s[(f_1\tens\sS{_Y}r_0)b_2 + (\sS{_X}r_0\tens g_1)b_2 + r_1b_1
+ b_1r_1]s^{-1} \\
= (\kf f\tens\sS{_Y}\kf r)m_2 - (\sS{_X}\kf r\tens\kf g)m_2
+ sr_1s^{-1}m_1 + m_1sr_1s^{-1}: \ca(X,Y) \to \cb(Xf,Yg)
\end{multline*}
shows that the following diagram commutes in $\ck$ for all objects $X$,
$Y$ of $\ca$
\begin{diagram}
\ca(X,Y) & \rTTo^{\kf f} & \cb(Xf,Yf) \\
\dTTo<{\kf g} && \dTTo>{(1\tens\sS{_Y}\kf r)m_2} \\
\cb(Xg,Yg) & \rTTo^{(\sS{_X}\kf r\tens1)m_2} & \cb(Xf,Yg)
\end{diagram}
Thus $\kf r$ is, indeed, a $\ck$\n-natural transformation.

One checks easily that the composition of functors is preserved, and
the both compositions of 1\n-morphisms and 2\n-morphisms are preserved.
The vertical composition of 2\n-morphisms is preserved due to the
property
\[ \sS{_X}\kf[(r\tens p)B_2] = \sS{_X}[(r\tens p)B_2]_0s^{-1}
= (\sS{_X}r_0\tens\sS{_X}p_0)b_2s^{-1}
= (\sS{_X}\kf r\tens\sS{_X}\kf p)m_2. \]

Let $\cb$ be a unital category. Then $\kf\cb$ is a unital
$\ck$\n-category. Indeed, for each object $X$ of $\cb$ consider the
corresponding element
$1_X=\sS{_X}\uni^\cb_0s^{-1}=\sS{_X}\kf\uni^\cb:\kk\to\cb^0(X,X)$. Then
for each pair $X$, $Y$ of objects of $\cc$ the following equations hold
in $\ck$
\begin{gather*}
(1\tens1_Y)m_2 = s(1\tens1_Ys)b_2s^{-1}
= s(1\tens\sS{_Y}\uni^\cb_0)b_2s^{-1} = ss^{-1} = 1:
\cb(X,Y) \to \cb(X,Y), \\
(1_X\tens1)m_2 = -s(1_Xs\tens1)b_2s^{-1}
= -s(\sS{_X}\uni^\cb_0\tens1)b_2s^{-1} = ss^{-1} = 1:
\cb(X,Y) \to \cb(X,Y).
\end{gather*}
That is, $1_X$ is the unit endomorphism of $X$.

If $f:\ca\to\cc$ is unital then, applying $\kf$ to the equivalence
$\uni^\ca f\equiv f\uni^\cc$, we find that
\( (1_{\id_{\kf\ca}})\kf f = (\kf f)1_{\id_{\kf\cc}} = 1_{\kf f} \),
that is, $\kf f$ is unital (it maps units into units).
\end{proof}

\begin{lemma}[Cancellation]\label{lem-cancel-phi}
Let $\phi:\cc\to\cb$ be an \ainf-functor, such that for all objects
$X$, $Y$ of $\cc$ the chain map
$\phi_1:(s\cc(X,Y),b_1)\to(s\cb(X\phi,Y\phi),b_1)$ is invertible in
$\ck$. Let $f,g:\ca\to\cc$ be \ainf-functors. Let
$y:f\phi\to g\phi:\ca\to\cb$ be a natural \ainf-transformation. Then
there is a unique up to equivalence natural \ainf-transformation
$t:f\to g:\ca\to\cc$ such that $y\equiv t\phi$.
\end{lemma}

\begin{proof}
First we prove the existence. We are looking for a 2\n-morphism
$t:f\to g:\ca\to\cc$ and a 3\n-morphism
$v:y\to t\phi:f\phi\to g\phi:\ca\to\cb$. We have to satisfy the
equations
\[ tb + bt = 0, \qquad y - t\phi = vb - bv. \]
Let us construct the components of $t$ and $v$ by induction. Given a
non-negative integer $n$, assume that we have already found components
$t_m$, $v_m$ of the sought $t$, $v$ for $m<n$, such that the equations
\begin{gather}
(tb)_m + (bt)_m = 0: s\ca(X_0,X_1)\tens\dots\tens s\ca(X_{m-1},X_m)
\to s\cc(X_0f,X_mg), \label{eq-tbm-btm-0-A-C} \\
y_m-(t\phi)_m=(vb-bv)_m: s\ca(X_0,X_1)\tens\dots\tens s\ca(X_{m-1},X_m)
\to s\cb(X_0f\phi,X_mg\phi), \label{eq-ym-tfm-vbbvm-A-C}
\end{gather}
are satisfied for all $m<n$. Introduce an $(f,g)$-coderivation
$\tilde{t}:Ts\ca\to Ts\cc$ of degree $-1$ by its components
$(t_0,\dots,t_{n-1},0,0,\dots)$ and an $(f\phi,g\phi)$-coderivation
$\tilde{v}:Ts\ca\to Ts\cb$ of degree $-2$ by its components
$(v_0,\dots,v_{n-1},0,0,\dots)$. Define an $(f,g)$-coderivation
$\lambda=\tilde{t}b+b\tilde{t}$ of degree $0$ and an
$(f\phi,g\phi)$-coderivation
$\nu=y-\tilde{t}\phi-\tilde{v}b+b\tilde{v}$ of degree $-1$. Then
equations \eqref{eq-tbm-btm-0-A-C}, \eqref{eq-ym-tfm-vbbvm-A-C} imply
that $\lambda_m=0$, $\nu_m=0$ for $m<n$. The identity
$\lambda b-b\lambda=0$ implies that
\begin{equation}
\lambda_nd = \lambda_nb_1 - \sum_{\alpha+1+\beta=n}
(1^{\tens\alpha}\tens b_1\tens1^{\tens\beta})\lambda_n = 0.
\label{eq-land-lnb-bla}
\end{equation}
The identity
\[ \nu b + b\nu = - \tilde{t}\phi b - b\tilde{t}\phi = - \lambda\phi \]
implies that
\begin{equation}
\nu_nd = \nu_nb_1 + \sum_{\alpha+1+\beta=n}
(1^{\tens\alpha}\tens b_1\tens1^{\tens\beta})\nu_n = - \lambda_n\phi_1.
\label{eq-nud-nb-bn--laf}
\end{equation}
Denote $N=s\ca(X_0,X_1)\tens_\kk\dots\tens_\kk s\ca(X_{n-1},X_n)$, and
consider the chain map
\[ u = \Hom(N,\phi_1): \Hom^\bull(N,s\cc(X_0f,X_ng)) \to
\Hom^\bull(N,s\cc(X_0f\phi,X_ng\phi)).
\]
Since $\phi_1$ is homotopy invertible, the map $u$ is homotopy
invertible as well. Therefore, the complex $\Cone(u)$ is acyclic.
Moreover, it is contractible by \lemref{lem-contractible}. Equations
\eqref{eq-land-lnb-bla} and \eqref{eq-nud-nb-bn--laf} in the form
$-\lambda_nd=0$, $\nu_nd+\lambda_n\phi_1=0$ imply that
\[ (\nu_n,\lambda_n) \in \Hom_\kk^{-1}(N,s\cb(X_0f\phi,X_ng\phi))
\oplus \Hom_\kk^0(N,s\cc(X_0f,X_ng)) = \Cone^{-1}(u) \]
is a boundary of some element
\[ (v_n,t_n) \in \Hom_\kk^{-2}(N,s\cb(X_0f\phi,X_ng\phi)) \oplus
\Hom_\kk^{-1}(N,s\cc(X_0f,X_ng)) = \Cone^{-2}(u), \]
that is, $v_nd+t_nu=\nu_n$ and $-t_nd=\lambda_n$. In other words,
equations \eqref{eq-tbm-btm-0-A-C}, \eqref{eq-ym-tfm-vbbvm-A-C} are
satisfied for $m=n$, and we prove the existence of $t$ with the
required properties by induction.

Now we prove the uniqueness of $t$. Assume that we have 2\n-morphisms
$t,t':f\to g:\ca\to\cc$ and 3\n-morphisms
$v:y\to t\phi:f\phi\to g\phi:\ca\to\cb$,
$v':y\to t'\phi:f\phi\to g\phi:\ca\to\cb$. We look for a 3\n-morphism
$w$ and a 4\n-morphism $x$:
\begin{gather*}
w: t \to t': f \to g: \ca \to \cc, \\
x: v' \to v+w\phi: y \to t'\phi: f\phi \to g\phi: \ca \to \cb.
\end{gather*}
They have to satisfy equations
\[ t - t' = wb - bw, \qquad v' - v - w\phi = xb + bx. \]
Let us construct the components of $w$ and $x$ by induction. Given a
non-negative integer $n$, assume that we have already found components
$w_m$ and $x_m$ of the sought $w$, $x$ for $m<n$, such that the
equations
\begin{gather}
t_m - t'_m = (wb-bw)_m : s\ca(X_0,X_1)\tens\dots\tens s\ca(X_{m-1},X_m)
\to s\cc(X_0f,X_mg), \label{eq-tbm-btm-wb-bwm-A-C} \\
v'_m - v_m - (w\phi)_m = (xb+bx)_m:
s\ca(X_0,X_1)\tens\dots\tens s\ca(X_{m-1},X_m)
\to s\cb(X_0f\phi,X_mg\phi), \label{eq-vm-wfm-xbbxm-A-B}
\end{gather}
are satisfied for all $m<n$. Introduce an $(f,g)$-coderivation
$\tilde{w}:Ts\ca\to Ts\cc$ of degree $-2$ by its components
$(w_0,\dots,w_{n-1},0,0,\dots)$ and an $(f\phi,g\phi)$-coderivation
$\tilde{x}:Ts\ca\to Ts\cb$ of degree $-3$ by its components
$(x_0,\dots,x_{n-1},0,0,\dots)$. Define an $(f,g)$-coderivation
$\lambda=t'-t+\tilde{w}b+b\tilde{w}$ of degree $-1$ and an
$(f\phi,g\phi)$-coderivation
$\nu=v'-v-\tilde{w}\phi-\tilde{x}b-b\tilde{x}$ of degree $-2$. Then
equations \eqref{eq-tbm-btm-wb-bwm-A-C}, \eqref{eq-vm-wfm-xbbxm-A-B}
imply that $\lambda_m=0$, $\nu_m=0$ for $m<n$. The identity
$\lambda b+b\lambda=0$ implies that
\begin{equation*}
\lambda_nd = \lambda_nb_1 + \sum_{\alpha+1+\beta=n}
(1^{\tens\alpha}\tens b_1\tens1^{\tens\beta})\lambda_n = 0.
\end{equation*}
The identity
\[ \nu b - b\nu = v'B_1 - vB_1 - \tilde{w}\phi b + b\tilde{w}\phi
= y - t'\phi - y + t\phi - \tilde{w}b\phi + b\tilde{w}\phi
= - \lambda\phi \]
implies that
\begin{equation*}
\nu_nd = \nu_nb_1 - \sum_{\alpha+1+\beta=n}
(1^{\tens\alpha}\tens b_1\tens1^{\tens\beta})\nu_n = - \lambda_n\phi_1.
\end{equation*}
Hence,
\[ (\nu_n,\lambda_n) \in \Hom_\kk^{-2}(N,s\cb(X_0f\phi,X_ng\phi))
\oplus \Hom_\kk^{-1}(N,s\cc(X_0f,X_ng)) = \Cone^{-2}(u) \]
is a cycle, therefore, it is a boundary of an element
\[ (x_n,w_n) \in \Hom_\kk^{-3}(N,s\cb(X_0f\phi,X_ng\phi)) \oplus
\Hom_\kk^{-2}(N,s\cc(X_0f,X_ng)) = \Cone^{-3}(u), \]
that is, $x_nd+w_n\phi_1=\nu_n$ and $-w_nd=\lambda_n$. In other words,
equations \eqref{eq-tbm-btm-wb-bwm-A-C}, \eqref{eq-vm-wfm-xbbxm-A-B},
are satisfied for $m=n$, and we prove the uniqueness of $t$, using
induction.
\end{proof}

A version of the following theorem is proved by Fukaya
\cite[Theorem~8.6]{Fukaya:FloerMirror-II} with a different notion of
unitality and under the additional assumption that the $\kk$\n-modules
$\cb(W,Z)$, $\cc(X,Y)$ are free.

\begin{theorem}\label{thm-phi-C-B-hinv-equiv}
Let $\cc$ be an \ainf-category and let $\cb$ be a unital
\ainf-category. Let $\phi:\cc\to\cb$ be an \ainf-functor such that
for all objects $X$, $Y$ of $\cc$ the chain map
$\phi_1:(s\cc(X,Y),b_1)\to(s\cb(X\phi,Y\phi),b_1)$ is invertible
in $\ck$. Let $h:\Ob\cb\to\Ob\cc$ be a mapping. Assume that for
each object $U$ of $\cb$ the $\kk$\n-linear maps
\begin{alignat*}3
\sS{_U}r_0 &: \kk \to (s\cb)^{-1}(U,Uh\phi), &\qquad
\sS{_U}p_0 &: \kk \to (s\cb)^{-1}(Uh\phi,U), \\
\sS{_U}w_0 &: \kk \to (s\cb)^{-2}(U,Uh\phi), &\qquad
\sS{_U}v_0 &: \kk \to (s\cb)^{-2}(Uh\phi,U)
\end{alignat*}
are given such that
\begin{gather}
\sS{_U}r_0b_1=0, \qquad \sS{_U}p_0b_1=0, \notag \\
(\sS{_U}r_0\tens\sS{_U}p_0)b_2 - \sS{_U}\uni^\cb_0 = \sS{_U}w_0b_1,
\label{eq-rp-i-pr-UHphi-i} \\
(\sS{_U}p_0\tens\sS{_U}r_0)b_2 - \sS{_{Uh\phi}}\uni^\cb_0
= \sS{_U}v_0b_1. \notag
\end{gather}
Then there is an \ainf-functor $\psi:\cb\to\cc$ such that
$\Ob\psi=h$, there are natural \ainf-transformations
$r:\id_\cb\to\psi\phi$, $p:\psi\phi\to\id_\cb$ such that their
0\n-th components are the given $\sS{_U}r_0$, $\sS{_U}p_0$.
Moreover, $r$ and $p$ are inverse to each other in the sense that
\[ (r\tens p)B_2 \equiv \uni^\cb,
\qquad (p\tens r)B_2 \equiv \psi\phi\uni^\cb. \]
There exist unique up to equivalence natural \ainf-transformations
$t:\id_\cc\to\phi\psi$, $q:\phi\psi\to\id_\cc$ such that
$t\phi\equiv\phi r:\phi\to\phi\psi\phi$ and
$q\phi\equiv\phi p:\phi\psi\phi\to\phi$.

Finally, $\cc$ is unital with the unit
\begin{equation*}
\uni^\cc = (t\tens q)B_2: \id_\cc \to \id_\cc: \cc \to \cc,
\end{equation*}
$\phi$ and $\psi$ are unital \ainf-equivalences, quasi-inverse to each
other via mutually inverse isomorphisms $r$ and $p$, $t$ and $q$ (in
particular, $(q\tens t)B_2\equiv\phi\psi\uni^\cc$).
\end{theorem}

\begin{proof}
We have to satisfy the equations
\[ \psi b = b\psi, \qquad rb + br = 0.
\]
We already know the map $\Ob\psi$ and the component $r_0$. Let us
construct the remaining components of $\psi$ and $r$ by induction.
Given a positive integer $n$, assume that we have already found
components $\psi_m$, $r_m$ of the sought $\psi$, $r$ for $m<n$, such
that the equations
\begin{gather}
(\psi b)_m + (b\psi)_m = 0:
s\cb(X_0,X_1)\tens\dots\tens s\cb(X_{m-1},X_m) \to s\cc(X_0h,X_mh),
\label{eq-psibm-bpsim-0} \\
(rb + br)_m = 0: s\cb(X_0,X_1)\tens\dots\tens s\cb(X_{m-1},X_m)
\to s\cb(X_0,X_mh\phi) \label{eq-rb-brm-0}
\end{gather}
are satisfied for all $m<n$. Introduce a cocategory homomorphism
$\tilde{\psi}:Ts\cb\to Ts\cc$ of degree $0$ by its components
$(\psi_1,\dots,\psi_{n-1},0,0,\dots)$ and a
$(\id_\cb,\tilde{\psi}\phi)$-coderivation $\tilde{r}:Ts\cb\to Ts\cb$ of
degree $-1$ by its components $(r_0,r_1,\dots,r_{n-1},0,0,\dots)$.
Define a $(\tilde{\psi},\tilde{\psi})$-coderivation
$\lambda=\tilde{\psi}b-b\tilde{\psi}$ of degree $1$ and a map
$\nu=-\tilde{r}b-b\tilde{r}+(\tilde{r}\tens\lambda\phi)\theta:
Ts\cb\to Ts\cb$ of degree $0$. The commutator $\tilde{r}b+b\tilde{r}$
has the following property:
\[ (\tilde{r}b+b\tilde{r})\Delta
= \Delta\bigl[ 1\tens(\tilde{r}b+b\tilde{r})
+ (\tilde{r}b+b\tilde{r})\tens\tilde{\psi}\phi
+ \tilde{r}\tens\lambda\phi \bigr].
\]
By \propref{pro-theta-D-D-theta-theta} the map
$(\tilde{r}\tens\lambda\phi)\theta$ has a similar property
\[ (\tilde{r}\tens\lambda\phi)\theta\Delta
= \Delta\bigl[ 1\tens(\tilde{r}\tens\lambda\phi)\theta
+ (\tilde{r}\tens\lambda\phi)\theta\tens\tilde{\psi}\phi
+ \tilde{r}\tens\lambda\phi \bigr].
\]
Taking the difference we find that $\nu$ is an
$(\id_\cb,\tilde{\psi}\phi)$-coderivation. Equations
\eqref{eq-psibm-bpsim-0}, \eqref{eq-rb-brm-0} imply that $\lambda_m=0$,
$\nu_m=0$ for $m<n$ (the image of $(\tilde{r}\tens\lambda\phi)\theta$
is contained in $T^{\ge2}s\cb$).

The identity $\lambda b+b\lambda=0$ implies that
\begin{equation}
\lambda_nd = \lambda_nb_1 + \sum_{\alpha+1+\beta=n}
(1^{\tens\alpha}\tens b_1\tens1^{\tens\beta})\lambda_n = 0.
\label{eq-lamd-lb-bl}
\end{equation}
The identity
\begin{equation*}
\nu b-b\nu = (\tilde{r}\tens\lambda\phi)\theta b
- b(\tilde{r}\tens\lambda\phi)\theta
\end{equation*}
implies that
\begin{equation}
\nu_nb_1 - \sum_{\alpha+1+\beta=n}
(1^{\tens\alpha}\tens b_1\tens1^{\tens\beta})\nu_n
= - (r_0\tens\lambda_n\phi_1)b_2 = - \lambda_n\phi_1(r_0\tens1)b_2.
\label{eq-nub-bnu-rlfb}
\end{equation}
Set $N=s\cb(X_0,X_1)\tens_\kk\dots\tens_\kk s\cb(X_{n-1},X_n)$,
and introduce a chain map
\[ u = \Hom(N,\phi_1(r_0\tens1)b_2): \Hom^\bull(N,s\cc(X_0h,X_nh)) \to
\Hom^\bull(N,s\cb(X_0,X_nh\phi)).
\]
Since $\phi_1$ and $(r_0\tens1)b_2$ are homotopy invertible by
\lemref{lem-r-p-inverse}, the map $u$ is homotopy invertible as well.
Therefore, the complex $\Cone(u)$ is contractible by
\lemref{lem-contractible}. Equations \eqref{eq-lamd-lb-bl} and
\eqref{eq-nub-bnu-rlfb} in the form $-\lambda_nd=0$,
$\nu_nd+\lambda_nu=0$ imply that
\[ (\nu_n,\lambda_n) \in \Hom_\kk^0(N,s\cb(X_0,X_nh\phi)) \oplus
\Hom_\kk^1(N,s\cc(X_0h,X_nh)) = \Cone^0(u) \]
is a cycle. Hence, it is a boundary of some element
\[ (r_n,\psi_n) \in \Hom_\kk^{-1}(N,s\cb(X_0,X_nh\phi)) \oplus
\Hom_\kk^0(N,s\cc(X_0h,X_nh)) = \Cone^{-1}(u), \]
that is, $r_nd+\psi_n\phi_1(r_0\tens1)b_2=\nu_n$ and
$-\psi_nd=\lambda_n$. In other words, equations
\eqref{eq-psibm-bpsim-0}, \eqref{eq-rb-brm-0} are satisfied for $m=n$,
and we prove the existence of $\psi$ and $r$ by induction.

Since $r_0$ and $p_0$ are homotopy inverse to each other in the sense
of \eqref{eq-rp-i-pr-UHphi-i}, we find by \propref{pro-rfg-p-r-1} that
there exists a natural \ainf-transformation $p:\psi\phi\to\id_\cb$ such
that $r$ and $p$ are inverse to each other.

The existence of $t$, $q$ such that $t\phi\equiv\phi r$ and
$q\phi\equiv\phi p$ follows by \lemref{lem-cancel-phi}. Let us prove
that $\uni^\cc=(t\tens q)B_2$ is a unit of $\cc$. Due to
\lemref{lem-r-p-inverse} the maps $(r_0\tens1)b_2$, $(1\tens r_0)b_2$,
$(p_0\tens1)b_2$, $(1\tens p_0)b_2$ are homotopy invertible. Let $f$
denote a homotopy inverse map of
$\phi_1:s\cc(X,Y)\to s\cb(X\phi,Y\phi)$. The identity
$t\phi\equiv\phi r$ implies that
$\sS{_X}t_0\phi_1=\sS{_{X\phi}}r_0+\kappa b_1$. Hence,
\[ (\sS{_{X\phi}}r_0\tens1)b_2 \sim (\sS{_X}t_0\phi_1\tens1)b_2 \sim
f(\sS{_X}t_0\tens1)b_2\phi_1. \]
Therefore,
$(\sS{_X}t_0\tens1)b_2\sim\phi_1(\sS{_{X\phi}}r_0\tens1)b_2f$ is
homotopy invertible. Similarly,
\[ (1\tens\sS{_{Y\phi}}r_0)b_2 \sim (1\tens\sS{_Y}t_0\phi_1)b_2 \sim
f(1\tens\sS{_Y}t_0)b_2\phi_1 \]
implies that
$(1\tens\sS{_Y}t_0)b_2\sim\phi_1(1\tens\sS{_{Y\phi}}r_0)b_2f$ is
homotopy invertible. Similarly, $(\sS{_X}q_0\tens1)b_2$ and
$(1\tens\sS{_Y}q_0)b_2$ are homotopy invertible.

The computation made in \eqref{eq-r1b-p1b-pr1b1b} shows that the
product of the above homotopy invertible maps
\[ (q_0\tens1)b_2(t_0\tens1)b_2 \sim
- (t_0\tens q_0\tens1)(b_2\tens1)b_2 = - (\uni^\cc_0\tens1)b_2 \]
is the map we are studying. Similarly,
\[ (1\tens t_0)b_2(1\tens q_0)b_2 \sim
(1\tens t_0\tens q_0)(1\tens b_2)b_2 = (1\tens\uni^\cc_0)b_2. \]
We conclude that both $(\uni^\cc_0\tens1)b_2$ and
$(1\tens\uni^\cc_0)b_2$ are homotopy invertible.

Let us prove that $(\uni^\cc\tens\uni^\cc)B_2\equiv\uni^\cc$. Due to
\propref{pro-1uni-n2uni-KAinf} we have
\begin{equation}
\uni^\cc\phi = (t\tens q)B_2\phi \equiv (t\phi\tens q\phi)B_2 \equiv
(\phi r\tens\phi p)B_2 = \phi(r\tens p)B_2 \equiv \phi\uni^\cb.
\label{eq-phi-unital}
\end{equation}
Using \propref{pro-1uni-n2uni-KAinf} again we get
\[ (\uni^\cc\tens\uni^\cc)B_2\phi \equiv
(\uni^\cc\phi\tens\uni^\cc\phi)B_2 \equiv
(\phi\uni^\cb\tens\phi\uni^\cb)B_2 = \phi(\uni^\cb\tens\uni^\cb)B_2
\equiv \phi\uni^\cb \equiv \uni^\cc\phi. \]
By \lemref{lem-cancel-phi} we deduce that
$(\uni^\cc\tens\uni^\cc)B_2\equiv\uni^\cc$, therefore, $\uni^\cc$ is a
unit of $\cc$.

Let us prove that $t$ and $q$ are inverse to each other. By definition,
$(t\tens q)B_2=\uni^\cc$. Due to \propref{pro-1uni-n2uni-KAinf}
\[ (q\tens t)B_2\phi \equiv (q\phi\tens t\phi)B_2 \equiv
(\phi p\tens\phi r)B_2 = \phi(p\tens r)B_2 \equiv \phi\psi\phi\uni^\cb
\equiv \phi\psi\uni^\cc\phi. \]
By \lemref{lem-cancel-phi} $(q\tens t)B_2\equiv\phi\psi\uni^\cc$.
Hence, $t$ and $q$ are inverse to each other, as well as $r$ and $p$.
Therefore, $\phi$ and $\psi$ are equivalences, quasi-inverse to each
other.

Relation~\eqref{eq-phi-unital} shows that $\phi$ is unital. Let us
prove that $\psi$ is unital. We know that $\psi\phi$ is isomorphic to
the identity functor. Thus, $\psi\phi$ is unital by
\eqref{eq-riBg-rgiC}. Hence,
$\uni^\cb\psi\phi\equiv\psi\phi\uni^\cb\equiv\psi\uni^\cc\phi$. By
\lemref{lem-cancel-phi} we have $\uni^\cb\psi\equiv\psi\uni^\cc$, and
$\psi$ is unital. The theorem is proven.
\end{proof}

\begin{corollary}\label{cor-equiv-unital}
Let $\cc$, $\cb$ be unital \ainf-categories, and let $\phi:\cc\to\cb$
be an equivalence in $\sS{^u}A_\infty$. Then $\phi$ is unital.
\end{corollary}

\begin{proof}
Since $\phi$ is an equivalence, $\kf\phi$ is an equivalence as well.
Hence, $\phi_1$ is invertible in $\ck$. There exists an \ainf-functor
$\psi:\cb\to\cc$ quasi-inverse to $\phi$, and mutually inverse
isomorphisms $r:\id_\cb\to\psi\phi$, $p:\psi\phi\to\id_\cb$. In
particular, the assumptions of \thmref{thm-phi-C-B-hinv-equiv} are
satisfied by $\phi$, $\Ob\psi:\Ob\cb\to\Ob\cc$, $r_0$ and $p_0$. The
theorem implies that $\phi$ is unital.
\end{proof}

\begin{corollary}
Let $\cc$ be an \ainf-algebra and let $\cb$ be a unital \ainf-algebra
(viewed as \ainf-categories with one object). Let $\phi:\cc\to\cb$ be
an \ainf-homomorphism such that $\phi_1:(s\cc,b_1)\to(s\cb,b_1)$ is
homotopy invertible. Then $\cc$ and $\phi$ are unital, and $\phi$ is an
\ainf-equivalence.
\end{corollary}

Existence of $\phi$ with the above property might be taken as an
equivalence relation on the class of unital \ainf-algebras.

\subsection{\texorpdfstring{Strictly unital $A_\infty$-categories.}
{Strictly unital A8-categories.}}
 \label{sec-Strict-unit-cat}
A \emph{strict unit} of an object $X$ of an \ainf-category $\ca$ is an
element $1_X\in\ca^0(X,X)$, such that $(f\tens1_X)m_2=f$,
$(1_X\tens g)m_2=g$, whenever these make sense, and
$(\dots\tens1_X\tens\dots)m_n=0$ if $n\ne2$ (see e.g.
\cite{FukayaOhOhtaOno:Anomaly,Fukaya:FloerMirror-II,math.RA/9910179}).
We may write it as a map $1_X:\kk\to\ca(X,X)$, $1\mapsto1_X$. Assume
that $\ca$ has a strict unit for each object $X$. For example, a
differential graded category $\ca$ has strict units. Then we introduce
a coderivation $\uni^\ca:\id_\ca\to\id_\ca:\ca\to\ca$, whose components
are $\uni^\ca_0:\kk\to s\ca(X,X)$, $1\mapsto1_Xs=\sS{_X}\uni^\ca_0$,
and $\uni^\ca_k=0$ for $k>0$. The conditions on $1_X$ imply that
$(1\tens \uni^\ca_0)b_2=1:s\ca(Y,X)\to s\ca(Y,X)$ and
$(\uni^\ca_0\tens1)b_2=-1:s\ca(X,Z)\to s\ca(X,Z)$. One deduces that
$\uni^\ca$ is a natural \ainf-transformation. If an \ainf-category
$\ca$ has two such transformations -- strict units $\uni$ and $\uni'$,
then they must coincide because of the above equations. We call $\ca$
\emph{strictly unital} if it has a strict unit $\uni^\ca$. Naturally, a
strictly unital \ainf-category is unital.

For any \ainf-functor $f:\cc\to\ca$ the natural \ainf-transformation
$1_fs=f\uni^\ca:f\to f:\cc\to\ca$ has the components
$\sS{_X}(f\uni^\ca)_0=\sS{_{Xf}}\uni^\ca_0:\kk\to s\ca(Xf,Xf)$ and
$(f\uni^\ca)_k=0$ for $k>0$. It is the unit 2\n-endomorphism of $f$.

If \ainf-category $\cb$ is strictly unital, then so is
$\cc=A_\infty(\ca,\cb)$ for an arbitrary \ainf-category $\ca$. Indeed,
for an arbitrary \ainf-functor $f:\ca\to\cb$ there is a unit
2\n-endomorphism $1_fs=f\uni^\cb:f\to f$. We set
$\uni^\cc_0:\kk\to[sA_\infty(\ca,\cb)]^{-1}(f,f)$, $1\mapsto1_fs$, and
$\uni^\cc_k=0$ for $k>0$. For any element $r\in\cc(g,f)$ we have
$(r\tens1_fs)B_2=r$. For any element $p\in\cc(f,h)$ we have
$p(1_fs\tens1)B_2=p((f\uni^\cb)_0\tens1)b_2=-p$. We have also
$\uni^\cc B_1=0$ and $(\dots\tens\uni^\cc\tens\dots)B_n=0$ if $n>2$,
due to \eqref{eq-Bn-components}. Therefore, $\uni^\cc$ satisfies the
required conditions.

Another approach to $\uni^\cc$ uses the \ainf-functor
$M: TsA_\infty(\ca,\cb)\tens TsA_\infty(\cb,\cb)
\to TsA_\infty(\ca,\cb)=\cc$.
We have $(1\tens\id_\cb)M=\id_\cc$ by \eqref{eq-Mn0-k-component}, and
the natural \ainf-transformations $(1\tens\uni^\cb)M$ and $\uni^\cc$ of
$\id_\cc$ coincide. Indeed,
$[(1\tens\uni^\cb)M]_0:\kk\to(s\cc)^{-1}(f,f)$,
$1\mapsto(f\mid\uni^\cb)M_{01}=f\uni^\cb=\sS{_f}\uni^\cc_0$. For all
$n\ge0$ we have
$[(1\tens\uni^\cb)M]_n:r^1\tens\dots\tens r^n\mapsto(r^1\tens\dots\tens
r^n\tens\uni^\cb)M_{n1}$.
By \eqref{eq-Mn1-k-component} the components
\[ [(r^1\tens\dots\tens r^n\tens\uni^\cb)M_{n1}]_k =
\sum_l (r^1\tens\dots\tens r^n)\theta_{kl}\uni^\cb_l=
(r^1\tens\dots\tens r^n)\theta_{k0}\uni^\cb_0 \]
vanish for $n>0$.

\subsection{\texorpdfstring{Other examples of unital $A_\infty$-categories.}
{Other examples of unital A8-categories.}}
 \label{sec-Other-examples-unit-cat}
More examples of unital categories might be obtained via
\thmref{thm-phi-C-B-hinv-equiv}. An \ainf-category with a homotopy unit
in the sense of Fukaya, Oh, Ohta and Ono
\cite[Definition~20.1]{FukayaOhOhtaOno:Anomaly} clarified by Fukaya
\cite[Definition~5.11]{Fukaya:FloerMirror-II} is also a unital category
in our sense. Indeed, these authors enlarge given \ainf-category $\cc$
to a strictly unital \ainf-category $\cb$ by adding extra elements to
$\cc(X,X)$, so that the natural embedding
$(\cc(X,Y),m_1)\hookrightarrow(\cb(X,Y),m_1)$ were a homotopy
equivalence. Setting $r_0=p_0=\uni^\cb_0$ we view the above situation
as a particular case of \thmref{thm-phi-C-B-hinv-equiv}.

If an \ainf-functor $\phi:\cc\to\cb$ to a unital \ainf-category $\cb$
is invertible, then $\cc$ is unital. Indeed, since there exists an
\ainf-functor $\psi:\cb\to\cc$ such that $\phi\psi=\id_\cc$ and
$\psi\phi=\id_\cb$, then the map $\phi_1$ is invertible with inverse
$\psi_1$. The remaining data are $\Ob\psi:\Ob\cb\to\Ob\cc$ and
$\sS{_X}r_0=\sS{_X}p_0=\sS{_X}\uni^\cb_0:\kk\to s\cb(X,X)$. Since
$(\uni^\cb\tens\uni^\cb)B_2\equiv\uni^\cb$ we have
$(\sS{_X}\uni^\cb_0\tens\sS{_X}\uni^\cb_0)b_2-\sS{_X}\uni^\cb_0\in\im
b_1$
and conditions~\eqref{eq-rp-i-pr-UHphi-i} are satisfied. The data
constructed in \thmref{thm-phi-C-B-hinv-equiv} will be precisely
$\psi:\cb\to\cc$ and $r=p=\uni^\cb$. Since $\phi$ is unital by
\thmref{thm-phi-C-B-hinv-equiv}, we may choose
$\uni^\cc=\phi\uni^\cb\psi$ as a unit of $\cc$.

If a unital \ainf-category $\cc$ is equivalent to a strictly unital
\ainf-category $\cb$ via an \ainf-functor $\phi:\cc\to\cb$, then
$(1\tens\phi)M:A_\infty(\ca,\cc)\to A_\infty(\ca,\cc)$ is also an
equivalence for an arbitrary \ainf-category $\ca$ as
\propref{pro-AA-Au-Au-2-functor} shows. Thus, a unital \ainf-category
$A_\infty(\ca,\cc)$ is equivalent to a strictly unital \ainf-category
$A_\infty(\ca,\cb)$. In particular, if $\phi$ is invertible, then
$(1\tens\phi)M$ is invertible as well.

\subsection{\texorpdfstring{Cohomology of $A_\infty$-categories.}
 {Cohomology of A8-categories.}}
 \label{sec-H0-cohom-cat}
Using a lax monoidal functor from $\ck$ to some monoidal category we
get another 2\n-functor, which can be composed with $\kf$. For
instance, there is a cohomology functor
$H^\bull:\ck\to\ZZ\grad\text-\kk\modul$, which induces a 2\n-functor
$H^\bull:\KCat\to\ZZ\grad\text-\kk\text-\Cat$. In practice we will use
the 0\n-th cohomology functor $H^0:\ck\to\kk\modul$, the corresponding
2\n-functor $H^0:\KCat\to\kk\text-\Cat$, and the composite 2\n-functor
\[ A_\infty^u \rTTo^\kf \KCat \rTTo^{H^0} \kk\text-\Cat, \]
which is also denoted by $H^0$. It takes a unital \ainf-category $\cc$
into a $\kk$\n-linear category $H^0(\cc)$ with the same class of
objects $\Ob H^0(\cc)=\Ob\cc$. Its morphism space between objects $X$
and $Y$ is $H^0(\cc)(X,Y)=H^0(\cc(X,Y),m_1)$, the 0\n-th cohomology
with respect to the differential $m_1=sb_1s^{-1}$. The composition in
$H^0(\cc)$ is induced by $m_2$, and the units by $\uni^\cc_0s^{-1}$.

For example, the homotopy category $\Kht(\ca)$ of complexes of objects
of an abelian category $\ca$ is the 0\n-th cohomology $H^0(\Com(\ca))$
of the differential graded category of complexes $\Com(\ca)$.

\appendix
\section{Enriched 2-categories}\label{ap-sec-enrich}
Recall that $\fu$ is a fixed universe. Let $\cv=(\cv,\tens,c,\1)$ be a
symmetric monoidal $\fu$\n-category (that is, all $\cv(X,Y)$ are
$\fu$\n-small sets). In this article we shall use
$(\cv,\tens,c,\1)=(\kk\modul,\tens_\kk,\sigma,\kk)$, where $\sigma$ is
the permutation isomorphism, or
$(\cv,\tens,c,\1)=(\ck,\tens_\kk,c,\kk)$, where $\ck$ is the category
of differential graded $\kk$\n-modules, whose morphisms are chain maps
modulo homotopy, and $c$ is its standard symmetry. There is a notion of
a category $\cc$ enriched in $\cv$ ($\cv$\n-categories,
$\cv$\n-functors, $\cv$\n-natural transformations), see Kelly
\cite{KellyGM:bascec}, summarized e.g. in \cite{KerLyu}: for all objects $X$,
$Y$ of $\cc$ $\cc(X,Y)$ is an object of $\cv$. Denote by
$\cv\text-{\mathcal C}at$ the category, whose objects are
$\cv$\n-categories and morphisms are $\cv$\n-functors. Since $\cv$ is
symmetric, the category $\cv\text-{\mathcal C}at$ is symmetric monoidal
with the tensor product $\ca\tens\cb$ of $\cv$\n-categories $\ca$,
$\cb$ defined via $\Ob(\ca\tens\cb)=\Ob\ca\times\Ob\cb$,
$\ca\tens\cb(X\times Y,U\times V)=\ca(X,U)\tens\cb(Y,V)$. Thus, we may
consider the 1\n-category $\cv\text-{\mathcal C}at\text-{\mathcal C}at$
of $\cv\text-{\mathcal C}at$-categories and
$\cv\text-{\mathcal C}at$-functors. We may interpret it in the same
way, as ${\mathcal C}at\text-{\mathcal C}at$ is interpreted as the
category of 2\n-categories. So we say that objects of
$\cv\text-{\mathcal C}at\text-{\mathcal C}at$ are $\cv$-2-categories
$\fA$, as defined below. To restore the definition of a usual
2\n-category, it suffices to take
$\cv=(\fu\text{-Sets},\times,\emptyset)$.

\begin{definition}[$\cv$-2-category]\label{def-2-cat-pack}
A 1\n-unital 2\n-unital $\cv$-2-category $\fA$ consists of
\begin{events}
\item[\sbull] a class of objects $\Ob\fA$;

\item[\sbull] for any pair of objects $\ca,\cb\in\Ob\fA$ a
$\cv$\n-category $\fA(\ca,\cb)$;

\item[\sbull] for any object $\ca\in\Ob\fA$ a $\cv$\n-functor
$\underline\1\to\fA(\ca,\ca)$, $1\mapsto\id_\ca$;

\item[\sbull] for any triple $\ca$, $\cb$, $\cc$ of objects of $\fA$ a
$\cv$\n-functor
\[ \fA(\ca,\cb)\tens\fA(\cb,\cc) \to \fA(\ca,\cc),\qquad (f,g)\mapsto
fg, \]
\end{events}
such that the following $\cv$\n-functors are equal (modulo the
associativity isomorphism in $\cv$): $f\id=f=\id f$, $f(gh)=(fg)h$.
\end{definition}

Here the unit $\cv$\n-category $\underline\1$ has the set of objects
$\Ob\underline\1=\{1\}$, and $\underline\1(1,1)=\1$ is the unit object
of $\cv$. The above definition has an equivalent unpacked form, namely,
\defref{def-unit-2cat}. We also need generalizations of the above
$\cv$-2-categories -- 1\n-unital non-2-unital $\cv$-2-categories, which
contain unit 1\n-morphisms, but do not contain unit 2\n-morphisms. An
expanded definition of the latter follows. It seems that it does not
have a concise version.

\begin{definition}[Non-2-unital $\cv$-2-category]
\label{def-non-unit-2cat}
A 1\n-unital non-2-unital $\cv$-2-category $\fA$ consists of
\begin{events}
\item[\sbull] a class of objects $\Ob\fA$;

\item[\sbull] a class of 1\n-morphisms $\fA(\ca,\cb)$ for any pair
$\ca$, $\cb$ of objects of $\fA$;

\item[\sbull] an object of 2\n-morphisms $\fA(\ca,\cb)(f,g)\in\Ob\cv$
for any pair of 1\n-morphisms $f,g\in\fA(\ca,\cb)$;

\item[\sbull] a strictly associative composition of 1\n-morphisms
$\fA(\ca,\cb)\times\fA(\cb,\cc)\to\fA(\ca,\cc)$, $(f,g)\mapsto fg$;

\item[\sbull] a strict two-sided unit 1\n-morphism
$\id_\ca\in\fA(\ca,\ca)$ for each object $\ca$ of $\fA$;

\item[\sbull] a right action of a 1\n-morphism $k:\cb\to\cc$ on
2\n-morphisms
$\cdot k:\fA(\ca,\cb)(f,g)\to\fA(\ca,\cc)(fk,gk)\in\Mor\cv$;

\item[\sbull] a left action of a 1\n-morphism $e:\cd\to\ca$ on
2\n-morphisms
$e\cdot:\fA(\ca,\cb)(f,g)\to\fA(\cd,\cb)(ef,eg)\in\Mor\cv$;

\item[\sbull] a vertical composition of 2\n-morphisms
$m_2:\fA(\ca,\cb)(f,g)\tens\fA(\ca,\cb)(g,h)\to\fA(\ca,\cb)(f,h)
\in\Mor\cv$,
\end{events}
such that
\begin{events}
\item[\sbull] $m_2$ is associative (in monoidal category $\cv$);

\item[\sbull] the right and the left actions
\begin{events}
\item[(a)] commute with each other:
\[ (e\cdot)(\cdot k) = (\cdot k)(e\cdot),\qquad \text{for }
\cd \rTTo^e \ca \pile{\rTTo^f\\ \rTTo_g} \cb \rTTo^k \cc,
\]
\item[(b)] are associative:
\begin{alignat*}2
(\cdot k)(\cdot k') &= \cdot(kk'), &&\qquad \text{for }
\ca \pile{\rTTo^f\\ \rTTo_g} \cb \rTTo^k \cc \rTTo^{k'} \cd, \\
(e'\cdot)(e\cdot) &= (e'e)\cdot, &&\qquad \text{for }
\cc \rTTo^{e'} \cd \rTTo^e \ca \pile{\rTTo^f\\ \rTTo_g} \cb,
\end{alignat*}
\item[(c)] and unital: $(\cdot\id_\cb)=\id$, $(\id_\ca\cdot)=\id$;
\end{events}

\item[\sbull] the right and the left actions of 1\n-morphisms on
2\n-morphisms preserve the vertical composition:
\[
\begin{diagram}[inline,nobalance]
\fA(\ca,\cb)(f,g)\tens\fA(\ca,\cb)(g,h) & \rTTo^{m_2} &
\fA(\ca,\cb)(f,h) \\
\dTTo<{(\cdot k)\tens(\cdot k)} & = & \dTTo>{\cdot k} \\
\fA(\ca,\cc)(fk,gk)\tens\fA(\ca,\cc)(gk,hk) & \rTTo^{m_2} &
\fA(\ca,\cc)(fk,hk)
\end{diagram}
\qquad \text{for }
\ca \pile{\rTTo^f\\ \rTTo~g\\ \rTTo_h} \cb \rTTo^k \cc,
\]
\[
\begin{diagram}[inline,nobalance]
\fA(\ca,\cb)(f,g)\tens\fA(\ca,\cb)(g,h) & \rTTo^{m_2} &
\fA(\ca,\cb)(f,h) \\
\dTTo<{(e\cdot)\tens(e\cdot)} & = & \dTTo>{e\cdot} \\
\fA(\cd,\cb)(ef,eg)\tens\fA(\cd,\cb)(eg,eh) & \rTTo^{m_2} &
\fA(\cd,\cb)(ef,eh)
\end{diagram}
\qquad \text{for }
\cd \rTTo^e \ca \pile{\rTTo^f\\ \rTTo~g\\ \rTTo_h} \cb;
\]

\item[\sbull] the distributivity law holds:
\begin{diagram}
&& \fA(\cb,\cc)(h,k)\tens\fA(\ca,\cb)(f,g) \\
& \ruTTo^c & \dTTo>{(f\cdot)\tens(\cdot k)} \\
\fA(\ca,\cb)(f,g)\tens\fA(\cb,\cc)(h,k) &&
\fA(\ca,\cc)(fh,fk)\tens\fA(\ca,\cc)(fk,gk) \\
\dTTo<{(\cdot h)\tens(g\cdot)} & = & \dTTo>{m_2} \\
\fA(\ca,\cc)(fh,gh)\tens\fA(\ca,\cc)(gh,gk) & \rTTo^{m_2} &
\fA(\ca,\cc)(fh,gk)
\end{diagram}
\[ \text{for }
\ca \pile{\rTTo^f\\ \rTTo_g} \cb \pile{\rTTo^h\\ \rTTo_k} \cc.
\]
\end{events}
\end{definition}

The following definition is equivalent to \defref{def-2-cat-pack}.

\begin{definition}[1\n-unital left-2-unital $\cv$-2-category]\label{def-unit-2cat}
A 1\n-unital {\color{blue} left-2-unital} $\cv$-2-category $\fA$
consists of the same data as in \defref{def-non-unit-2cat} plus a
morphism $1_f:\1\to\fA(\ca,\cb)(f,f)$ for any 1\n-morphism $f$, which
is a two-sided unit with respect to $m_2$, such that homomorphisms
$e\cdot$ preserve the units $1_-$.
 {\color{blue} Moreover, if homomorphisms $\cdot k$ also preserve the
units $1_-$, such $\fA$ is the same as a 1\n-unital 2\n-unital
$\cv$-2-category.}
\end{definition}

\section{Contractibility}\label{ap-Contractibility}
One can avoid using the following lemma in this article. However, it
might be used in order to replace inductive constructions with
recurrent formulas.

\begin{lemma}\label{lem-contractible}
Let a chain map $u:A\to C$ be homotopically invertible. Then $\Cone(u)$
is contractible.
\end{lemma}

\begin{proof}
The homotopy category $\ck=\Kht(\kk\modul)$ is triangulated and it has
a distinguished triangle
$A \rTTo^u C \rTTo^p \Cone(u) \rTTo^q A[1] \rTTo^{u[1]}$
(e.g. \cite[Corollaire~5.13]{Grivel-derivees}). Since $up=0$, $qu[1]=0$
(e.g. \cite[Proposition~2.8]{Grivel-derivees}), and $u$ is invertible
in $\ck$, we deduce that $p=0$ and $q=0$ in $\ck$. Since
$\ck(\Cone(u),\_)$ is a homological functor (e.g.
\cite[Proposition~2.10]{Grivel-derivees}), we have
$\ck(\Cone(u),\Cone(u))=0$, that is, $\Cone(u)\simeq0$ in $\ck$.
\end{proof}

Let us construct an explicit homotopy between $\id_{\Cone(u)}$ and
$0_{\Cone(u)}$. There exists a chain map $v:C\to A$ homotopically
inverse to $u$. That is, there are maps $h':A\to A$, $h'':C\to C$ of
degree $-1$ such that $uv=1+h'd^A+d^Ah':A\to A$,
$vu=1+h''d^C+d^Ch'':C\to C$. Using the notation at the end of
\secref{Conventions} we define a map $h:\Cone(u)\to\Cone(u)$ of degree
$-1$ by the formula
$(c,a)h=(-ch'',cv+ah')$, $(c,a)\in C^k\oplus A^{k+1}=\Cone^k(u)$.
Let us compute the boundary of $h$:
\begin{align*}
(c,a)(hd+dh) &= (-ch'',cv+ah')d + (cd^C+au,-ad^A)h \\
&= (-ch''d^C+cvu+ah'u,-cvd^A-ah'd^A) \\
&\quad + (-cd^Ch''-auh'',cd^Cv+auv-ad^Ah') \\
&= (c+ah'u-auh'',a).
\end{align*}
Hence, $hd+dh=1-f:\Cone(u)\to\Cone(u)$, where the map
$f:\Cone(u)\to\Cone(u)$ is defined via $(c,a)f=(auh''-ah'u,0)$. We
conclude that $f$ is a chain map homotopic to the identity map. A
sequence of equivalences
$\id_{\Cone(u)}\sim f=\id_{\Cone(u)}f\sim f^2=0$ proves that $\Cone(u)$
is contractible. It gives also an explicit homotopy -- the map
$\underline{h}=h+hf:\Cone(u)\to\Cone(u)$ of degree $-1$, which
satisfies $\id_{\Cone(u)}=\underline{h}d+d\underline{h}$. This homotopy
might be used to replace inductive constructions in this paper with
recurrent formulas.

\begin{acknowledgement}
I am grateful to all the participants of the \ainf-category seminar
at the Institute of Mathematics, Kyiv, especially to Yu.~Bespalov,
O.~Manzyuk, S.~Ovsienko.
\end{acknowledgement}

\end{document}